\documentclass[reqno]{memo-l} 

\usepackage{amssymb}

\DeclareFontFamily{OT1}{rsfs}{}
\DeclareFontShape{OT1}{rsfs}{m}{n}{ <-7> rsfs5 <7-10> rsfs7 <10-> rsfs10}{} 
\DeclareMathAlphabet{\mathscr}{OT1}{rsfs}{m}{n}

\def\mysavedown#1{\edef\mysubs{\mysubs#1}}
\def\mysaveup#1{\edef\mysups{\mysups#1}}
\def\mydown#1{{\mytensor}_{\vphantom{\mysubs}#1}}
\def\myup#1{{\mytensor}^{\vphantom{\mysups}#1}}
\def\tensor#1#2{
  #1
  \def\mytensor{\vphantom{#1}}
  \def\mysubs{\relax}
  \def\mysups{\relax}
  \let\down=\mysavedown
  \let\up=\mysaveup
  #2
  \let\down=\mydown
  \let\up=\myup
  #2
  }

\newcommand{\letters}
  {\renewcommand{\theenumi}{\alph{enumi}}
   \renewcommand{\labelenumi}{(\theenumi)}}

\newcommand{\defname}[1]{{\it #1}}

\newcommand{\Ker}{\operatorname{Ker}}
\newcommand{\Image}{\operatorname{Im}}
\newcommand{\Real}{\operatorname{Re}}

\newcommand{\Tr}{\operatorname{Tr}}
\newcommand{\Id}{\operatorname{Id}}
\newcommand{\SO}{\operatorname{SO}}
\newcommand{\C}{\mathbb C}
\newcommand{\R}{\mathbb R}
\newcommand{\B}{\mathbb B}
\renewcommand{\S}{\mathbb S}
\renewcommand{\setminus}{\smallsetminus}
\renewcommand{\emptyset}{\varnothing}
\renewcommand{\to}{\rightarrow}

\renewcommand{\centerdot}{\mathbin{\text{\protect\raisebox{-.3ex}[1ex][0ex]{\Large{$\cdot$}}}}}
\newcommand{\cross}{\mathbin{\times}}
\newcommand{\intersect}{\mathbin{\cap}}
\newcommand{\tprod}{\mathbin{\otimes}}
\newcommand{\Del}{\nabla}
\newcommand{\DeltaL}{\Delta_{\mathrm L}}

\DeclareMathOperator{\Hom}{Hom}

\DeclareMathOperator{\Vol}{Vol}

\DeclareMathOperator{\Ortho}{O}
\newcommand{\half}{{\tfrac12}}
\newcommand{\mapsinto}{\mathrel{\hookrightarrow}}
\renewcommand{\phi}{\varphi}
\renewcommand{\epsilon}{\varepsilon}
\renewcommand{\hat}{\widehat}
\let\scr=\mathscr

\def\crn#1#2{{\vcenter{\vbox{
        \hbox{\kern#2pt \vrule width.#2pt height#1pt
           }
          \hrule height.#2pt}}}}
\newcommand{\intprod}{\mathchoice\crn54\crn54\crn{3.75}3\crn{2.5}2}
\newcommand{\into}{\mathbin{\intprod}}

\newcommand{\tw}{\widetilde}

\newcommand{\hyp}{{\breve g}}
\newcommand{\Phyp}{{\breve P}}

\newcommand{\<}{\langle}
\renewcommand{\>}{\rangle}
\newcommand{\del}{\partial}
\renewcommand{\H}{\mathbb H}

\newcommand{\grad}{\operatorname{grad}}
\newcommand{\supp}{\operatorname{supp}}
\newcommand{\dm}{{\partial M}}

\newcommand{\Mbar}{{\overline M}}
\newcommand{\gbar}{{\overline g}}
\renewcommand{\hbar}{{\overline h}}
\newcommand{\ghat}{{\widehat g}}
\newcommand{\hhat}{{\widehat h}}

\newcommand{\rade}{R{\aa}de}

\newcommand{\gcirc}{\overset{\circ}{g}}

\renewcommand{\)}{\textup{)}}

\theoremstyle{plain}
\newtheorem{theorem}{Theorem}[chapter]
\newtheorem{lemma}[theorem]{Lemma}
\newtheorem{proposition}[theorem]{Proposition}
\newtheorem{corollary}[theorem]{Corollary}
\newtheorem{bigtheorem}{Theorem}[chapter]

\newtheorem{bigproposition}[bigtheorem]{Proposition}

\theoremstyle{definition}

\theoremstyle{remark}
\newtheorem*{remark}{Remark}

\numberwithin{equation}{chapter}

\begin{document}

\frontmatter
\title[Fredholm Operators and Einstein Metrics]
{Fredholm Operators and Einstein Metrics
on Conformally Compact Manifolds}
\author {John M. Lee}
\address{
  University of Washington \\
  Department of Mathematics\\
  Box 354350\\
  Seattle, WA 98195-4350}
\email{lee@math.washington.edu}
\urladdr{http://www.math.washington.edu/\~{}lee}
\thanks{Research supported in part by National Science Foundation grants
DMS-8901493 and DMS-9404107.}

\date{February 16, 2004}



\begin{abstract}
The main purpose of this monograph is to give an
elementary and self-contained account of the 
existence of asymptotically hyperbolic 
Einstein metrics with prescribed
conformal infinities sufficiently close to that of a given
asymptotically hyperbolic Einstein metric with nonpositive curvature.
The proof is based on an elementary derivation
of sharp Fredholm theorems for self-adjoint
geometric linear elliptic operators on asymptotically
hyperbolic manifolds.  
\end{abstract}

\maketitle

% Insert an empty page after the title page, to get the
% odd pages to print on right-hand sides.

\newpage

\setcounter{page}{4}
\tableofcontents

\mainmatter

%%
%%
%% Fredholm operators and Einstein metrics
%% on conformally compact manifolds
%%
%% by John M. Lee
%% 
%% Chapter 1
%%

\chapter{Introduction}\label{intro-section}

In 1985, Charles Fefferman and Robin Graham \cite{FG} introduced a new
and powerful approach to the study of local invariants of conformal
structures, based on the observation that the group of conformal
automorphisms of the $n$-sphere (the M\"obius group) is essentially the
same as the group of Lorentz transformations of $(n+2)$-dimensional
Minkowski space.  Concretely, one sees this by noting that the conformal
structure of the $n$-sphere can be obtained by viewing the sphere as a
cross-section of the forward light cone in Minkowski space.  Fefferman
and Graham attempted to embed an arbitrary conformal $n$-manifold into
an $(n+2)$-dimensional Ricci-flat Lorentz manifold.  They showed that
such a Lorentz metric (which they called the \defname{ambient metric})
can be constructed by formal power series for any conformal Riemannian
metric, to infinite order when $n$ is odd and to order $n/2$ when $n$ is
even; and that this formal metric is a conformal invariant of the
original conformal structure.

Once this Lorentz metric (or at least its formal power series) is
constructed, its pseudo-Riemannian invariants then automatically give
conformal invariants of the original conformal structure.  Combined with
later work of Bailey, Eastwood, and Graham \cite{BEG}, this construction
produces all local scalar conformal invariants of odd-di\-men\-sion\-al
conformal Riemannian structures.

In \cite{Graham-Lee}, Graham and I adapted this idea to the global setting.  We
began with the observation that there is a third natural realization of
the M\"obius group, as the set of isometries of the hyperbolic metric on
the interior of the unit ball $\B^{n+1}\subset\R^{n+1}$.  As 
noted by Fefferman and Graham in \cite{FG}, given a compact Riemannian
manifold $(S,\ghat)$, the problem of finding a Ricci-flat ambient
Lorentz metric for the conformal structure $[\ghat]$ is equivalent to
that of finding an asymptotically hyperbolic Einstein metric $g$ on the
interior of an $(n+1)$-dimensional manifold-with-boundary $\overline M$
that has $S$ as boundary and $[\ghat]$ as \defname{conformal
infinity} in the following sense: For any smooth, positive
defining function
$\rho$ for $S=\dm$, the metric $\rho^2 g$ extends continuously to $\overline M$
and restricts to a metric in the conformal class $[\ghat]$ on $S$.  This
suggests the following natural problem: Given 
a compact manifold-with-boundary $\overline M$ and a
conformal structure
$[\ghat]$ on $\del M$,
can one find a complete, asymptotically hyperbolic Einstein metric $g$
on the interior manifold $M$ that has $[\ghat]$ as conformal infinity?
Such a metric is said to be \defname{conformally compact}.
To the extent that the answer is yes and the resulting
Einstein metric is unique up to boundary-fixing diffeomorphisms,
one would thereby obtain a correspondence between
global conformal invariants of $[\ghat]$ and global
Riemannian invariants of $g$.
Graham and I showed in \cite{Graham-Lee} 
that every conformal structure on $\S^n$
sufficiently close to that of the round metric is the conformal
infinity of an Einstein metric close to the hyperbolic metric.  

In recent years, interest in asymptotically hyperbolic Einstein
metrics has risen dramatically, in no small part 
because of the role they
play in physics.  In fact, the notion of conformal
infinity for a (pseudo-) Riemannian metric was originally 
introduced by Roger Penrose \cite{Penrose} in order to
analyze the behavior of gravitational energy in
asymptotically flat space-times.  
More recently, 
asymptotically hyperbolic Einstein metrics have
begun to play a central role in 
the ``AdS/CFT correspondence'' of quantum field theory.
This is not the place for a complete survey of
the relevant literature, but let me just refer the reader
to \cite{Anderson-closed,Anderson-L2,Anderson-adscft,Biquard-survey,Biquard-selfdual,Biquard,Cai-Galloway,Graham-Lee,Graham-Witten,Hitchin-SU2,Hitchin,LeBrun-H-space,Lee,Min-Oo,Pedersen,Petersen-Maldacena,Wang,Witten-holography}.

The principal purpose of this monograph is to give an
elementary and self-contained account of the 
following generalization
of the perturbation result of 
\cite{Graham-Lee}.
First, a couple of definitions: The 
\defname{Lichnerowicz Laplacian} (cf.~\cite{Besse})
is the operator $\DeltaL$ 
defined on symmetric
$2$-tensors by $\DeltaL=\Del^*\Del + 2 \mathring{Rc}
- 2 \mathring{Rm}$,
where $\mathring{Rc}$ and $\mathring{Rm}$
are the natural actions of the 
Ricci and Riemann curvature tensors on symmetric 
$2$-tensors given in coordinates by
\begin{equation}\label{eq:def-Rc-Rm}
\begin{aligned}
\mathring {Rc}(u)_{ij} &= \half (R_{ik}u_j{}^k
+ R_{jk} u_i{}^k),\\
\mathring {Rm}(u)_{ij} &= R_{ikjl} u^{kl}.
\end{aligned}
\end{equation}
For any conformal class of Riemannian metrics 
on $\del M$, the \defname{Yamabe invariant} is defined as 
the infimum of the total scalar curvature
$\int_{\del M} S_{\ghat}\, dV_{\ghat}$ over unit-volume 
metrics $\ghat$
in the conformal class.

\begin{bigtheorem}\label{thm:einstein}
Let $M$ be the interior of a smooth, 
compact, $(n+1)$-dimensional
manifold-with-boundary $\Mbar$, $n\ge 3$, and let $h$ be an
Einstein metric
on $M$ that 
is conformally compact of class $C^{l,\beta}$ with 
$2\le l\le n-1$ and
$0<\beta<1$.  Let $\rho$ be a smooth defining function
for $\del M$, and let $\hhat = \rho^2 h|_{\del M}$.
Suppose the operator $\DeltaL+2n$ associated with $h$
has trivial $L^2$ kernel on the space of 
trace-free symmetric $2$-tensors.
Then there is a constant $\epsilon>0$ such that for any
$C^{l,\beta}$ Riemannian
metric $\ghat$ on $\dm$ with $\|\ghat-\hhat\|_{C^{l,\beta}}<\epsilon$,
there exists an Einstein metric $g$ on $M$ that
has $[\ghat]$ as conformal infinity and
is conformally compact
of class $C^{l,\beta}$.
In particular, this is the case if 
either of the following hypotheses is satisfied:
\begin{enumerate}\letters
\item\label{part:neg-curv}
$h$ has nonpositive sectional
curvature.
\item\label{part:yamabe}
The Yamabe invariant of $[\hhat]$ is nonnegative and
$h$ has sectional curvatures bounded above by
$(n^2-8n)/(8n-8)$.
\end{enumerate}
\end{bigtheorem}

An important feature of this result 
is that the conformal compactification 
of the Einstein
metric has optimal 
H\"older regularity 
up to the boundary (in terms of standard H\"older spaces on $\overline M$), 
at least when $n$ is even.
The work of Fefferman and Graham showed that generically
for $n$ even there will be a $\rho^n\log\rho$ term
in the asymptotic expansion of $\overline g$ if $g=\rho^{-2}\overline g$ is
Einstein; thus we
cannot expect to find Einstein metrics with $C^n$ conformal
compactifications in general.  If $h$ has a sufficiently 
smooth conformal
compactification and $\hat g$ is close to  $\hat h$ with
sufficiently many derivatives, then 
Theorem \ref{thm:einstein} gives an 
Einstein metric that is conformally compact of class
$C^{n-1,\beta}$ for any $0<\beta<1$, which is the best H\"older regularity that
can be expected when $n$ is even.  
The results
of \cite{CDLS} show that any Einstein metric with 
smooth conformal infinity and $C^2$ conformal compactification
actually has an infinite-order
asymptotic expansion in powers of $\rho$ and $\log\rho$,
and in fact is smooth when $n$ is odd.

A version of Theorem \ref{thm:einstein}
(but with
less boundary regularity of the resulting Einstein metrics)
was also proved independently 
by Olivier Biquard \cite{Biquard}
around the same time as this monograph was originally completed.  
In addition, 
analogous perturbation results starting with
K\"ahler-Einstein metrics on bounded domains
have been proved independently by Biquard \cite{Biquard} and John Roth
\cite{Roth-thesis}, and for
Einstein metrics asymptotic to
the quaternionic and octonionic
hyperbolic metrics by Biquard \cite{Biquard}.
See
below for more on this.
See also 
\cite{Anderson-prescribed,Anderson-closed} for 
some recent results by Michael Anderson on existence and uniqueness of 
asymptotically hyperbolic Einstein metrics on $4$-manifolds, and 
\cite{Delay,Delay-Herzlich} for results
by Erwann Delay and Marc Herzlich
on the related problem of prescribing the Ricci curvatures of metrics
close to the Einstein models.

The basic approach in this monograph is similar to that of
\cite{Graham-Lee}.  Because the Einstein equation 
is invariant under the full diffeomorphism group of $M$, it
is not elliptic as it stands.
We obtain
an elliptic equation 
by adding a gauge-breaking term:  Fixing a
conformally compact ``reference metric'' 
$g_0$, let $\mathbf{\Delta}_{gg_0}(\Id)$ denote the
harmonic map Laplacian of the identity map, considered as a map from
$(M,g)$ to $(M,g_0)$, and let
$\delta_g$ be the divergence operator associated with
$g$.  Then, as is by now familiar, the nonlinear equation
\begin{equation}\label{regularized-einstein}
Q(g,g_0) := Rc_{g} + n g - \delta_g^*(\mathbf{\Delta}_{gg_0}(\Id)) = 0
\end{equation}
is elliptic as a function of $g$.  
Under relatively mild assumptions on $g$ (see \cite[Lemma
2.2]{Graham-Lee}), the solutions to (\ref{regularized-einstein}) are exactly
those Einstein metrics $g$ such that $\mathbf{\Delta}_{gg_0}(\Id)=0$, i.e.,
such that $\Id\colon (M,g)\to (M,g_0)$ is harmonic.

The linearization of the left-hand side of \eqref{regularized-einstein}
with respect to $g$ at a conformally compact Einstein metric $h$ is
\begin{equation}\label{linearized-Einstein}
D_1Q_{(h,h)} = \tfrac12(\DeltaL + 2n ),
\end{equation}
where $\DeltaL$ is the
Lichnerowicz Laplacian defined above.
In \cite{Graham-Lee}, we proved 
that when $h$ is the hyperbolic metric on the unit ball $\B^{n+1}$, 
$\DeltaL+2n$ is an isomorphism between certain
weighted H\"older spaces.  
However, our methods in that paper
were insufficient to prove a sharp isomorphism
theorem, with the consequence that our boundary
regularity results were not optimal,
so even in the case of hyperbolic space
Theorem \ref{thm:einstein} is an improvement on the results of
\cite{Graham-Lee}. 

Much of this monograph is devoted to an elementary proof of some 
general
sharp Fredholm and isomorphism theorems for self-adjoint geometric
linear elliptic operators on asymptotically hyperbolic
manifolds.  
These operators are, in particular,
\defname{uniformly degenerate}, in the
terminology introduced in \cite{Graham-Lee} (with perhaps less smoothness of the
coefficients than we insisted on there):  A partial
differential operator $P$
on a manifold with
boundary is said to be uniformly degenerate if
in local coordinates 
$(\theta^1,\dots,\theta^n,\rho)$ such that 
$\rho=0$ defines the boundary,
$P$ can be written locally as a a system of partial
differential operators that are polynomials in
the vector
fields $(\rho\partial/\partial \theta^1,\dots,
\rho\partial/\partial \theta^n,\rho\partial/\partial \rho)$ with coefficients
that are at least continuous up to the boundary.  

To get an idea
of what can be expected, consider a simple example: the scalar Laplacian
$\Delta =d ^*d$ on functions.  When applied to a function of
$\rho$ alone, this becomes an ordinary differential operator with a regular
singular point at $\rho=0$.  From classical ODE theory, we know that a
second-order ordinary differential operator $L$ with a regular singular
point at $0$ has two \defname{characteristic exponents} $s_1$ and $s_2$,
defined by $L(\rho^{s_i}) = O(\rho^{s_i+1})$; if $s_1\ne s_2$, the homogeneous
equation $Lu=0$ has two independent solutions $u_1,u_2$ with $u_i =
O(\rho^{s_i})$, and the inhomogeneous equation $Lu=f=O(\rho^s)$ can be solved
with $u=O(\rho^s)$ whenever $s$ is not a characteristic exponent.

To see how this generalizes to systems of
partial differential operators, let $P\colon
C^\infty(M;E)\to C^\infty(M;F)$
be a uniformly degenerate operator 
of order $m$ acting between 
(real or complex) tensor bundles $E$ and $F$,
and let $s$ be any complex number.  
Define the \defname{indicial
map} $I_s(P)\colon E|_{\dm}\to F|_{\dm}$ by setting
\begin{displaymath}
I_s(P)\overline u := \rho^{-s} P(\rho^s \overline u)|_{\dm},
\end{displaymath}
where $\overline u$ is any $C^m$ section of $E|_{\dm}$, extended arbitrarily
to a $C^m$ section of $E$ near $\dm$.  
It follows easily from the definition
of uniformly degenerate operators that $I_s(P)$ is a (pointwise) bundle
endomorphism whose coefficients in any local coordinates are polynomials
in $s$ with continuous coefficients depending on $\theta\in \dm$.  (In the
special case of a single ODE with a regular singular point, $I_s(P)$ is just
multiplication by a number depending polynomially on $s$, called the
{\it indicial polynomial} of the ODE.)  We say a number $s$ is a
\defname{characteristic exponent} for $P$ if $I_s(P)$ is singular somewhere
on $\dm$.  
Just as in the ODE case, one can construct formal series
solutions to $Pu=f$ in which the characteristic exponents correspond to
powers of $\rho$ whose coefficients are arbitrary.

Now consider a 
formally self-adjoint operator acting on sections of a 
tensor 
bundle $E$ of
type $\binom{r_1}{r_2}$ (i.e., covariant rank $r_1$ and contravariant
rank $r_2$).  In this case, 
the difference $r=r_1-r_2$, which we call the \defname{weight}
of $E$, is of central importance, and 
the exponent $n/2-r$ plays a special role.  
If $\overline u$ is a section of $E$ 
that is continuous up to the boundary, then 
$\rho^{n/2-r}\overline u$ 
is just on the borderline of being in $L^2$ (see Lemma
\ref{lemma:being-in-Lp}\eqref{part:overlineu-in-Lp}).  
If $P$ is formally self-adjoint,
the set of characteristic
exponents turns out to be 
symmetric about the line $\Real s = n/2-r$ 
(Corollary \ref{cor:symmetric-char-exp}).
Therefore, we define the
{\it indicial radius} of $P$ to be the smallest real
number $R\ge 0$ such that $P$ has a characteristic
exponent whose real part is equal to $n/2-r+R$.
A little experimentation leads one to expect that $Pu=f$
should be well-posed roughly 
when $u$ and $f$ behave like $\rho^s \overline
u$, where $\overline u$ is continuous up to the boundary and 
$n/2-r-R<s<n/2-r+R$, because for lower values of $s$
one expects $P$ to have an infinite-dimensional kernel, and
for higher values one expects an infinite-dimensional cokernel.
(This will be made precise in Chapter \ref{section:fredholm}.)
Therefore, a fundamental necessary condition for all
of our Fredholm results will be that $P$ has positive
indicial radius.

To  control  the order of vanishing or singularity of $u$ and $f$
at the boundary, we work in weighted Sobolev and H\"older spaces, with
weights given by powers of the defining function $\rho$.  Precise
definitions are given in Chapter \ref{spaces-section}, but the spaces we work
in can be roughly defined as follows: The weighted Sobolev space
$H^{k,p}_\delta$ is just the space of tensor fields of the form $\rho^\delta u$
for $u$ in the usual intrinsic Sobolev space $H^{k,p}$ (tensor fields
with $k$ covariant derivatives in $L^p$ with respect to $g$); and the
weighted H\"older space $C^{k,\alpha}_\delta$ consists of tensor fields of
the form $\rho^\delta u$ for $u$ in the usual H\"older space $C^{k,\alpha}$.
In each case, the Sobolev and H\"older norms are defined with respect to
a fixed conformally compact metric on $M$.  

An elementary criterion for operators to be Fredholm 
on $L^2$ is given in the following proposition.  Terms
used in the statements of these results 
will be defined in Chapters \ref{spaces-section}
and \ref{elliptic-operator-section}.

\begin{bigproposition}\label{prop:L2-Fredholm}
Let $(M,g)$ be a connected
asymptotically hyperbolic $(n+1)$-manifold of 
class $C^{l,\beta}$, with $n\ge 1$, $l\ge 2$, 
and $0\le\beta< 1$, and 
let $E\to M$ be a
geometric tensor bundle over $M$.  Suppose 
$P\colon C^\infty(M;E)\to C^\infty(M;E)$ is an 
elliptic, formally self-adjoint, geometric partial differential 
operator of order $m$, $0<m\le l$.  
As an unbounded operator, 
$P\colon L^2(M;E)\to L^2(M;E)$ is Fredholm if and only if 
there exist a compact set $K\subset M$ and a 
positive constant $C$ such that 
\begin{equation}\label{eq:asymptotic-L2-estimate}
\|u\|_{L^2}\le C\|Pu\|_{L^2} \text{ for all $u\in C^\infty_c(M\setminus K;E)$.}
\end{equation}
\end{bigproposition}

The main analytic result we need for Theorem \ref{thm:einstein}
is the following sharp Fredholm theorem in
weighted Sobolev and H\"older spaces for geometric elliptic
operators.  (See below for references to 
other proofs of these and similar results using different approaches.) 

\begin{bigtheorem}\label{thm:main-fredholm}
Let $(M,g)$ be a connected
asymptotically hyperbolic $(n+1)$-manifold of 
class $C^{l,\beta}$, with $n\ge 1$, $l\ge 2$, 
and $0\le\beta< 1$, and 
let $E\to M$ be a
geometric tensor bundle over $M$.  Suppose 
$P\colon C^\infty(M;E)\to C^\infty(M;E)$ is an 
elliptic, formally self-adjoint, geometric partial differential 
operator of order $m$, $0<m\le l$, and   
assume $P$ satisfies 
\eqref{eq:asymptotic-L2-estimate}.
\begin{enumerate}\letters
\item\label{part:pos-indicial}
The indicial radius $R$ of $P$ is positive.
\item\label{sobolev-fredholm-result}
If $1<p<\infty$ and $m\le k\le l$, the natural extension
\begin{displaymath}
P\colon H^{k,p}_\delta(M;E)\to H^{k-m,p}_\delta(M;E)
\end{displaymath}
is Fredholm if and only if $|\delta +n/p-n/2|<R$.
In that case, its index is zero, and its kernel is equal to
the $L^2$ kernel of $P$.
\item\label{holder-fredholm-result}
If $0<\alpha<1$ and 
$m< k+\alpha\le l+\beta$, the natural extension
\begin{displaymath}
P\colon C^{k,\alpha}_\delta(M;E)\to C^{k-m,\alpha}_\delta(M;E)
\end{displaymath}
is Fredholm if and only if $|\delta -n/2|<R$.
In that case, its index is zero, and its kernel is equal to
the $L^2$ kernel of $P$.
\end{enumerate}
\end{bigtheorem}

The statement and proof of this theorem extend easily
to operators on spinor bundles when $M$ is a spin manifold.  
We restrict attention here
to the case of tensor bundles
mainly for simplicity of exposition.

We define a \defname{Laplace operator} to be a 
formally self-adjoint second-order 
geometric operator of the form $\Del^*\Del + \scr K$,
where $\scr K$ is a bundle endomorphism.
Of course, the ordinary Laplacian on functions and the
covariant Laplacian on tensor fields are obvious examples
of Laplace operators, as are the
Lichnerowicz Laplacian $\DeltaL$ defined above, and the 
Laplace-Beltrami operator on differential forms by virtue of
Bochner's formula.

The most important example to which we will apply these
results is the Lichnerowicz Laplacian on 
symmetric $2$-tensors.

\begin{bigproposition}\label{prop:lichnerowicz}
If $c\in \R$, the 
operator $\DeltaL+c$ acting
on symmetric $2$-tensors 
satisfies the hypotheses of Theorem \ref{thm:main-fredholm}
if and only if 
$c>2n-n^2/4$, in which case the indicial radius of $\DeltaL+c$ is
\begin{displaymath}
R = \sqrt{ \frac{n^2}{4} - 2n + c}.
\end{displaymath}
The essential $L^2$ spectrum of $\DeltaL$ 
is $[n^2/4-2n,\infty)$.
\end{bigproposition}

(This characterization of the essential spectrum of $\DeltaL$
has also been proved by
E. Delay in \cite{Delay-spectrum},
using some of the ideas from
an earlier draft of this monograph.)

Another example is the covariant Laplacian
on trace-free symmetric tensors of any rank.

\begin{bigproposition}\label{prop:cov-Laplacian}
If $c\in \R$, 
the 
operator $\Del^*\Del+c$
acting on trace-free covariant symmetric $r$-tensors 
satisfies the hypotheses of Theorem \ref{thm:main-fredholm}
if and only if $c>-r-n^2/4$, 
in which case the indicial radius of $\Del^*\Del+c$ is
\begin{displaymath}
R = \sqrt{ \frac{n^2}{4} + r + c}.
\end{displaymath}
The essential $L^2$ spectrum of $\Del^*\Del$ 
is $[n^2/4+r,\infty)$.
\end{bigproposition}

Also, for completeness, 
we point out the following result,
originally proved in the $L^2$
case by Rafe Mazzeo 
(cf.~\cite{Mazzeo-thesis,Mazzeo-Hodge,Andersson}).

\begin{bigproposition}\label{prop:hodge-laplacian}
The Laplace-Beltrami operator acting
on $q$-forms 
satisfies the hypotheses of Theorem \ref{thm:main-fredholm}
in the following cases:
\begin{enumerate}\letters
\item 
When $0\le q<n/2$, with $R=n/2-q$.
\item
When $n/2+1<q\le n+1$, with $R = q-n/2-1$.
\end{enumerate}
In each case, the essential $L^2$ spectrum
of $\Delta$ is $[R^2,\infty)$.
\end{bigproposition}

Finally, we describe 
one significant non-Laplace
operator to which Theorem \ref{thm:main-fredholm}
applies.
The \defname{conformal Killing operator} is the operator
$L\colon C^\infty(M;TM)\to C^\infty(M;\Sigma^2_0M)$
(where $\Sigma^2_0M$ is the bundle of trace-free symmetric
$2$-tensors)
defined by letting $LV$ be the trace-free part of
the symmetrized covariant derivative of $V$
(see Chapter \ref{L2-section} for details).
A vector field $V$ satisfies $LV=0$ if and only if the
flow of $V$ preserves the conformal class of $g$.
The operator
$L^*L$, sometimes called the 
\defname{vector Laplacian}, plays an important role in 
the conformal approach to constructing initial data 
for the Einstein field equations of general relativity;
see \cite{AC,Isenberg-Park,Park-thesis}.

\begin{bigproposition}\label{prop:vector-laplacian}
The vector Laplacian $L^*L$
satisfies the hypotheses of Theorem 
\ref{thm:main-fredholm}, with $R=n/2+1$.
\end{bigproposition}

Most 
of the analytic results in 
Theorem \ref{thm:main-fredholm} 
and Propositions 
\ref{prop:lichnerowicz}--\ref{prop:vector-laplacian}
are not really new.  
The systematic  treatment of elliptic uniformly degenerate operators
dates back to the work of Mazzeo, 
building on earlier work of Richard Melrose and others
\cite{Melrose,Melrose-Mendoza,Mazzeo-Melrose,Mazzeo-thesis,Mazzeo-Hodge},
and Theorem \ref{thm:main-fredholm} can also be
derived from Mazzeo's microlocal ``edge calculus'' (cf. \cite[Theorem
6.1]{Mazzeo-edge}).
Also, $L^2$-Fredholm criteria for a general class of operators including the ones
considered here have been obtained by 
Robert Lauter, Bertrand Monthubert, and Victor Nistor
\cite[Thm.\ 4]{LMN}, using the theory of pseudodifferential operators
on groupoids.
For many purposes in geometric analysis, however, 
it is useful to have a more
``low-tech'' approach that does not use pseudodifferential operators.
An elementary approach to uniformly degenerate operators
based on a priori estimates has been used 
by Graham and Lee \cite{Graham-Lee}, 
Lars Andersson and Piotr Chrusciel \cite{Andersson,AC}, 
Johan \rade\  \cite{Rade}, and
Michael Anderson \cite{Anderson-closed,Anderson-prescribed,Anderson-remarks}.

The exposition I  present here is based on this low-tech approach, and consists of
three main ingredients: sharp a priori $L^2$ estimates,
a decay estimate for the hyperbolic Green kernel using
the spherical symmetry of the ball model, and a technique 
due to \rade\  \cite{Rade} 
for piecing together a parametrix out of this
model Green kernel.  
One advantage of this approach is that it 
deals quite naturally with operators whose coefficients are
not smooth up to the boundary, a feature that is crucial for the 
application to Einstein metrics, because the metrics around which
we linearize generally have only finite 
boundary regularity.  
It is worth noting that for many Laplace operators,
the Green kernel estimates are needed only for extending
the sharp Fredholm results to H\"older and
$L^p$ spaces; sharp Fredholm results 
in weighted $L^2$
spaces can be obtained by a much more elementary proof
based solely on a priori estimates.  See Chapter \ref{L2-section}
for details.

The main new results in this monograph 
are the Bochner-type formula 
\eqref{integral-formula} for Laplacians on tensor-valued 
differential forms, which yields
a completely elementary proof of sharp Fredholm theorems in weighted $L^2$
spaces and 
the identification of the 
essential spectrum
(Chapter \ref{L2-section}); 
new sufficient curvature conditions for the invertibility of
the Lichnerowicz Laplacian 
(part \eqref{part:yamabe} of Theorem \ref{thm:einstein});  
and the construction of Einstein metrics in the
optimal $C^{n-1,\beta}$ H\"older class up to the boundary when the
conformal infinity is sufficiently smooth.

Theorem \ref{thm:main-fredholm}
can be viewed
as complementary to the regularity results of Andersson and Chru{\'s}ciel; in
particular, for an operator of the type considered here, 
Proposition \ref{prop:optimal-regularity} 
below can be used to obtain sharp ``regularity
intervals'' in the sense defined in \cite{AC}, and then 
the results of \cite{AC} can be used to derive
asymptotic expansions for solutions to equations of the form
$Pu=f$ when $f$ and the metric are sufficiently smooth.  
We leave it to
the interested reader to work out the details.
See \cite{CDLS} for example,
where these ideas play a central role in the proof of boundary regularity for 
asymptotically hyperbolic Einstein metrics.
Similar results can also be obtained using the edge calculus
of \cite{Mazzeo-edge}.

The assumption of a weak a priori $L^2$ 
estimate \eqref{eq:asymptotic-L2-estimate}
near the boundary is used primarily to rule out
$L^2$ kernel of the model operator on hyperbolic space.  
For almost all the examples treated here, this hypothesis
is equivalent to positive indicial radius
(see Propositions
\ref{prop:lichnerowicz}, 
\ref{prop:cov-Laplacian},
\ref{prop:hodge-laplacian},
and \ref{prop:vector-laplacian}
above).  However, it is possible for \eqref{eq:asymptotic-L2-estimate}
to fail even when the indicial radius is positive---an 
example is given by the Laplace-Beltrami operator
on differential forms of degree $k$ on a manifold of dimension $2k$,
which is not Fredholm despite the fact that 
its indicial radius is $1/2$
(see 
\cite{Mazzeo-thesis,Mazzeo-Hodge}
and Lemma \ref{lemma:laplace-beltrami-ind-radius}
below).  Thanks are due to Robin Graham for pointing out
the significance of this example.

A word about the history of this monograph is in order.
The result of Theorem \ref{thm:einstein}
was originally
announced in preliminary form at an AMS meeting in 1991. 
Shortly after that meeting, I discovered a gap in the 
proof, and 
stopped work on the paper until I found a way to
bridge the gap, sometime around 1998.  By that time,
administrative responsibilities kept me away from 
writing until early 2001.  The first complete version of this monograph
was posted on {\tt www.arxiv.org} in May 2001; the present version
is a minor modification of that one.

After this monograph was nearly finished, the
beautiful monograph \cite{Biquard} by 
Olivier Biquard came to my attention.  
There is substantial overlap between the
results of \cite{Biquard} and this monograph---in 
particular, Biquard proves a perturbation result
for Einstein metrics similar to 
Theorem \ref{thm:einstein}, as well as
many of the results of Theorem 
\ref{thm:main-fredholm} in the special case
of Laplace operators, using methods very similar
to those used here.  (He also proves much more,
extending many of the same results to 
Einstein metrics that are asymptotic to
the complex, quaternionic, and octonionic
hyperbolic metrics.)  
On the other hand, many of the results of this monograph are stronger
than those of \cite{Biquard}, notably the
curvature assumptions of in Theorem \ref{thm:einstein} and 
the boundary regularity of the resulting Einstein metrics.

In Chapter \ref{mobius-section} of this monograph, we give the
main definitions, and 
describe  special ``M\"obius coordinate charts'' on an
asymptotically hyperbolic manifold that relate the geometry
to that of hyperbolic space.  In Chapter \ref{spaces-section}, we
introduce our weighted Sobolev and H\"older spaces and
prove some of their basic properties, and in Chapter
\ref{elliptic-operator-section} we prove some basic
mapping properties of geometric elliptic operators on
these spaces.   The core of the analysis begins in Chapter
\ref{model-section}, where we undertake a detailed study
of the Green kernels for  elliptic geometric operators on
hyperbolic space.  This is then applied in Chapter
\ref{section:fredholm} to construct a parametrix for
an arbitrary operator satisfying the hypotheses of 
Theorem \ref{thm:main-fredholm}.  Using this analysis,
we prove Theorem \ref{thm:main-fredholm}
and Proposition \ref{prop:L2-Fredholm}.
In Chapter  \ref{L2-section}, we explore how these
results apply in detail to Laplace operators,
and prove Propositions \ref{prop:lichnerowicz},
\ref{prop:cov-Laplacian}, 
\ref{prop:hodge-laplacian}, and
\ref{prop:vector-laplacian}.  Finally, in Chapter 
\ref{section:einstein} we
construct asymptotic solutions to 
\eqref{regularized-einstein} using a delicate procedure
that does not lose any boundary regularity, and use
these together with Theorem \ref{thm:main-fredholm} to
prove Theorem \ref{thm:einstein}.

Among the many people to whom I am indebted for
inspiration and good ideas while this work was 
in progress, I would particularly like to 
express my thanks
to Lars Andersson, Piotr
Chru{\'s}ciel, Robin Graham, Jim Isenberg, Rafe Mazzeo, 
Dan Pollack, and John Roth.  
I also would like to apologize to
them and to all who expressed interest in this work for the long delay
between the first announcement of these results and the appearance of
this monograph.
Finally, I am indebted to the referee for a number of useful suggestions.

%%
%%
%% Fredholm operators and Einstein metrics
%% on conformally compact manifolds
%%
%% by John M. Lee
%% 
%% Chapter 2
%%

\chapter{M\"obius Coordinates}\label{mobius-section}

Let $\Mbar$ be a 
smooth, compact, $(n+1)$-dimensional manifold-with-boundary,
$n\ge 1$, and $M$ its interior.  A \defname{defining function} will mean
a function $\rho\colon\Mbar\to \R$ of class at least $C^1$ that  is
positive in $M$, vanishes on $\dm$, and has nonvanishing differential
everywhere on $\dm$.  
We choose a fixed smooth defining function $\rho$ once and for all.
For any
$\epsilon>0$, let $A_\epsilon\subset M$ denote the open subset where
$0<\rho<\epsilon$.  

A Riemannian metric $g$ on $M$ is said to be \defname{conformally
compact of class ${C^{l,\beta}}$} for a nonnegative integer
$l$ and $0\le\beta< 1$ if for any smooth defining function $\rho$,
the conformally rescaled metric $\rho^2 g$ has a $C^{l,\beta}$
extension, denoted by $\gbar$, to a positive definite tensor field on
$\Mbar$.  For such a metric $g$, the induced boundary metric $\ghat :=
{\left. \gbar\right|}_{T\dm}$ is a $C^{l,\beta}$ Riemannian metric on
$\dm$ whose conformal class $[\ghat]$ is independent of the choice of
smooth defining function $\rho$; this conformal class is called the
\defname{conformal infinity} of $g$. 

Throughout this monograph, we will use the Einstein summation convention,
with Roman indices $i,j,k,\dots$ running from $1$ to $n+1$ and 
Greek indices $\alpha,\beta,\gamma,\dots$ running from $1$ to $n$.
We indicate components of 
covariant derivatives of a tensor field by indices
preceded by a semicolon, as in $u_{ij;kl}$.
In component calculations, 
we will always assume that a fixed conformally compact metric $g$
has been chosen, and all covariant
derivatives and index raising and lowering operations will be
understood to be with respect to $g$ unless otherwise
specified, except that $\overline g^{ij}$
denotes the {\it inverse} of the metric $\overline g = \rho^2 g$,
not its raised-index version.
Our convention for the components of the curvature tensor is
chosen so that the Ricci tensor is given by the contraction
$R_{ik} = R_{ijk}{}^j$.

An important fact about conformally compact metrics is that their
local geometry near the boundary looks asymptotically very much like that of
hyperbolic space.  Mazzeo \cite{Mazzeo-thesis,Mazzeo-Hodge} showed, for
example, that if $g$ is conformally compact of class at least $C^{2,0}$,
then $g$ is complete and has sectional curvatures uniformly approaching
$-{\left|d\rho\right|}^2_{\gbar}$ near $\dm$.  Thus, if $g$ is conformally
compact of class $C^{l,\beta}$ with $l\ge 2$, and
${\left|d\rho\right|}^2_{\gbar}=1$ on $\dm$, we say $g$ is
\defname{asymptotically hyperbolic of class ${C^{l,\beta}}$}.

In fact, the relationship with hyperbolic space can be made even more
explicit by constructing special coordinate charts near the boundary.
Throughout the rest of this chapter, we assume given a fixed
metric $g$ on $M$ that is asymptotically hyperbolic 
of class $C^{l,\beta}$, with
$l\ge 2$ and $0\le\beta< 1$.  Let $\gbar=\rho^2 g$, a $C^{l,\beta}$ metric
on $\Mbar$, and let $\ghat$ denote the restriction of $\gbar$ to
$T\dm$.  

We begin by choosing a covering of a neighborhood of  
$\dm$ in $\overline M$ by finitely many smooth
coordinate charts $(\Omega,\Theta)$, where 
each coordinate map $\Theta$
is of the form $\Theta=(\theta,\rho)=(\theta^1,\dots,\theta^n,\rho)$
and extends to a neighborhood of $\overline \Omega$ in $\overline M$.  
In keeping with our index convention, we will sometimes
denote $\rho$ by 
$\theta^{n+1}$, and a symbol with a Roman
index such as 
$\theta^i$ can 
refer to any of the coordinates $\theta^1,\dots,\theta^n,\rho$.

We fix once and for all finitely many such charts
covering a neighborhood $\scr W$ 
of $\dm$ in $\overline M$.
We will call any of 
these charts ``background coordinates'' for
$\Mbar$.  By compactness, there is a positive number $c$ such that
$A_{c}\subset\scr W$, and such that
every point $p\in A_{c}$ is 
contained in a background coordinate chart containing a set of
the form 
\begin{equation}\label{eq:bkg-chart-set}
\{(\theta,\rho): |\theta-\theta(p)|<c, 0\le \rho<c\}.
\end{equation}

We will use two models of hyperbolic space, depending on context.
In the upper half-space model, we regard hyperbolic space
as the open upper half-space $\H = \H^{n+1}\subset\R^{n+1}$, with
coordinates $(x,y)=(x^1,\dots,x^n,y)$, and with the hyperbolic
metric $\hyp$ given in coordinates
by $\hyp = y^{-2}\sum_i (dx^i)^2$.  
(As above, $x^i$ can denote any of the coordinates
$x^1,\dots,x^n,x^{n+1}=y$.)
The other model is the Poincar\'e ball model, in which we
regard hyperbolic space as the open unit ball 
$\B = \B^{n+1}\subset\R^{n+1}$,
with coordinates
$(\xi^1,\dots,\xi^{n+1})$,
and with the hyperbolic metric (still denoted by $\hyp$)
given by 
$\hyp=4(1-|\xi|)^{-2}\sum_i (d\xi^i)^2$,
where $|\xi|$
denotes the Euclidean norm.

In this chapter, we will work exclusively with the upper
half-space model.
For any $r>0$, we let $B_r\subset\H$ denote the hyperbolic
geodesic ball of
radius $r$ about the point $(x,y) = (0,1)$:
\begin{equation*}
B_r = \{(x,y)\in\H: d_{\hyp}((x,y),(0,1))<r\}.
\end{equation*}
It is easy to check by direct computation that 
\begin{equation*}
B_r \subset \{(x,y): |x|<\sinh r,\ e^{-r} < y < e^r\},
\end{equation*}
where $|x|$ denotes the Euclidean norm of $x\in\R^n$.

If $p_0$ is any point in $A_{c/8}$, choose such a background
chart containing $p_0$, and define a 
map $\Phi_{p_0}\colon B_2\to M$, called a
\defname{M\"obius chart} centered at $p_0$, by
\begin{displaymath}
(\theta,\rho) = \Phi_{p_0}(x,y) = (\theta_0+\rho_0 x,\rho_0 y),
\end{displaymath}
where $(\theta_0,\rho_0)$ are the background coordinates
of $p_0$.
(It is more convenient in this context to consider a
``chart'' to be a mapping from $\H\subset\R^{n+1}$ into $M$,
rather than from $M$ to $\R^{n+1}$ as is more
common.)
Because $\rho_0 < c/8$ and 
$e^2<8$, 
$\Phi_{p_0}$ maps
$B_2$ diffeomorphically onto 
a neighborhood of $p_0$ in $A_{c}$.
In these coordinates, $p_0$ corresponds to the point 
$(x,y)=(0,1)\in \H$.
For each $0<r\le 2$, let $V_r(p_0)\subset A_c$
be the neighborhood of $p_0$ defined by
\begin{equation*}
V_r(p_0) = \Phi_{p_0}(B_r).
\end{equation*}
We also choose finitely many smooth coordinate charts
$\Phi_i\colon B_2\to M$ such that the sets $\{\Phi_i(B_1)\}$
cover a neighborhood of
$M\setminus A_{c/8}$, and such that each chart $\Phi_i$
extends smoothly to 
a neighborhood of $\overline B_2$.  
For consistency, 
we will also call these
``M\"obius charts.''

The following lemma shows that the geometry
of $(M,g)$ is uniformly bounded in M\"obius charts.

\begin{lemma}\label{properties-of-mobius-coordinates}
There exists a constant $C > 0$ such that
if $\Phi_{p_0}\colon B_2\to M$ is any M\"obius chart,
\begin{align}
\|\Phi_{p_0}^*g - \hyp\|_{C^{l,\beta}(B_2)} &\le C,
\label{eq:g-est-in-mobius}\\
\sup_{B_2} |(\Phi_{p_0}^*g)^{-1}\hyp| &\le C.
\label{eq:g-inv-est-in-mobius}
\end{align}
\end{lemma}

(The H\"older and sup norms 
in this estimate are the usual norms applied to the
components of a tensor in coordinates; since $\overline B_2$ is
compact, these are equivalent to the intrinsic H\"older
and sup norms
on tensors with respect to the hyperbolic metric.)

\begin{proof}
The estimate is immediate
for the finitely many charts covering the interior of $M$,
so we need only consider M\"obius charts near the boundary.
In background coordinates,
$g$ can be written
\begin{displaymath}
g = \rho^{-2}\gbar_{ij}(\theta,\rho)d\theta^i d\theta^j.
\end{displaymath}
Pulling
back to $\H$, we obtain
\begin{align*}
\Phi_{p_0}^{*} g - \hyp &= (\rho_0 y)^{-2}\gbar_{ij} 
(\theta_0+\rho_0 x, \rho_0 y) d(\rho_0 x^i)\,
d(\rho_0 x^j)-y^{-2} \delta_{ij} dx^i\,dx^j\\
&= y^{-2}(\gbar_{ij} (\theta_0+\rho_0 x, \rho_0 y)-\delta_{ij}) dx^i
dx^j.
\end{align*}
Since $y$ and $\delta_{ij}$ are
smooth functions bounded above and below
together with all derivatives 
on $\overline B_2$, it suffices to estimate
$\gbar_{ij} (\theta_0+\rho_0 x, \rho_0 y)$.
Our choice of background coordinates ensures
that the eigenvalues of 
$\overline g_{ij}$ are uniformly bounded above and below
by a global constant.  
Uniform estimates on 
the derivatives of $\Phi_{p_0}^* g$ follow  from 
the fact that 
\begin{displaymath}
\del_{x^{i_1}} \dotsm \del_{x^{i_m}} (\Phi_{p_0}^* \overline g)_{ij}
= \rho_0^m \Phi_{p_0}^*(\del_{\theta^{i_1}} \dotsm \del_{\theta^{i_m}}
\overline g_{ij}).
\end{displaymath}
Finally, an easy computation yields
uniform H\"older estimates for the
$l$th derivatives of $\Phi_{p_0}^*g$.
(See Lemma \ref{lemma:bdry-mobius-charts} below for a sharper
estimate.)
\end{proof}

The next lemma is a version of the well-known Whitney covering lemma
adapted to the present situation.

\begin{lemma}\label{covering-by-mobius-coordinates}
There exists a countable collection of points
$\{p_i\}\subset M$ and corresponding M\"obius 
charts $\Phi_i=\Phi_{p_i}\colon B_2\to V_2(p_i)\subset M$
such that the sets $\{V_1(p_i)\}$ cover $M$
and the sets $\{V_2(p_i)\}$ are uniformly locally finite:
There exists an integer $N$ such that for each $i$,
$V_{2}({p_i})$ 
has nontrivial intersection with 
$V_{2}({p_j})$ for at most $N$ values of $j$.
\end{lemma}

\begin{proof}
We only need to show that there exist points $\{p_i\} \subset A_{c/8}$ such that $\{V_1(p_i)\}$
cover $A_{c/8}$ and $\{V_2(p_i)\}$ are uniformly locally finite, for then we can choose finitely
many additional charts for the interior without disturbing the uniform local finiteness.

By the preceding lemma, there are positive numbers $r_0 < r_1$ such that each
set $V_1(p)$ contains the $g$-geodesic ball of radius $r_0$ about $p$, 
and each set $V_2(p)$ is
contained in the geodesic ball of radius $r_1$. Let $\{p_i\}$ be any maximal collection of
points in $A_{c/8}$ such that the open geodesic balls $\{B_{r_0/2}(p_i)\}$ are disjoint. (Such a
maximal collection exists by an easy application of Zorn's lemma.) If $p$ is any point
in $A_{c/8}$, by the maximality of the set $\{p_i\}$, $B_{r_0/2}(p)$ must intersect at least one of the
balls $B_{r_0/2}(p_i)$ nontrivially, which implies that $p\in B_{r_0} (p_i) \subset V_1(p_i)$
by the triangle
inequality. Therefore the sets $\{V_1(p_i)\}$ cover $A_{c/8}$.
To bound the number of sets $\{V_2(p_i)\}$ that can intersect, it suffices to bound the
number of geodesic balls of radius $r_1$ around points $p_i$ that can intersect. Let $i$ be
arbitrary and suppose $B_{r_1} (p_j)\cap B_{r_1} (p_i) \ne \emptyset$ 
for some $j$. By the triangle inequality
again, $B_{r_0/2}(p_j) \subset B_{2r_1+r_0/2}(p_i)$. Since 
$M$ has bounded sectional curvature, standard
volume comparison theorems (see, for example, 
\cite[Theorems 3.7 and 3.9]{Chavel}
yield uniform volume estimates
\begin{displaymath}
\Vol(B_{r_0/2}(p_j)) \ge C_1, \qquad \Vol(B_{2r_1+r_0/2}(p_j)) \le C_2.
\end{displaymath}
Since the sets $B_{r_0/2}(p_j)$ are disjoint for different values of $j$, there can be at most
$C_2/C_1$ such points $p_j$.
\end{proof}

%%
%%
%% Fredholm operators and Einstein metrics
%% on conformally compact manifolds
%%
%% by John M. Lee
%% 
%% Chapter 3
%%

\chapter{Function Spaces}\label{spaces-section}

In this chapter,
we will define weighted Sobolev and
H\"older spaces of tensor fields
that are well adapted to the geometry
of asymptotically hyperbolic manifolds.
Similar spaces have been defined by other
authors; for some examples, 
see \cite{AC,Andersson,Delay,Delay-Herzlich,Graham-Lee}.

Throughout this chapter, we 
assume $\overline M$ is a 
connected smooth $(n+1)$-manifold,
$g$ is a metric on $M$ that is
asymptotically hyperbolic of class $C^{l,\beta}$,
with $l\ge 2$ and $0\le\beta< 1$, and $\rho$ is a fixed smooth
defining function for $\del M$.  (It is easy to verify
that choosing another smooth defining function will replace the norms
we define below by equivalent ones, and will leave the 
function spaces unchanged.)

A \defname{geometric tensor bundle} over $\overline M$
is a subbundle $E$ of some tensor bundle $T^{r_1}_{r_2}\overline M$ 
(tensors of covariant rank $r_1$ and contravariant rank $r_2$)
associated to a direct summand (not necessarily
irreducible) of the standard 
representation of $\Ortho(n+1)$ (or $\SO(n+1)$ if $M$
is oriented)
on tensors of type $\binom{r_1}{r_2}$ over $\R^{n+1}$.
We will also use the same symbol $E$ to denote the restriction
of this bundle to $M$.
We define the \defname{weight} of such a bundle
$E\subset T^{r_1}_{r_2}\overline M$ to be $r=r_1-r_2$.  
The significance of this definition lies in the way tensor
norms scale conformally: 
With $\overline g = \rho^2 g$, an
easy computation shows that
\begin{displaymath}
|T|_g = \rho^r |T|_{\overline g}\quad\text{for all $T\in T^{r_1}_{r_2}$}.
\end{displaymath}

We begin by defining H\"older spaces of functions
and tensor fields 
that are continuous up to the boundary.
For $0\le \alpha<1$ and $k$ a nonnegative integer,
we let $C^{k,\alpha}_{(0)}(\overline M)$ denote the usual
Banach space of functions on $\overline M$ with $k$ derivatives
that are H\"older continuous of degree $\alpha$ up to the
boundary in each
background coordinate chart, with the obvious norm.
(When $\alpha=0$, this is just the usual space of
functions that are $k$ times continuously differentiable on
$\overline M$.)
If $s$ is a real number such that 
$0\le s \le k+\alpha$,
we define a subspace $C^{k,\alpha}_{(s)}(\overline M)
\subset
C^{k,\alpha}_{(0)}(\overline M)$
by 
\begin{displaymath}
C^{k,\alpha}_{(s)}(\overline M) = \{u\in 
C^{k,\alpha}_{(0)}(\overline M): u = O(\rho^s)\}.
\end{displaymath}
The next lemma gives some elementary properties of these
spaces.

\begin{lemma}\label{lemma:properties-of-vanishing-spaces}
Suppose $0\le\alpha<1$ and 
$0\le s \le k+\alpha$.  
\begin{enumerate}\letters
\item\label{part:vanishing=derivs}
$C^{k,\alpha}_{(s)}(\overline M) = \{u\in 
C^{k,\alpha}_{(0)}(\overline M): \del^i u/\del\rho^i |_{\del M} = 0
\text{ for $0\le i < s$}\}$.
\item\label{part:vanishing-closed}
$C^{k,\alpha}_{(s)}(\overline M)$ is a
closed subspace of $C^{k,\alpha}_{(0)}(\overline M)$.
\item\label{part:vanishing-loc-const}
If $j$ is a positive integer
and $j-1<s\le j\le k$, then $C^{k,\alpha}_{(s)}(\overline M)
=
C^{k,\alpha}_{(j)}(\overline M)$.  
\item\label{part:vanishing-loc-const-alpha}
If $k<s\le k+\alpha$, then $C^{k,\alpha}_{(s)}(\overline M)
=
C^{k,\alpha}_{(k+\alpha)}(\overline M)$.
\item\label{part:derivs-map-vanishing}
If $u\in C^{k,\alpha}_{(s)}(\overline M)$ for $s\ge 1$, then any 
background coordinate derivative $\del_iu$ is in 
$C^{k-1,\alpha}_{(s-1)}(\overline M)$.
\item\label{part:vanishing-mult-pos}
If $\delta$ is a positive real number 
such that $s+\delta\le k+\alpha$, then
$\rho^\delta C^{k,\alpha}_{(s)}(\overline M)\subset
C^{k,\alpha}_{(s+\delta)}(\overline M)$.
\item\label{part:vanishing-mult-neg}
If $j$ is an integer such that $0<j\le s$, then
$\rho^{-j}C^{k,\alpha}_{(s)}(\overline M)\subset C^{k-j,\alpha}_{(s-j)}(\overline M)$.
\end{enumerate}
\end{lemma}

\begin{proof}
All of these claims are local, so we fix one background
coordinate chart $(\theta,\rho)$ and do
all of our computations there.  
Let $m$ be any nonnegative integer less than or equal to $k$.
Applying the one-variable version of 
Taylor's formula to $u\in C^{k,\alpha}_{(0)}(\overline M)$, we obtain
\begin{equation}\label{eq:Taylor}
u(\theta,\rho) 
=\sum_{i=0}^{m-1}\frac{1}{i!}\rho^i\frac{\del^i u}{\del\rho^i} (\theta,0)
+ \frac{1}{(m-1)!} \rho^m \int_0^1 (1-t)^{m-1}\frac{\del^m u}{\del\rho^m}(\theta,t\rho)\,dt.
\end{equation}
The integral above is easily shown to define a 
function of $(\theta,\rho)$ that is in
$C^{k-m,\alpha}_{(0)}$ up to the boundary and agrees with
$\del^m u/\del\rho^m$ on $\del M$; therefore the last term in \eqref{eq:Taylor} is
$O(\rho^m)$ in general, and if $\alpha>0$, it is
$O(\rho^{m+\alpha})$ if and only if 
$\del^m u/\del\rho^m$ vanishes on $\del M$.
Part \eqref{part:vanishing=derivs} follows easily 
from this, 
and 
\eqref{part:vanishing-closed}, \eqref{part:vanishing-loc-const},
\eqref{part:vanishing-loc-const-alpha}, and  
\eqref{part:derivs-map-vanishing}
follow from
\eqref{part:vanishing=derivs}.
Part \eqref{part:vanishing-mult-pos}
is an immediate consequence of the definition,
and \eqref{part:vanishing-mult-neg} follows by setting
$m=j$ in \eqref{eq:Taylor} and 
multiplying through by $\rho^{-j}$.
\end{proof}

Because of part \eqref{part:vanishing-closed}
of this lemma, we can consider
$C^{k,\alpha}_{(s)}(\overline M)$
as a Banach space with the norm inherited from
$C^{k,\alpha}_{(0)}(\overline M)$.

If $E$ is a geometric tensor bundle over $\overline M$,
we extend this definition to spaces of tensor fields
by letting $C^{k,\alpha}_{(s)}(\overline M;E)$ be the
space of tensor fields whose components in each background
coordinate chart are in 
$C^{k,\alpha}_{(s)}(\overline M)$.  All the 
claims of Lemma 
\ref{lemma:properties-of-vanishing-spaces}
extend immediately to tensor fields.

Next we define some spaces of tensor fields over the
interior manifold $M$ associated with the asymptotically
hyperbolic metric $g$.  Let us start
with the Sobolev spaces, for which the definitions are a bit
simpler.  
First, for $1< p < \infty$ and $k$ a nonnegative integer
less than or equal to $l$,
we define $H^{k,p}(M;E)$ to be the usual (intrinsic) $L^p$ Sobolev
space determined by the metric $g$: That is, $H^{k,p}(M;E)$ is the
Banach space
of all locally integrable sections $u$ of $E$ such that $\Del^j
u$ (interpreted in the distribution sense) is in $L^p(M;E\tprod
T^jM)$ for $0\le j \le k$, with the norm
\begin{displaymath}
\|u\|_{k,p} = \bigg(\sum_{j=0}^k \int_M |\Del^j u|^p\,
  dV_g\bigg)^{1/p}.
\end{displaymath}
(Here $dV_g$ is the Riemannian density.)
In the special case of $H^{0,2} (M;E) = L^2(M;E)$, we will denote the
norm $\|\cdot\|_{0,2}$ and its associated inner product simply by
$\|\cdot\|$ and $(\cdot,\cdot)$, respectively.  (We reserve the
notations $|\cdot|_g$ and $\left <\cdot,\cdot\right>_g$ for the pointwise
norm and inner product on tensors.)

For each real number $\delta$, we
define the \defname{weighted Sobolev space} $H^{k,p}_\delta(M;E)$ by
\begin{displaymath}
H^{k,p}_\delta(M;E) := \rho^\delta H^{k,p}(M;E) =
\{\rho ^\delta u: u \in H^{k,p}(M;E)\}
\end{displaymath}
with norm
\begin{displaymath}
\|u\|_{k,p,\delta} := \|\rho^{-\delta}u\|_{k,p}.
\end{displaymath}
These are easily seen to be Banach spaces, and
$H^{k,2}_\delta(M;E)$ is a Hilbert space.  
Note that H\"older's inequality implies that when $1/p+1/p^*=1$,
$H^{0,p^*}_{-\delta}(M;E)$ is naturally isomorphic
to the dual space of $H^{0,p}_\delta(M;E)$
under the $L^2$ pairing $(u,v) = \int_M\<u,v\>_g\,dV_g$.

The following lemma is elementary but often useful.

\begin{lemma}\label{lemma:being-in-Lp}
Let $(M,g)$ be a conformally compact
$(n+1)$-manifold, let $\rho$
be a defining function for $M$, and let $u$ be a continuous
section of a natural tensor bundle $E$ of weight $r$ on $M$.
\begin{enumerate}\letters
\item
If $|u|_{g}\le C\rho^s$ with $s$ real and greater than $\delta+n/p$,
then $u\in H^{0,p}_\delta(M;E)$.
\item
If $|u|_{g}\ge C\rho^s>0$ on the complement of
a compact set,
with $s\in\R$ and $s\le \delta+n/p$,
then $u\notin H^{0,p}_\delta(M;E)$.
\item\label{part:overlineu-in-Lp}
If $u=\rho^s\overline u$, where $s\in\C$ and
$\overline u$ is 
continuous on $\overline M$
and does not vanish identically on $\dm$, then
$u\in H^{0,p}_\delta(M;E)$
if and only if $\Real s>\delta+n/p-r$.
\end{enumerate}
\end{lemma}

\begin{proof}
Let $\overline g = \rho^2 g$, which is a continuous
Riemannian metric on $\overline M$.
The lemma follows from the easily-verified facts that
$|\rho^s \overline u|_g = \rho^{\Real s+r}|\overline u|_{\overline g}$,
$dV_g = \rho^{-n-1}dV_{\overline g}$, and
$\int_M \rho^s \, dV_{g}<\infty$
if and only if $\Real s>n$.
\end{proof}

Next we turn to the weighted H\"older spaces.
To define H\"older norms for tensor fields on a manifold, one is always faced
with the problem of comparing values of a tensor field at nearby points
in order to make sense of
quotients of the form $|u(x)-u(y)|/d(x,y)^\alpha$.
There are various intrinsic ways to do this, such as parallel
translating $u(x)$ along a geodesic from $x$ to $y$, or comparing the
components of $u(x)$ and $u(y)$ in Riemannian normal coordinates
centered at $x$; on a noncompact manifold, these can yield different
spaces depending on the behavior of the metric near infinity.  We will
adopt a definition in terms of the M\"obius coordinates constructed
above which, though not obviously
intrinsic to the geometry of $(M,g)$, has the
advantage that estimates for elliptic operators in these norms follow
very easily from standard 
local elliptic estimates in
M\"obius coordinates.

Let $\alpha$ be a real number such that $0\le \alpha <1$,
and let $k$ be a nonnegative integer such that $k+\alpha\le l+\beta$.
For any tensor field $u$ with locally
$C^{k,\alpha}$ coefficients, define the norm $\|u\|_{k,\alpha}$ by
\begin{equation}\label{eq:def-holder-norm}
\|u\|_{k,\alpha} :=
\sup_\Phi\|\Phi^*u\|_{C^{k,\alpha}(B_2)},
\end{equation}
where $\|v\|_{C^{k,\alpha}(B_2)}$ is just the usual Euclidean H\"older
norm of the components of $v$ on $B_2\subset\H$, 
and the supremum is over all 
M\"obius charts defined on $B_2$.  
Let $C^{k,\alpha}(M;E)$ be the space of sections of
$E$ for which this norm is finite.
(We distinguish notationally between the
H\"older and Sobolev norms as follows: A Greek subscript in the second
position takes values in the interval $[0,1)$, and the notation
$\|u\|_{k,\alpha}$ denotes a H\"older norm; a Roman subscript takes
values in $(1,\infty)$, and the notation $\|u\|_{k,p}$ denotes a Sobolev
norm.  We will avoid the ambiguous cases $\alpha=1$ and $p=1$.)
The corresponding
\defname{weighted H\"older spaces} are defined for $\delta\in\R$ by
\begin{displaymath}
C^{k,\alpha}_\delta(M;E) := \rho^\delta C^{k,\alpha}(M;E) =
\{\rho ^\delta u: u \in C^{k,\alpha}(M;E)\}
\end{displaymath}
with norms
\begin{displaymath}
\|u\|_{k,\alpha,\delta} := \|\rho^{-\delta}u\|_{k,\alpha}.
\end{displaymath}

If $U\subset M$ is a subset, the restricted norms are denoted by
$\|\cdot\|_{k,p,\delta;U}$ and $\|\cdot\|_{k,\alpha,\delta;U}$, and the spaces
$H^{k,p}_\delta(U;E)$ and $C^{k,\alpha}_\delta(U;E)$ are the spaces of sections
over $U$ for which these norms are finite.  When $E$ is the trivial line
bundle (i.e.\ when the tensor fields in question are just scalar
functions), we omit the bundle from the notation: For example,
$H^{k,p}(M)$ is the Sobolev space of scalar functions on $M$ with $k$
covariant derivatives in $L^p(M)$.

It is obvious from the definitions that $\rho^\delta\in C^{k,\alpha}_\delta(M)$
for every $k$.  More importantly, we have the following lemma.

\begin{lemma}\label{delrho}
If $0\le j\le l$, then $\Del^j\rho\in C^{l-j,\beta}_1(M;T^{j}M)$.
\end{lemma}

\begin{proof}
From the definition of $C^{k,\beta}_1(M;T^jM)$, we need to show 
for any M\"obius chart $\Phi_{p_0}$
that the coefficients of
$\Phi_{p_0}^*(\rho^{-1}\Del^j\rho)$
are in
$C^{l-j,\beta}(B_2)$, with norm bounded independently of $p_0$.
Since 
$\Phi_{p_0}^*\rho = \rho(p_0) y$,
this follows immediately from the 
coordinate expression for $\Phi_{p_0}^*(\rho^{-1}\Del^j\rho) =
y^{-1}(\Del_{\Phi_{p_0}^*g})^j y$ and
the fact that the Christoffel symbols of $\Phi_{p_0}^*g$ in M\"obius
coordinates are uniformly bounded in $C^{l-1,\beta}(B_2)$.
\end{proof}

We have defined the weighted spaces by multiplying the standard Sobolev
and H\"older spaces by powers of $\rho$.  In many circumstances, it is
more convenient to use an alternative characterization in terms of
weighted norms of covariant derivatives of $u$, as given in the
following lemma.

\begin{lemma}\label{characterize-weighted-spaces}
Let $u$ be a locally integrable section of a tensor bundle $E$ over an
open subset $U\subset M$.
\begin{enumerate}\letters
\item\label{sobolev-case}
For $1<p<\infty$ and
$0\le k\le l$,
$u\in H^{k,p}_\delta(U;E)$ if and only if
$\rho^{-\delta}\Del^j u\in L^{p}(U;E\tprod T^jM)$ for $0\le j\le k$, and the
$H^{k,p}_\delta$ norm is equivalent to
\begin{displaymath}
\sum_{0\le j\le k} \|\rho^{-\delta}\Del^j u\|_{0,p;U}.
\end{displaymath}
\item\label{holder-case} If $0\le\alpha<1$ and $0< k+\alpha\le l+\beta$,
$u\in C^{k,\alpha}_\delta(U;E)$ if and only if $\rho^{-\delta}\Del^j u\in
C^{0,\alpha}(U;E\tprod T^jM)$ for $0\le j\le k$, and the
$C^{k,\alpha}_\delta$ norm is equivalent to
\begin{displaymath}
\sum_{0\le j\le k} \sup_U |\rho^{-\delta}\Del^j u| + 
  \|\rho^{-\delta} \Del^k u\|_{0,\alpha;U}
\end{displaymath}
\end{enumerate}
\end{lemma}

\begin{proof}
First consider part \eqref{sobolev-case}.  By definition,  $u\in
H^{k,p}_\delta(U;E)$ iff $\Del^j(\rho^{-\delta}u)\in L^{p}(U;E\tprod T^jM)$
for $0\le j\le k$.  By the product rule and induction, we can write
\begin{equation}\label{product-rule}
\Del^j (\rho^{-\delta}u) =
\rho^{-\delta}
\hspace{-2ex}
\sum_{\substack{
0\le i \le j\\
{i+j_1+\dots+j_p = j}
}}
\hspace{-2ex}
C(\delta,i,j_1,\dots,j_p)
  \Del^{i} u\tprod
  \frac{\Del^{j_1}\rho}{\rho}\tprod\dots\tprod
  \frac{\Del^{j_{p}}\rho}{\rho}
  ,
\end{equation}
where $C(\delta,i,j_1,\dots,j_p)$ is a constant, equal to 1 when $i=j$.
Since, by the preceding lemma, $|\Del^{j}\rho|/\rho$ is bounded as long
as $j\le l$, the result follows easily by induction.

For part \eqref{holder-case}, the case $\delta=0$ follows by inspecting the
coordinate expression for $\Del^ju$ in M\"obius coordinates, recalling
that the Christoffel symbols of $g$ are uniformly bounded in
$C^{l-1,\beta}$ in these coordinates.  The
general case follows as above from \eqref{product-rule} and the fact
that $\Del^{j}\rho/\rho\in C^{l-j,\beta}(U;T^jM)$.
\end{proof}

It follows from Lemma \ref{properties-of-mobius-coordinates} that the
norm $\|\cdot\|_{k,0}$
is equivalent to the usual intrinsic $C^k$ norm $\sum_{0\le i\le
k}\sup_M |\Del^iu|$ for $0\le k \le l$.
Moreover, if $\Phi$ and $\tw\Phi$ are any M\"obius charts
whose images intersect, the overlap map $\Phi^{-1}\circ\tw\Phi$
induces an isomorphism on the space
$C^{l,\beta}(U)$ on its domain
of definition $U$, with norm bounded by a 
uniform constant.  Therefore, for $0\le k+\alpha\le l+\beta$,
it is immediate that 
the pullback by any M\"obius chart 
induces isomorphisms on the unweighted H\"older 
spaces, with 
norms bounded above and below
by constants independent of $i$.
Moreover, 
if $\Phi_i$ is a M\"obius chart centered at $p_i$,
$\Phi_i^*\rho = \rho(p_i)y$, which 
is uniformly bounded above and below on $B_2$ 
by constant multiples of $\rho(p_i)$. 
Therefore
the weighted norms scale as follows:
for any $r\le 2$,
\begin{equation}\label{eq:mobius-uniform-holder}
C^{-1}\rho(p_i)^{-\delta}\|\Phi_i^*u\|_{k,\alpha;B_r}
\le \|u\|_{k,\alpha,\delta;V_r(p_0)}
\le C\rho(p_i)^{-\delta}\|\Phi_i^*u\|_{k,\alpha;B_r}.
\end{equation}
Similarly, it follows directly from Lemma 
\ref{properties-of-mobius-coordinates} that
\begin{equation}\label{eq:mobius-uniform-sobolev}
C^{-1}\rho(p_i)^{-\delta}\|\Phi_i^*u\|_{k,p;B_r}
\le \|u\|_{k,p,\delta;V_r(p_0)}
\le C\rho(p_i)^{-\delta}\|\Phi_i^*u\|_{k,p;B_r}.
\end{equation}
The next lemma is a strengthening of this.

\begin{lemma}\label{lemma:norms-by-covering}
Suppose $\{\Phi_i=\Phi_{p_i}\}$ is a uniformly locally finite
cover of $M$ by M\"obius charts as in Lemma
\ref{covering-by-mobius-coordinates}.
Then we have the following norm equivalences
for any $r$ with $1\le r\le 2$:
\begin{align}
C^{-1}\sum_i\rho(p_i)^{-\delta}\|\Phi_i^*u\|_{k,p;B_r}
&\le \|u\|_{k,p,\delta}
\le C\sum_i\rho(p_i)^{-\delta}\|\Phi_i^*u\|_{k,p;B_r},
\label{eq:sobolev-norm-equiv}
\\
C^{-1}\sup_i\rho(p_i)^{-\delta}\|\Phi_i^*u\|_{k,\alpha;B_r}
&\le \|u\|_{k,\alpha,\delta}
\le C\sup_i\rho(p_i)^{-\delta}\|\Phi_i^*u\|_{k,\alpha;B_r},
\label{eq:holder-norm-equiv}
\end{align}
\end{lemma}

\begin{proof}
If $\Phi_i$ is a M\"obius chart centered at $p_i$,
$\Phi_i^*\rho = \rho(p_i)y$, which 
is uniformly bounded above and below
by constant multiples of $\rho(p_i)$. 
Thus if $u\in H^{k,p}_\delta(M;E)$,
\eqref{eq:mobius-uniform-sobolev} yields
\begin{align*}
\sum_i\rho(p_i)^{-\delta}\|\Phi_i^*u\|_{k,p;B_r}
&\le C\sum_i\|\Phi_i^*(\rho^{-\delta}u)\|_{k,p;B_r}\\
&\le C'\sum_i\|\rho^{-\delta}u\|_{k,p;V_r(p_i)}\\
&\le C'N \|\rho^{-\delta}u\|_{k,p}\\
&=  C'N \|u\|_{k,p,\delta},
\end{align*}
where $N$ is an upper bound on the number of 
sets $V_2(p_i)$ that can intersect nontrivially.
Conversely, if $\sum_i\rho(p_i)^{-\delta}\|\Phi_i^*u\|_{k,p;B_r}$
is finite,
then 
\begin{align*}
\|u\|_{k,p,\delta}
&\le \sum_i\|\rho^{-\delta}u\|_{k,p;V_r(p_i)}\\
&\le C\sum_i\|\Phi_i^*(\rho^{-\delta}u)\|_{k,p;B_r}\\
&\le C'\sum_i\rho(p_i)^{-\delta}\|\Phi_i^*u\|_{k,p;B_r}.
\end{align*}
The argument for the H\"older case is similar but
simpler, because we do not need the uniform local 
finiteness in that case.
\end{proof}

The following results follow easily from 
Lemma
\ref{lemma:norms-by-covering}
together with standard
facts about Sobolev and H\"older spaces 
on $B_r\subset\H$, so the proofs are left to the reader
(cf.\ \cite{Graham-Lee,Andersson}).

\begin{lemma}\label{properties-of-spaces}
Let $U$ be any open subset of $M$, and let $E,E_1,E_2$ be geometric
tensor bundles over $M$.
\begin{enumerate}\letters
\item\label{part:mult-cts}
If $1<p<\infty$,
$0\le\alpha<1$, $\delta,\delta'\in\R$, 
$0\le k+\alpha\le l+\beta$, 
and $1\le k'+\alpha\le l+\beta$,
the pointwise tensor product induces  continuous maps 
\begin{align*}
C^{k,\alpha}_\delta(U;E_1)
\cross H^{k,p}_{\delta'}(U;E_2)&\to
H^{k,p}_{\delta+\delta'}(U;E_1\otimes E_2),\\
C^{k',\alpha}_\delta(U;E_1)
\cross C^{k',\alpha}_{\delta'}(U;E_2)&\to
C^{k',\alpha}_{\delta+\delta'}(U;E_1\otimes E_2).
\end{align*}
\item\label{inclusions}
If $1<p,p'<\infty$, $0\le\alpha<1$, and
$0\le k+\alpha\le l+\beta$, 
we have continuous inclusions:
\begin{align*}
H^{k,p}_\delta(U;E) &\mapsinto H^{k,p'}_{\delta'}(U;E),
&p\ge p',\quad \delta+\frac n p &> \delta' + \frac{n}{p'};\\
C^{k,\alpha}_\delta(U;E) &\mapsinto H^{k,p'}_{\delta'}(U;E),
& \delta &> \delta' + \frac{n}{p'}.
\end{align*}
\item{\sc (Sobolev embedding)}
\label{prop:Sobolev}
If $1<p<\infty$, $0<\alpha<1$, $1\le k \le l$, $k+\alpha\le l+\beta$,
and $\delta\in\R$, we have
continuous inclusions
\begin{align*}
H^{k,p}_\delta(U;E) &\mapsinto C^{j,\alpha}_\delta(U;E),&&
k-\frac{n+1}{p}\ge j+\alpha,\\
H^{k,p}_\delta(U;E) &\mapsinto H^{j,p'}_\delta(U;E),&&
k-\frac{n+1}{p}\ge j-\frac{n+1}{p'}.
\end{align*}
\item{\sc (Rellich lemma)}
If $1<p<\infty$, $0<\alpha<1$, $0\le k \le l$, 
and $0<j+\alpha\le l+\beta$,
then the following inclusions
are compact operators:
\begin{align*}
H^{k,p}_\delta(M;E)&\mapsinto H^{k',p}_{\delta'}(M;E),&k>k',\quad\delta&>\delta',\\
C^{j,\alpha}_\delta(M;E)&\mapsinto C^{j',\alpha}_{\delta'}(M;E),&j>j',\quad\delta&>\delta'.
\end{align*}
\end{enumerate}
\end{lemma}

The following relationships between the H\"older
spaces on $M$ and those on $\overline M$ will play an 
important role in Chapter \ref{section:einstein}.

\begin{lemma}\label{lemma:spaces-m-and-mbar}
Let $E$ be a geometric tensor bundle of weight $r$
over $\overline M$, and 
suppose $0<\alpha<1$, $0<k+\alpha\le l+\beta$,
and $0\le s \le k+\alpha$.  The following
inclusions are continuous.
\begin{enumerate}\letters
\item\label{part:holder-mbar-into-m}
$C^{k,\alpha}_{(s)}(\overline M;E)
\mapsinto
C^{k,\alpha}_{s+r}(M;E)$.
\item\label{part:holder-m-into-mbar}
$C^{k,\alpha}_{k+\alpha+r}(M;E)
\mapsinto
C^{k,\alpha}_{(0)}(\overline M;E)$.
\end{enumerate}
\end{lemma}

\begin{proof}
We will prove  \eqref{part:holder-mbar-into-m} by
induction on $k$.  
Suppose $k=0$, and 
consider first the case of scalar functions, so that
$E$ is the trivial line bundle.
By
Lemma 
\ref{lemma:properties-of-vanishing-spaces}\eqref{part:vanishing-loc-const}
there are only two distinct cases: $s=0$ and $s=\alpha$.
For $s=0$, let $\Phi_{p_0}$ be a M\"obius chart
and let $(\rho_0,\theta_0)$
be the background coordinates of $p_0$.
We estimate
\begin{align}
|\Phi_{p_0}^*u(x,y)|
&\le \sup_M |u|\notag\\
&\le \|u\|_{C^{0,\alpha}_{(0)}(\overline M)}.
\label{eq:pullback-sup}\\
|\Phi_{p_0}^*u(x,y) - \Phi_{p_0}^*u(x',y')|
&= |u(\theta_0+\rho_0 x,\rho_0 y) - u(\theta_0+\rho_0 x',\rho_0 y')|\notag\\
&\le \|u\|_{C^{0,\alpha}_{(0)}(\overline M)}
|(\rho_0 x,\rho_0 y) - (\rho_0 x',\rho_0 y')|^\alpha\notag\\
&\le \|u\|_{C^{0,\alpha}_{(0)}(\overline M)}
\rho_0^\alpha|( x, y) - ( x', y')|^\alpha.
\label{eq:pullback-holder}
\end{align}
Since $\rho_0^\alpha$ is uniformly bounded on $M$, the result
follows.  

When  $s=\alpha$,
we need to show that $\Phi_{p_0}^*(\rho^{-\alpha} u)$
is uniformly bounded in $C^{0,\alpha}(B_2)$.
The H\"older estimate for $u$ in background coordinates, together
with the fact that $u$ vanishes on $\del M$,
shows that 
\begin{displaymath}
|u(\theta,\rho)|\le \|u\|_{C^{0,\alpha}_{(0)}(\overline M)}\rho^\alpha
\end{displaymath}
near $\del M$,
from which the 
zero-order
estimate 
\begin{displaymath}
|\Phi_{p_0}^*(\rho^{-\alpha}u)|\le C\|u\|_{C^{0,\alpha}_{(0)}(\overline M)}
\end{displaymath}
follows immediately.  The H\"older estimate is proved
by noting that $\Phi_{p_0}^*(\rho^{-\alpha}u)
=  y^{-\alpha}\rho_0^{-\alpha}\Phi_{p_0}^*u$; since
$y^{-\alpha}$ is uniformly bounded together with
all derivatives on $B_2$, it suffices to show that 
$\rho_0^{-\alpha}\Phi_{p_0}^*u$ is uniformly bounded in $C^{0,\alpha}(B_2)$.
This is proved
as follows:
\begin{align*}
|\rho_0^{-\alpha}\Phi_{p_0}^*u(x,y) -
\rho_0^{-\alpha}\Phi_{p_0}^*u(x',y')|
&= \rho_0^{-\alpha} |u(\theta_0+\rho_0 x,\rho_0 y) -u(\theta_0+\rho_0 x,\rho_0 y)|\\
&\le \rho_0^{-\alpha}\|u\|_{C^{0,\alpha}_{(0)}(\overline M)}
|(\rho_0 x,\rho_0 y) -(\rho_0 x,\rho_0 y)|^\alpha\\
&= \|u\|_{C^{0,\alpha}_{(0)}(\overline M)}|(x,y)-(x',y')|^\alpha.
\end{align*}

Now let $E$ be 
a tensor bundle of type $\binom p q$ (and thus of weight $r=p-q$).
In any background coordinate domain 
$\Omega$, 
basis tensors of the form $d\theta^{i_1}\otimes \dotsm \otimes
d\theta^{i_p} \otimes \del_{\theta^{j_1}}\otimes
\dotsm \otimes \del_{\theta_{j_q}}$ 
are easily seen to be uniformly bounded in
$C^{l,\beta}_{r}(\Omega;E)$.  Since any
$u\in C^{0,\alpha}_{(s)}(\overline M;E)$ can be
written locally as a linear combination of such tensors multiplied
by functions in $C^{0,\alpha}_{(s)}(\Omega)$,
the result follows from the scalar case together
with Lemma
\ref{properties-of-spaces}\eqref{part:mult-cts}.

Now suppose $1<k+\alpha\le l+\beta$ and $0\le s \le k+\alpha$,
and assume claim \eqref{part:holder-mbar-into-m}
is true for $C^{k_0,\alpha}_{(s_0)}(\overline M;E)$
when $0\le k_0<k$ and $0\le s_0\le k_0+\alpha$.
If $u\in C^{k,\alpha}_{(s)}(\overline M;E)$, then
$|u|_g = \rho^r|u|_{\overline g} = O(\rho^{s+r})$ 
by definition, so it suffices to
prove that $\Del u\in C^{k-1,\alpha}_{s+r}(M;E\otimes T^*M)$
with norm bounded by $\|u\|_{C^{k,\alpha}_{(s)}(\overline M;E)}$.
If $D=\Del-\overline\Del$
denotes the 
difference tensor between the Levi-Civita
connections of $g$ and $\overline g$,
a computation shows that
$D$ is the $3$-tensor whose components
are given in any coordinates by
\begin{equation}\label{eq:difference-tensor-components}
D_{ij}^k = - \rho^{-1} (\delta^k_i \del_j \rho + \delta^k_j \del_i\rho
- \overline g^{kl} \overline g_{ij} \del_l\rho).
\end{equation}
Working in background coordinates, we find that
$\rho D u \in C^{k,\alpha}_{(s)}(M;E\otimes T^*M)$ since
the coefficients of $\rho D$ are in $C^{l,\beta}_{(0)}(\overline M)$.
Since 
$\overline \Del u\in C^{k-1,\alpha}_{(s-1)}(M;E\otimes T^*M)$
by Lemma 
\ref{lemma:properties-of-vanishing-spaces}\eqref{part:derivs-map-vanishing},
we use Lemma 
\ref{lemma:properties-of-vanishing-spaces}\eqref{part:vanishing-mult-pos}
to 
conclude that 
\begin{align*}
\rho\Del u &= \rho\overline\Del u + \rho D u\\
&\in 
\rho C^{k-1,\alpha}_{(s-1)}(\overline M;E\otimes T^*M)
+C^{k,\alpha}_{(s)}(\overline M;E\otimes T^*M)\\
&\subset C^{k-1,\alpha}_{(s)}(\overline M;E\otimes T^*M)
\end{align*}
with norm bounded by a multiple of $\|u\|_{C^{k,\alpha}_{(s)}(\overline M;E)}$.
Therefore, since $\rho\Del u$ is
a tensor field of weight $r+1$, the inductive hypothesis implies that
$\rho\Del u\in C^{k-1,\alpha}_{s+r+1}(M;E\otimes T^*M)$,
which implies in turn that
$\Del u\in C^{k-1,\alpha}_{s+r}(M;E\otimes T^*M)$
as desired.  

For part 
\eqref{part:holder-m-into-mbar}, 
we  begin with the scalar case, and
proceed  by induction on $k$.
Let $k=0$, and 
suppose 
$u\in C^{0,\alpha}_{\alpha}(M)$.  
Given any M\"obius chart $\Phi_{p_0}\colon
B_2\to V_2(p_0)$,
let $(\theta_0,\rho_0)$ be
the background coordinates of $p_0$, and 
let $v$ be the function $\Phi^*(\rho^{-\alpha}u)$ on $B_2$,
so that 
\begin{displaymath}
u(\theta,\rho) = \rho^\alpha v((\theta-\theta_0)/\rho_0, \rho/\rho_0).
\end{displaymath}
The hypothesis means that $v\in C^{0,\alpha}(B_2)$, with
norm bounded by $\|u\|_{0,\alpha,\alpha}$.
For $(\theta,\rho)\in V_2(p_0)$, we
estimate
\begin{equation}\label{eq:sup-u}
|u(\theta,\rho)|
= |\rho^\alpha v((\theta-\theta_0)/\rho_0, \rho/\rho_0)|
\le \rho^\alpha\|u\|_{0,\alpha,\alpha}.
\end{equation}
Since $\rho$ is bounded above on $M$, this shows in particular that $\sup|u|$ is bounded
by a multiple of $\|u\|_{0,\alpha,\alpha}$.
Similarly, if $(\theta,\rho)$ and $(\theta',\rho')$ are in $V_2(p_0)$, we have
\begin{equation}\label{eq:holder-close}
\begin{aligned}
|u(\theta,\rho) - u(\theta',\rho')|
&= |\rho^\alpha v((\theta-\theta_0)/\rho_0, \rho/\rho_0)
 - {\rho'}^\alpha v((\theta'-\theta_0)/\rho_0, \rho'/\rho_0)|\\
&\le \rho^\alpha|v((\theta-\theta_0)/\rho_0, \rho/\rho_0)
 - v((\theta'-\theta_0)/\rho_0, \rho'/\rho_0)|\\
&\qquad+(\rho^\alpha - {\rho'}^\alpha)|v((\theta'-\theta_0)/\rho_0,
\rho'/\rho_0)|\\
&\le C\|u\|_{0,\alpha,\alpha} |(\theta,\rho)-(\theta',\rho')|^\alpha.
\end{aligned}
\end{equation}

To complete the $k=0$ case, it suffices to extend \eqref{eq:holder-close}
to any $(\theta,\rho)$ and $(\theta',\rho')$ that lie
in the same background chart, 
for then $u$  extends continuously to the
boundary as an element of $C^{0,\alpha}_{(0)}(\overline M)$.
Note that there is a real number $\gamma\in(1,2)$ such that
whenever $|(\theta,\rho)-(\theta',\rho')|\le \gamma\rho$, 
the points $(\theta,\rho)$ and $(\theta',\rho')$ lie in the 
image of the same M\"obius chart.
We estimate as follows:
\begin{equation}\label{eq:breakup-u}
|u(\theta,\rho) - u(\theta',\rho')|
\le 
|u(\theta,\rho) - u(\theta',\rho)| +
|u(\theta',\rho) - u(\theta',\rho')|.
\end{equation}
For the first term, the case in which $|\theta-\theta'|\le\gamma\rho$ is
taken care of by \eqref{eq:holder-close}.  On the other hand,
if $|\theta-\theta'|\ge \gamma\rho$, \eqref{eq:sup-u} gives
\begin{align*}
|u(\theta,\rho) - u(\theta',\rho)|
&\le |u(\theta,\rho)| + |u(\theta',\rho)|\\
&\le 2\rho^\alpha \|u\|_{0,\alpha,\alpha}\\
&\le 2\gamma^{-\alpha}\|u\|_{0,\alpha,\alpha} |\theta-\theta'|^\alpha\\
&\le C\|u\|_{0,\alpha,\alpha} |(\theta,\rho)-(\theta',\rho')|^\alpha.
\end{align*}
To estimate the second term of \eqref{eq:breakup-u}, 
let $N$ be a positive
integer such that $\rho'$ lies in the interval $[\gamma^{N}\rho,\gamma^{N+1}\rho]$.  
Then, since $(\theta',\gamma^i\rho)$ and $(\theta',\gamma^{i+1}\rho)$ lie in 
the image of a single M\"obius chart, as do $(\theta',\gamma^N\rho)$ and $(\theta',\rho')$, 
\eqref{eq:holder-close} gives
\begin{align*}
|u(\theta',\rho') - u(\theta',\rho)|
&\le |u(\theta',\rho') - u(\theta',\gamma^N\rho)| +
\sum_{i=0}^{N-1} |u(\theta',\gamma^{i+1}\rho) - u(\theta',\gamma^{i}\rho)|\\
&\le C\|u\|_{0,\alpha,\alpha} \left(|\rho'-\gamma^N\rho|^\alpha
+
\sum_{i=0}^{N-1} |\gamma^{i+1}\rho - \gamma^{i}\rho|^\alpha\right)\\
&=  C\|u\|_{0,\alpha,\alpha} \left(|\rho'-\gamma^N\rho|^\alpha
+ (\gamma-1)^\alpha\rho^\alpha
\sum_{i=0}^{N-1} \gamma^{\alpha i}\right)\\
&=  C\|u\|_{0,\alpha,\alpha} \left(|\rho'-\gamma^N\rho|^\alpha
+ (\gamma-1)^\alpha\rho^\alpha
\frac{\gamma^{\alpha N} -1}{\gamma^\alpha-1}\right)\\
&\le  C'\|u\|_{0,\alpha,\alpha} \left(|\rho'-\gamma^N\rho|^\alpha
+ |\gamma^N\rho - \rho|^\alpha\right),
\end{align*}
where the last inequality follows from the fact that 
$\gamma^{\alpha N} -1 \le C(\gamma^N-1)^\alpha$ when $N\ge 1$.
Since both terms in parentheses above are bounded by $|\rho - \rho'|^\alpha$
and thus by $|(\theta,\rho)-(\theta',\rho')|^\alpha$,
this completes the argument for the $k=0$ case. 
(I am indebted to Eric Bahuaud for pointing out a gap in an earlier
version of this proof, and providing helpful suggestions for fixing it.)

For $k\ge 1$, to show that 
$u\in C^{k,\alpha}_{(0)}(\overline M)$,
it suffices to show that $u$ is  bounded
and $Vu\in C^{k-1,\alpha}_{(0)}(\overline M)$
whenever $V$ is a smooth vector field on $\overline M$.   
It is straightforward to check that 
any such vector field
maps $C^{k,\alpha}_{k+\alpha}( M)$
into $C^{k-1,\alpha}_{k-1+\alpha}( M)$.
Thus if  $u\in
C^{k,\alpha}_{k+\alpha}(M)$, then 
$V u \in C^{k-1,\alpha}_{k-1+\alpha}(M)
\subset C^{k-1,\alpha}_{(0)}(\overline M)$ by induction,
with norm bounded by $\|u\|_{k,\alpha,k+\alpha}$.
Since $|u|$ is obvious uniformly bounded by $\|u\|_{k,\alpha,k+\alpha}$,
it follows that $u\in C^{k,\alpha}_{(0)}(\overline M)$.

Finally, let $E$ be a $\binom p q$ tensor bundle.
If $u\in C^{k,\alpha}_{k+\alpha+r}(M;E)$, 
to show that $u\in C^{k,\alpha}_{(0)}(\overline M;E)$ we have to
show that the components of $u$ in any background coordinate
domain $\Omega$
are in $C^{k,\alpha}_{(0)}(\Omega)$.
 These components are given by complete contractions of
tensor products of $u$
with tensors of the form $\del_{\theta^{i_1}}\otimes
\dotsm \otimes \del_{\theta^{i_p}} \otimes d\theta^{j_1}
\otimes \dotsm \otimes d\theta^{j_q}$, each of which is
in $C^{l,\beta}_{-r}(\Omega;T^q_pM)$.  By 
Lemma \ref{properties-of-spaces}\eqref{part:mult-cts},
these tensor products are in 
$C^{k,\alpha}_{k+\alpha}(\Omega; E\otimes T^q_pM)$.
Since complete contraction clearly maps 
this space into 
$C^{k,\alpha}_{k+\alpha}(\Omega)$,
the result for tensor fields follows from the scalar case.
\end{proof}

\begin{lemma}\label{psi-epsilon}
Let $\psi\colon  \R^+ \to [0,1]$ be a smooth function that is 
equal to $1$ on
$[0,\half]$ and supported in $[0,1)$, and set $\psi_\epsilon(q) =
\psi(\rho(q)/\epsilon)$ for $q\in M$.  
If  $0\le\alpha<1$ and $0\le k+\alpha \le l+\beta$,
then $\psi_\epsilon\in C^{k,\alpha}(M)$, with norm bounded independently of
$\epsilon$.
\end{lemma}

\begin{proof}
Working directly with the definition of $C^{k,\alpha}(M)$,
for any M\"obius coordinate chart $\Phi_{p_0}$, we have to show that
$\Phi_{p_0}^*\psi_\epsilon(x,y) = \psi(\rho(p_0)y/\epsilon)$ is uniformly
bounded in $C^{k,\alpha}(B_2)$.  Since 
$e^{-2} < y < e^2$ on $B_2$,
$\Phi_{p_0}^*\psi_\epsilon$ will be identically equal to 0 or 1 on $B$
unless $\half\epsilon e^{-2} < \rho(p_0) < \epsilon e^2$.  Under these restrictions,
it is easy to verify that $\psi(\rho(p_0)y/\epsilon)$ is uniformly bounded
in $C^{k,\alpha}(B_2)$.
\end{proof}

\begin{lemma}
If $1<p<\infty$, 
$\delta\in \R$, and $0\le k \le l$, the set of compactly supported smooth
tensor fields is dense in $H^{k,p}_\delta(M;E)$.
\end{lemma}

\begin{proof}
Suppose $u \in H^{k,p}_\delta(M;E)$.
First we show that $u$ can be approximated in the
$H^{k,p}_\delta$ norm by compactly
supported elements of $H^{k,p}_\delta(M;E)$.

Let $\psi_\epsilon$ be as in the preceding lemma.  We will show that
$(1-\psi_\epsilon) u \to u$ in $H^{k,p}_\delta$ as $\epsilon
\to 0$, which by Lemma \ref{characterize-weighted-spaces} is the same as
$\Del ^j(\psi_\epsilon u) \to 0$ in $H^{0,p}_{\delta}$ whenever $0 \le j
\le k$.  By the product rule,
\begin{displaymath}
\Del ^j(\psi_\epsilon u) = \sum_{i=0}^j C_i \Del ^{j-i}
\psi_\epsilon\tprod \Del ^i u.
\end{displaymath}
Since $|\Del ^{j-i} \psi_\epsilon|_g$ is bounded and supported where
$\rho \le 2\epsilon$, we have
\begin{displaymath}
\|\Del ^{j-i}
\psi_\epsilon \otimes \Del ^i u\|^p_{0,p,\delta} \le
C\int\limits_{\{\rho \le 2\epsilon\}}
| \rho^{-\delta} \Del ^i u|_g ^p
\,dV_g .
\end{displaymath}
By Lemma \ref{characterize-weighted-spaces}, $| \rho^{-\delta} \Del ^i
u|_g ^p$ is integrable, so the integral on the right-hand side above goes
to zero as $\epsilon \to 0$ by the dominated convergence theorem.

Next we must check that if $u \in H^{k,p}_\delta(M;E)$ is compactly supported
in $M$, it can be approximated in the $H^{k,p}_\delta$ norm by smooth,
compactly supported tensor fields.  The classical argument involving
convolution with an approximate identity shows that $u$ can be
approximated in the standard Sobolev $H^{k,p}$ norm on a slightly larger
compact set by tensor fields in $C^\infty_c(M;E)$.  However, on any
fixed compact subset of $M$, it is easy to see that the $H^{k,p}_\delta$ norm
is equivalent to the $H^{k,p}$ norm, thus completing the proof.
\end{proof}

%%
%%
%% Fredholm operators and Einstein metrics
%% on conformally compact manifolds
%%
%% by John M. Lee
%% 
%% Chapter 4
%%

\chapter{Elliptic Operators}\label{elliptic-operator-section}

In this chapter, we collect some basic facts regarding elliptic
operators on asymptotically hyperbolic manifolds.  Throughout this chapter we
assume $(M,g)$ is a connected
asymptotically hyperbolic $(n+1)$-manifold of class $C^{l,\beta}$ for some
$l\ge 2$ and $0\le\beta< 1$, and $\rho$ is a fixed smooth defining
function.

Let $E$ and $F$ be geometric 
tensor bundles over $M$.  We will say a linear
partial differential operator 
$P\colon C^\infty(M;E)\to C^\infty(M;F)$
is a \defname{geometric operator of order $m$} 
if the components of $Pu$ in any coordinate
frame are given 
by linear functions of 
the components of $u$
and their partial derivatives, 
whose coefficients
are universal 
polynomials in the 
components of the 
metric $g$, their partial derivatives, and the function
$(\det g_{ij})^{-1}$ (or 
$(\det g_{ij})^{-1/2}$ if $M$ is oriented),
such that the coefficient of any $j$th
derivative of $u$ involves at most the first $m-j$ derivatives of 
the metric.  
A geometric operator is, in particular, an example of
a \defname{regular natural differential operator} in the 
terminology introduced by P. Stredder \cite{Stredder}.
We note that such an operator is automatically invariant
under (orientation-preserving) isometries.
By the results of \cite{Stredder}, $P$ is geometric of order $m$
if and only if $Pu$ can be written as a sum of 
contractions of tensors of the form
\begin{multline}\label{eq:P-tensor-prod}
\Del^j u
 \otimes 
\Del^{k_1}Rm \otimes \dots \otimes \Del^{k_l} Rm
 \otimes\\
\underbrace{g \otimes \dots \otimes g}_{\text{$p$ factors}} 
 \otimes 
\underbrace{g^{-1} \otimes\dots \otimes g^{-1}}_{\text{$q$ factors}}
 \otimes 
\underbrace{dV_g \otimes\dots \otimes dV_g}_{\text{$s$ factors}}
,
\end{multline}
(possibly after reordering indices),
with $0\le j \le m$ and $0\le k_i \le m-j-2$,
and with $s=0$ unless $M$ is oriented.
Here $Rm$ is the (covariant) Riemann curvature
tensor of $g$, and $dV_g$ is its Riemannian volume form.

If $E$ and $F$ are tensor bundles over $\overline M$,
 a differential operator $P\colon C^\infty(M;E)\to 
C^\infty(M;F)$ 
is said to be 
\defname{uniformly degenerate} if
it can
be written 
in background coordinates 
as a system of operators that are polynomials in
$\rho\partial/\partial \theta^i$ with coefficients that are at
least continuous
up to the boundary.  

\begin{lemma}
Let $(M,g)$ be an
asymptotically hyperbolic $(n+1)$-manifold
of class $C^{l,\beta}$,
$0\le\beta<1$.
Suppose $E$ is a geometric 
tensor bundle over $M$, and $P\colon C^\infty(M;E)
\to C^\infty(M;E)$ is a geometric operator of order $m\le l$.  
Then $P$ is
uniformly degenerate.
\end{lemma}

\begin{proof}
As noted above, for any section $u$ of $E$, 
$Pu$ can be written as 
a sum of contractions of terms like 
\eqref{eq:P-tensor-prod}.
If $u$ 
is covariant of degree $r_1$ and contravariant of degree $r_2$,
then this tensor product has 
$r_2 + 2q$ upper indices and $r_1 + j + 2p + (n+1)s  + 4l+ \sum_i k_i$
lower indices.  
(The $4l$ lower indices are the undifferentiated indices of the $l$ copies of $Rm$.)
Because we are assuming that $Pu$ is the
same type of tensor as $u$,  
$2q$ of the upper indices must be contracted against
$j + 2p + (n+1)s  + 4l + \sum_i k_i$ of the lower indices, so
in particular we must have 
\begin{equation}\label{eq:add-up-indices}
2q = j + 2p + (n+1)s  + 4l + \sum_i k_i.
\end{equation}

It is obvious that 
tensoring with $\overline g=\rho^2 g$, 
$\overline g^{-1}=\rho^{-2}g^{-1}$, and
$dV_{\overline g} = \rho^{n+1}dV_{g}$ are all uniformly degenerate
operators.
Using formula \eqref{eq:difference-tensor-components}
for the components of the difference tensor
$D=\Del-\overline\Del$, we see that the components of
$\rho D_{ij}^k$ in background coordinates
are $C^{l,\beta}$ up to the boundary.
It follows that  
$(\rho \Del)^j = (\rho\overline\Del + \rho D)^j$ 
is a uniformly degenerate operator for $0\le j\le l$.
Since $[\Del,\rho]u = u\otimes d\rho$ is also
uniformly degenerate, it follows by induction
that $\rho^j\Del^j$
is uniformly degenerate as well.

A straightforward 
computation (cf.\ \cite{Mazzeo-Hodge,Graham-Lee}) shows that 
the components of $Rm$ are given by
\begin{equation}\label{eq:Rijkl}
R_{ijkl} = - |d\rho|_{\overline g}^2 (g_{ik}g_{jl} - g_{il}g_{jk})
+ \rho^{-3}p_1(\overline g, \overline g^{-1},\del \overline g) 
+ \rho^{-2}p_2(\overline g, \overline g^{-1},\del \overline g,\del^2 \overline g),
\end{equation}
where $p_1$ and $p_2$ are universal polynomials,
so $\rho^4 Rm\in C^{l-2,\beta}_{(0)}(\overline M,T^4\overline M)$.
Moreover, an easy induction shows that for $k\le l-2$,
$\rho^{4+k}\Del^k Rm \in 
C^{l-2-k,\beta}_{(0)}(\overline M,T^4\overline M)$.
It follows that
tensoring with $\rho^{4+k}\Del^k Rm$ is a uniformly degenerate operator.
By virtue of \eqref{eq:add-up-indices}, therefore,
we can rewrite 
\eqref{eq:P-tensor-prod}
in the following manifestly uniformly degenerate form:
\begin{multline}\label{eq:P-tc-tensor-prod}
\rho^j\Del^j u 
 \otimes 
\rho^{4+k_1}\Del^{k_1}Rm \otimes \dots \otimes \rho^{4+k_l}\Del^{k_l} Rm
 \otimes
\rho^2 g \otimes \dots \otimes \rho^2 g
 \otimes 
\\
\rho^{-2} g^{-1} \otimes\dots \otimes \rho^{-2} g^{-1}
 \otimes 
\rho^{n+1} dV_g \otimes\dots \otimes\rho^{n+1} dV_g .
\end{multline}
Since contraction of a lower index against an upper
one is also uniformly degenerate, the result follows.
\end{proof}

If $P$ is a
uniformly degenerate operator,
for each $s\in \C$ we define the 
\defname{indicial map} of $P$ to be the
bundle map $I_s(P)\colon E|_{\del M} \to E|_{\del M}$
defined by
\begin{displaymath}
I_s(P)(\overline u) = (\rho^{-s} P(\rho^s\overline u))|_{\dm},
\end{displaymath}
where $\rho$ is any smooth defining function, and
$\overline u$ is extended arbitrarily to a $C^m$
section of $E$ in a neighborhood of $\dm$.
It is easy to check that the indicial map
of any uniformly
degenerate operator is a continuous
bundle map, which is independent of the
extension or the choice of defining function.
For geometric operators, we can say more.

\begin{lemma}\label{lemma:Clbeta-indicial-map}
Suppose $(M,g)$ is an 
asymptotically hyperbolic $(n+1)$-manifold of class $C^{l,\beta}$,
and $P\colon C^\infty(M;E)\to C^\infty(M;E)$ is a geometric
operator of order $m\le l$.  Then for each $s\in \C$,
$I_s(P)\colon E|_{\del M}\to E|_{\del M}$ is a
$C^{l,\beta}$ bundle map.
\end{lemma}

\begin{proof}
It is an immediate consequence of the definition of the
indicial map that if $P_1$ and $P_2$ are uniformly
degenerate operators with sufficiently smooth
coefficients, then 
\begin{align*}
I_s(P_1\circ P_2) &= I_s({P_1})\circ I_s({P_2}),\\
I_s({P_1+P_2}) &= I_s({P_1}) + I_s({P_2}).
\end{align*}
Since a sum or composition of two $C^{l,\beta}$ bundle
maps is again of class $C^{l,\beta}$, it suffices to
display $P$ as a sum of compositions of uniformly degenerate
operators
with $C^{l,\beta}$ indicial maps.  

Writing $P$ as a sum of contractions of terms
of the form \eqref{eq:P-tc-tensor-prod},
we see that it suffices
to show that each of the following uniformly degenerate
operators
has $C^{l,\beta}$ indicial map:
\begin{enumerate}\letters
\item\label{part:delmap}
$u\mapsto \rho^j\Del^j u$;
\item
$u\mapsto u\otimes \rho^{4+k}\Del^k Rm$.
\item
$u\mapsto u\otimes \rho g$;
\item
$u\mapsto u\otimes \rho^{-1} g^{-1}$;
\item
$u\mapsto u\otimes \rho^{n+1} dV_g$.
\end{enumerate}

The last three operators above 
are themselves $C^{l,\beta}$
bundle maps over $\overline M$, whose indicial maps are just their
restrictions to $\del M$, so there is nothing to
prove in those cases.  
For \eqref{part:delmap}, observe that
$\rho^j\Del^j u$ can be written as a sum of compositions
of the operators $[\rho,\Del]$ and $\rho\Del$.
Since the commutator $[\rho,\Del]u = -u\otimes d\rho$
is a smooth bundle map and therefore
has smooth indicial map,
we need only consider the operator $\rho\Del$.
Let 
$D=\Del-\overline\Del$
be the 
difference tensor 
as in the preceding proof, and observe that 
for any $\overline u\in C^1_{(0)}(\overline M;E)$,
\begin{align*}
\rho^{-s}(\rho\Del(\rho^s\overline u))
&=
\rho^{-s}(\rho\overline \Del (\rho^s \overline u) 
+ \rho D(\rho^s \overline u))\\
&= s\, \overline u \otimes d\rho+ 
\rho D(\overline u) + O(\rho).
\end{align*}
As noted above, 
$\overline u\mapsto \overline u\otimes d\rho $
is a smooth bundle map; and it follows 
from \eqref{eq:difference-tensor-components}
that
$\rho D$ is a $C^{l,\beta}$ bundle map.
Therefore the indicial map of $\rho\Del$ is $C^{l,\beta}$
as claimed.

To analyze the remaining operator, 
$u\mapsto u\otimes \rho^{4+k}\Del^k Rm$,
let $K$ be the $4$-tensor field on $M$ 
whose components are
\begin{displaymath}
K_{ijkl} = -(g_{ik}g_{jl}
- g_{il}g_{jk}),
\end{displaymath}
and let $\overline K$ be 
\begin{displaymath}
\overline K_{ijkl} = -(\overline g_{ik}\overline g_{jl}
- \overline g_{il}\overline g_{jk}) = \rho^4 K_{ijkl}.
\end{displaymath}
Because the assumption that $g$ is asymptotically
hyperbolic means that $|d\rho|_{\overline g}=1$
on $\del M$, we can write $|d\rho|_{\overline g} = 1
+\rho v$, where $v\in C^{l-1,\beta}_{(0)}(\overline M)$.
Using 
\eqref{eq:Rijkl}, therefore, we can write
\begin{equation}\label{eq:RM=K+V}
Rm = K + \rho^{-3} \overline V = \rho^{-4}\overline K + \rho^{-3} \overline V,
\end{equation}
where $\overline V\in C^{l-2,\beta}_{(0)}(\overline M,T^4\overline M)$.
It follows that $I_s(u\mapsto u\otimes \rho^4 Rm)$ is just tensoring
with $\overline K|_{\del M}$,
which is a $C^{l,\beta}$ bundle map.
On the other hand, Lemma
\ref{lemma:spaces-m-and-mbar}\eqref{part:holder-mbar-into-m}
shows that 
\begin{align*}
Rm - K &= \rho^{-3}\overline V 
\in \rho^{-3}C^{l-2,\beta}_{(0)}(\overline M;T^4\overline M)\\
&\subset \rho^{-3}C^{l-2,\beta}_{4}( M;T^4 M)
= C_{1}^{l-2,\beta}(M;T^4M). 
\end{align*}
 Because
$K$ is parallel,
$\Del^k Rm = \Del^k(Rm - K) \in 
C_{1}^{l-2-k,\beta}(M;T^4M)$.
In particular, this implies that 
\begin{equation}\label{eq:delrm-vanishing}
|\rho^{4+k}\Del^k Rm|_{\overline g} = |\Del^k Rm|_{g} = O(\rho),
\end{equation}
so $\rho^{4+k}\Del^k Rm$ has vanishing indicial map for $k>0$.
\end{proof}

A complex number $s$ is called a \defname{characteristic exponent} of
$P$ at $\hat p\in \dm$ if $I_s(P)$ is singular at $\hat p$. 
Its \defname{multiplicity} 
is the dimension of the kernel
of $I_s(P)\colon E_{\hat p}\to E_{\hat p}$.
If $s$ is a characteristic exponent of $P$ somewhere 
on $\dm$, we just say it is a characteristic exponent of $P$.

Let $(\B,\hyp)$
be the unit ball model of hyperbolic space
as described in Chapter \ref{mobius-section}, let
$\breve E$ denote the tensor 
bundle over $\B $ that corresponds to $E$ (i.e., associated
to the same representation of $\Ortho(n+1)$ or 
$\SO(n+1)$), and let
$\breve P\colon C^\infty(\B ;\breve E)\to C^\infty(\B ;\breve E)$
be the geometric operator on $\breve E$ whose coordinate
formula is the same as that of $P$.
Because $\breve P$ is invariant under (orientation-preserving)
isometries of $\B $, and Euclidean
rotations in the unit ball model of hyperbolic space
are isometries of $\hyp$ that act transitively on 
the sphere $\S^n = \del \B $, it follows that the
characteristic exponents of $\breve P$ and their
multiplicities are constant on $\del\B $.
The next lemma shows that the analogous statement holds for $P$
as well.

\begin{lemma}\label{lemma:const-char-exps}
Let $(M,g)$ be a connected asymptotically 
hyperbolic $(n+1)$-manifold of class $C^{l,\beta}$,
and let $P\colon C^\infty(M;E)\to C^\infty(M;E)$ be a
geometric operator of order $m\le l$.  
The characteristic exponents of $P$
and their multiplicities are constant on $\del M$,
and are the same as those of 
$\breve P$.
\end{lemma}

\begin{proof}
Let $\hat p\in \del M$ be arbitrary, and let $(\theta,\rho)$
be any background coordinates
defined on a neighborhood of $\hat p$ in $\overline M$.  By an affine
change of background coordinates, we may arrange that 
$\hat p = (0,0)$ and 
$\overline g_{ij} = \delta_{ij}$ at $\hat p$.
Let $\Phi\colon V\to M$ be the M\"obius chart given by
$(\theta,\rho) = \Phi(x,y) = (x,y)$ in terms of these background
coordinates, defined on some neighborhood $V$ of $(0,0)$ in $\overline\H$.
Let $\tw g$ be the Riemannian metric 
$\Phi^* g$ on $V\cap \H$, let $\tw P = \Phi^* \circ P \circ \Phi^{-1*}$,
and let $I_s({\tw P})(\overline u) = y^{-s}\tw P(y^s\overline u)$ be
the indicial map of $\tw P$.  
Then $I_s({\tw P}) = \Phi^* \circ I_s(P) \circ \Phi^{-1*}$,
so the characteristic exponents and multiplicities of
$\tw P$ at $(0,0)\in \del\H$ are the same as those
of $P$ at $\Phi(x,0)$.  Thus it suffices to show that
$I_s({\tw P}) = I_s({\breve P})$ at $(0,0)$.

For convenience, here is a summary of the several metrics being considered in this proof:
\begin{align*}
g & && \text{(the given asymptotically hyperbolic metric on $M$),}\\
\overline g &= \rho^2 g &&\text{(on $\overline M$)},\\
\tw g &= \Phi^* g &&\text{(on an open subset of $\H$)},\\
\breve g &&&\text{(the hyperbolic metric on $\H$)}.
\end{align*}

Arguing as in the proof of the preceding lemma, 
it suffices to show that each of the following operators
has vanishing indicial map at $(0,0)$:
\begin{enumerate}\letters
\item\label{factor:rhodel}
$u\mapsto \rho\tw \Del u-\rho \breve \Del u$,
\item\label{factor:g}
$u\mapsto u\otimes(\rho^2\tw g - \rho^2\hyp)$,
\item\label{factor:ginv}
$u\mapsto u\otimes( \rho^{-2}\tw g^{-1} - \rho^{-2}\hyp^{-1})$,
\item\label{factor:dv}
$u\mapsto  u\otimes(\rho^{n+1}dV_{\tw g} - \rho^{n+1}dV_{\hyp})$,
\item\label{factor:delrm}
$u\mapsto  u\otimes(\rho^{4+k}\tw\Del^k\tw{Rm} -\rho^{4+k} \breve\Del^k\breve{Rm})$.
\end{enumerate}

If we let $E= \tw\Del - \breve\Del$ be the
difference tensor between the two connections,
a computation based on \eqref{eq:difference-tensor-components}
yields
\begin{equation*}
\rho E^i_{jk} = 
+ \del_j\rho(\overline g^{il} \overline g_{kl} - \delta^{il}\delta_{kl}) + O(\rho).
\end{equation*}
Because $\overline g_{ij} = \delta_{ij}$ at $\hat p$,
it follows that 
$\rho\tw\Del -\rho\breve\Del$ is a zero-order operator
that vanishes at $(0,0)$.

The next three operators, tensoring with
\eqref{factor:g}--\eqref{factor:dv}, are obviously
zero-order operators vanishing 
at $(0,0)$, so their indicial maps have the same property.
To complete the proof, we need to show 
that the indicial map of \eqref{factor:delrm}
vanishes at $(0,0)$.  This is obvious from the
argument in the preceding proof when $k>0$, because
both 
$\rho^{4+k}\tw\Del^k\tw{Rm}$ and 
$\rho^{4+k} \breve\Del^k\breve{Rm}$ individually have 
vanishing indicial maps.
On the other hand, the preceding proof showed that 
the restriction of 
$\rho^4\tw{Rm}$
to $\del M$ is equal to $\overline K$, which
in turn is equal to $\rho^4\breve{Rm}$
at $(0,0)$.
This completes the proof.
\end{proof}

The next observation we need to make is that if
$P$ is a formally self-adjoint geometric operator acting on a
geometric tensor bundle of weight $r$, its
set of characteristic exponents is symmetric about the line 
$\Real s = n/2-r$.  This will follow easily from the next lemma.

\begin{proposition}\label{prop:Is*}
Let $(M,g)$ be an 
asymptotically hyperbolic $(n+1)$-manifold of class $C^{l,\beta}$.
Suppose that $m\le l$ and 
$P\colon C^\infty(M;E)\to C^\infty(M;E)$ is an $m$th order geometric
operator acting on sections of a geometric tensor bundle $E$
of weight $r$, and let $P^*$ denote its formal adjoint.
Then 
\begin{displaymath}
I_s({P^*}) = I_{n-2r-\overline s}(P)^*.
\end{displaymath}
\end{proposition}

\begin{proof}
Choose any point $\hat p\in \dm$, and let $(\theta,\rho)$
be background coordinates in a neighborhood $\Omega$ of $\hat p$ in $\overline M$.
Then we can write $P$ locally as
\begin{displaymath}
Pu(\theta,\rho) = \sum_{\substack{0\le k\le m\\1\le j_i,\dots,j_k\le n+1}}
a^{j_1\dots j_k}(\theta,\rho) (\rho\del_{j_1}) \dotsm (\rho\del_{j_k}) u(\theta,\rho),
\end{displaymath}
where each coefficient $a^{j_1\dots j_k}(\theta,\rho)$ is a 
matrix-valued $C^{l-m+k,\beta}$ 
function (the coordinate representation of a 
bundle map from $E$ to itself).  
Suppose now that $\overline u$ is a section of
$E$ that has a $C^m$ extension to 
$\overline M$.  Then $\rho\del_i(\rho^s\overline u)
= O(\rho^{s+1})$ if $i\ne n+1$, so 
$I_s(P)$ is given locally by
\begin{align*}
I_s(P) \overline u &= \lim_{\rho\to 0}\sum_{\substack{0\le k\le m\\j_i=\dots=j_k= n+1}}
\rho^{-s} a^{j_1\dots j_k}(\theta,\rho) (\rho\del_{n+1}) \dotsm 
(\rho\del_{n+1})(\rho^s\overline u(\theta,\rho))\\
&= \sum_{\substack{0\le k\le m\\j_i=\dots=j_k= n+1}}
s^k a^{j_1\dots j_k}(\theta,0) \overline u(\theta,0).
\end{align*}

To compute $I_s(P^*)$, we first compute the formal adjoint of 
a monomial of the form $\rho\del_i$ as follows.  Let $\overline g
= \rho^2 g$, which is a $C^{l,\beta}$ metric on $\overline M$.
The 
inner products on $E$ defined by $g$ and $\overline g$ 
are related by
$\<u,v\>_g = \rho^{2r} \<u,v\>_{\overline g}$, and
the volume elements by $dV_g = \rho^{-n-1} dV_{\overline g}$.
If $u,v$ are
smooth sections of $E$ compactly supported in $\Omega\cap M$,
then
\begin{align*}
\int_M \< u, \rho \del_i v\>_g dV_g
&= \int_M \rho^{2r-n-1}\<u, \rho\del_i v\>_{\overline g} (\det \overline g)^{1/2} \,d\theta^1\dots d\theta^{n+1}\\
&= -\int_M \<\del_i(\rho^{2r-n} u b(\overline g)),v\>_{\overline g}(\det \overline g)^{1/2}\,d\theta^1\dots d\theta^{n+1}\\
&= -\int_M \<\rho^{n+1-2r}\del_i(\rho^{2r-n}  u b(\gbar)),v\>_{g}\,dV_g,
\end{align*}
where $b(\gbar)$ is a constant-coefficient polynomial in the 
components of $\gbar$, $\gbar^{-1}$, and $(\det \overline g)^{1/2}$.
From this, it follows easily that
\begin{displaymath}
(\rho\del_i)^* = -\rho\del_i + (n-2r)\delta_i^{n+1}- \rho b_i,
\end{displaymath}
where $b_i = \del_i b(\gbar)$ is $C^{l-1,\beta}$ up to $\del M$.
Therefore, 
\begin{multline*}
P^*u(\theta,\rho) = \sum_{\substack{0\le k\le m\\1\le j_i,\dots,j_k\le n+1}}
(-\rho\del_{j_k} + (n-2r)\delta_{j_k}^{n+1}-\rho b_{j_k})
\dotsm\\
\qquad
(-\rho\del_{j_1} + (n-2r)\delta_{j_1}^{n+1}-\rho b_{j_1})
(a^{j_1\dots j_k}(\theta,\rho )^*u(\theta,\rho )),
\end{multline*}
and we conclude that
\begin{align*}
I_s(P^*)\overline u 
&= \lim_{\rho \to 0} \sum_{\substack{0\le k\le m\\j_i=\dots=j_k= n+1}}
\rho ^{-s}(-\rho \del_{n+1 } + (n-2r)-\rho b_{n+1})
\dotsm\\
&\qquad\qquad
(-\rho \del_{n+1 } + (n-2r)-\rho b_{n+1})
(a^{j_1\dots j_k}(\theta,\rho )^*(\rho ^s\overline u(\theta,\rho )))\\
&= \sum_{\substack{0\le k\le m\\j_i=\dots=j_k= n+1}}
(n-2r-s)^k a^{j_1\dots j_k}(\theta,0)^* \overline u(\theta,0)\\
&= I_{n-2r-\overline s}(P)^*\overline u,
\end{align*}
which was to be proved.
\end{proof}

\begin{corollary}\label{cor:symmetric-char-exp}
If $P$ is a formally self-adjoint geometric operator of order
$m\le l$, then the set of characteristic exponents
of $P$ is symmetric about the line $\Real s = n/2-r$.
\end{corollary}

\begin{proof}
The preceding proposition shows that $I_s(P) = I_s({P^*}) = I_{n-2r-\overline s}(P)^*$.
Thus if $s$ is a characteristic exponent of $P$, then
$I_{n-2r-\overline s}(P)^*$, and hence also $I_{n-2r-\overline s}(P)$, is singular, 
which means that 
$s'=n-2r-\overline s$ is also a characteristic exponent.
Since $\Image s = \Image s'$ and 
$\half(\Real s + \Real s') = n/2-r$, the result follows.
\end{proof}

This corollary shows that a geometric self-adjoint operator
$P$ must have at least one
characteristic exponent whose real part is greater than or equal 
to $n/2-r$.  
We define the 
\defname{indicial radius} of $P$ to be 
the smallest nonnegative number $R$ such that $P$ has
a characteristic exponent whose real part is $n/2-r+R$.
(This is well-defined because the set of characteristic
exponents is finite by 
Lemma \ref{lemma:const-char-exps}.)

Next we investigate the mapping properties of 
geometric operators on our weighted spaces.

\begin{lemma}\label{mapping-properties}
Let $P\colon C^\infty(M;E)\to C^\infty(M;F)$ be a geometric operator of
order $m$.
\begin{enumerate}\letters
\item\label{mapping-sobolev-case}
If $\delta\in \R$, $1< p<\infty$, and $m\le k\le l$, then $P$ extends
naturally to a bounded mapping
\begin{displaymath}
P\colon H^{k,p}_\delta(M;E)\to H^{k-m,p}_\delta(M;F).
\end{displaymath}
\item\label{mapping-holder-case}
If $\delta\in\R$, $0\le\alpha<1$, and $m\le k+\alpha\le l+\beta$, then $P$
extends naturally to a bounded mapping
\begin{displaymath}
P\colon C^{k,\alpha}_{\delta}(M;E)\to C^{k-m,\alpha}_{\delta}(M;F).
\end{displaymath}
\end{enumerate}
\end{lemma}

\begin{proof}
Let $\{\Phi_i\}$ be a uniformly locally finite covering
of $M$ by M\"obius charts as in 
Lemma \ref{covering-by-mobius-coordinates}.  
Let $\breve E$ be the bundle over hyperbolic space
associated to the same $\Ortho(n+1)$ or 
$SO(n+1)$ representation as $E$, and 
for each
$i$, let $P_i\colon C^\infty(B_2;E)\to C^\infty(B_2;E)$
be the operator defined by $P_i = \Phi_i^* \circ P \circ
(\Phi_i^{-1})^*$.
Since the metric $g_i = \Phi_i^* g$ is uniformly $C^{l,\beta}$
equivalent
to $\hyp$ in M\"obius coordinates
by Lemma \eqref{properties-of-mobius-coordinates}, it follows that
the coefficients of the $j$th derivatives of $u$ that appear
in $P_iu$ are uniformly bounded in 
$C^{l-m+j,\beta}(B_2)$.  Therefore, using
Lemmas
\ref{lemma:norms-by-covering} and 
\ref{properties-of-spaces}\eqref{part:mult-cts},
we have 
\begin{align*}
\|Pu\|_{k,p,\delta}
&\le C\sum_i \rho(p_i)^{-\delta} \|\Phi_i^* Pu\|_{k,p;B_2}\\
&=  C\sum_i \rho(p_i)^{-\delta} \|P_i (\Phi_i^*u)\|_{k,p;B_2}\\
&\le  C'\sum_i \rho(p_i)^{-\delta} \|\Phi_i^*u\|_{k-m,p;B_2}\\
&\le C''\|u\|_{k-m,p,\delta},
\end{align*}
with an analogous estimate in the H\"older case.
\end{proof}

Recall from Chapter \ref{spaces-section} that when $p^* = p/(p-1)$, 
$H^{0,p^*}_{-\delta}(M;E)$ is naturally dual to $H^{0,p}_\delta(M;E)$
under the standard $L^2$ pairing.  The next lemma shows that
$P$ is symmetric with respect to this pairing.

\begin{lemma}\label{lemma:P-symmetric}
Suppose $P$ satisfies the hypotheses of Theorem \ref{thm:main-fredholm}.
If $1<p<\infty$, $p^* = p/(p-1)$, and $\delta\in \R$, then
$(Pu,v) = (u,Pv)$ for all $u\in H^{0,p}_{\delta}(M;E)$, $v\in 
H^{0,p^*}_{-\delta}(M;E)$.
\end{lemma}

\begin{proof}
This is true for $u,v\in C^\infty_c(M;E)$ by the fact that
$P$ is formally self-adjoint.  The general result follows
by density.
\end{proof}

The following lemma is a standard application of rescaling techniques
and classical interior elliptic regularity (cf.\ \cite[Prop.\ 3.4]{Graham-Lee},
\cite[Lemma 2.4]{Andersson}, and \cite[Lemma 4.1.1]{AC}).

\begin{lemma}\label{rescaling-estimates}
Let $P\colon C^\infty(M;E)\to C^\infty(M;F)$ be a geometric elliptic
operator of order $m$.
\begin{enumerate}\letters
\item\label{rescaling-sobolev-case}
Suppose $\delta\in \R$, $1<p<\infty$, and $m\le k\le l$.
If $u \in H^{0,p}_{\delta}(M;E)$ is such that $Pu\in
H^{k-m,p}_{\delta}(M;F)$, then $u \in H^{k,p}_{\delta}(M;E)$ and
\begin{equation}\label{rescaling-sobolev-estimate}
\|u\|_{k,p,\delta} \le C(\|Pu\|_{k-m,p,\delta} + \|u\|_{0,p,\delta}).
\end{equation}
\item\label{rescaling-holder-case}
Suppose $\delta\in\R$, $0<\alpha<1$, and $m< k+\alpha\le l+\beta$.  If $u \in
C^{0,0}_{\delta}(M;E)$ is such that $Pu\in C^{k-m,\alpha}_{\delta}(M;F)$, then $u \in
C^{k,\alpha}_{\delta}(M;E)$ and
\begin{equation}\label{rescaling-holder-estimate}
\|u\|_{k,\alpha,\delta} \le C(\|Pu\|_{k-m,\alpha,\delta} + \|u\|_{0,0,\delta}).
\end{equation}
\end{enumerate}
\end{lemma}

\begin{proof}
Under the hypotheses of case \eqref{rescaling-sobolev-case}, $u$ is
locally in $H^{k,p}$ by interior elliptic regularity, so only the
estimate \eqref{rescaling-sobolev-estimate} needs to be proved.  
Let $\{\Phi_i\}$ be a uniformly locally finite covering
of $M$ by M\"obius charts as in 
Lemma \ref{covering-by-mobius-coordinates}, and let
$P_i = \Phi_i^* \circ P \circ
(\Phi_i^{-1})^*$ as in the proof of Lemma
\ref{mapping-properties}.
Since
the coefficients of the  highest-order terms in $P_i$ are
constant-coefficient polynomials in 
$g_i=\Phi_i^*g$ and $(\det g_i)^{-1/2}$,
and $g_i$ is 
uniformly equivalent to the 
hyperbolic metric by Lemma \ref{properties-of-mobius-coordinates},
$P_{i}$ is uniformly elliptic on $B_2$.
Moreover,
since the coefficients of $P_i$ are uniformly bounded
in $C^{k-m,\alpha}(B_2)$ by the same lemma,
we have the following
standard local elliptic
estimate \cite{Nirenberg,ADN}:
\begin{displaymath}
\|u\|_{k,p;B_{1}} \le C(\|P_{i}u\|_{k-m,p;B_2} +
\|u\|_{0,p;B_2}),
\end{displaymath}
where the constant $C$ depends on $P$, $k$, $p$, and $\delta$, but is
independent of $u$ and $i$.  
Thus, using  Lemma \ref{lemma:norms-by-covering} again,
we have
\begin{align*}
\|u\|_{k,p,\delta}
&\le C \sum_i \rho(p_i)^{-\delta} \|\Phi_i^*u\|_{k,p;B_1}\\
&\le C' \sum_i \rho(p_i)^{-\delta} (\|P_i(\Phi_i^*u)\|_{k-m,p;B_2}
+ \|\Phi_i^*u\|_{0,p;B_2})\\
&= C' \sum_i \rho(p_i)^{-\delta} \|\Phi_i^*(Pu)\|_{k-m,p;B_2}
+ C' \sum_i \rho(p_i)^{-\delta} \|\Phi_i^*u\|_{0,p;B_2}\\
&\le C'' (\|Pu\|_{k-m,p,\delta} + \|u\|_{0,p,\delta}).
\end{align*}

The argument for case \eqref{rescaling-holder-case} is similar, using
interior H\"older estimates \cite{DN}.
\end{proof}

\begin{lemma}
If $P$ satisfies the hypotheses of Theorem \ref{thm:main-fredholm},
then $P$ is self-adjoint as an unbounded operator on $L^2(M;E)$.
\end{lemma}

\begin{proof}
Because of the density of $C^\infty_c(M;E)$ in $H^{m,2}(M;E)$,
$P$ is densely defined, and  Lemma \ref{rescaling-estimates}
shows that its domain is exactly $H^{m,2}(M;E)$.  
Clearly the domain of its Hilbert space 
adjoint contains $H^{m,2}(M;E)$.
On the other hand, if $v$ is in the domain of the adjoint,
then there exists $w\in L^2(M;E)$ such that
$(v,Pu) = (w,u)$ for all $u\in L^2(M;E)$.  This means
in particular that $Pv = w$ as distributions, which
by Lemma \ref{rescaling-estimates} implies that 
$v\in H^{m,2}(M;E)$.  Thus the domain of the adjoint is
equal to the domain of $P$.
\end{proof}

Next we present some elementary preliminary results
about Fredholm
operators on these weighted spaces.
Our first
lemma reduces the problem of proving $L^p$
Fredholm theorems to 
estimates near the boundary.  
A partial differential
operator $P$ acting on sections of a vector bundle
over a connected manifold $M$
is said to have the \defname{weak
unique continuation property} if any solution to $Pu=0$
that vanishes on an open set must vanish on all of $M$.
It has been shown recently by 
G. Nakamura, G. Uhlmann, and J.-N. Wang \cite{NUW}
that every 
strongly elliptic operator has this property.
We say an operator is {\it semi-Fredholm} if it
has finite-dimensional kernel and closed range.

\begin{lemma}\label{lemma:semi-fred}\label{lemma:fred}
Let $P\colon C^\infty(M;E)\to C^\infty(M;E)$ be a
formally self-adjoint geometric elliptic
operator of order $m\le l$.
Suppose $1<p<\infty$
and $\delta\in\R$.
Then 
$P\colon H^{k,p}_\delta(M;E)\to H^{k-m,p}_\delta(M;E)$
is semi-Fredholm for $m\le k\le l$
if and only if 
there exist 
a compact subset $K\subset M$ and a constant $c>0$
such that the following estimate holds for all
$u\in C^\infty_c(M\setminus K;E)$:
\begin{equation}\label{eq:p-estimate}
c\|u\|_{0,p,\delta}
\le \|Pu\|_{0,p,\delta}.
\end{equation}
If $P$ is semi-Fredholm and has the weak unique continuation property,
then \eqref{eq:p-estimate}
holds for every compact subset $K\subset M$
with nonempty interior.
If both \eqref{eq:p-estimate} and
\begin{equation}\label{eq:p*-estimate}
c\|u\|_{0,p^*,-\delta}
\le \|Pu\|_{0,p^*,-\delta}
\end{equation}
hold for $u\in C^\infty_c(M\setminus K;E)$,
where $p^* = p/(p-1)$,
then $P$ is Fredholm.
\end{lemma}

\begin{proof}
The argument in the 
forward direction
is fairly standard;
see, for example, \cite[Thm.\ 1.10]{Bartnik} and
\cite[Prop.\ 2.6]{Andersson}.
Let $U$ and $V$ be precompact open subsets of $M$ 
such that $K\subset U\subset \overline U \subset V$, and
let $\psi\in C^\infty(M)$ be a smooth
bump function that is 
equal to $1$ on $K$
and supported in $\overline U$.
It is easy to check that 
multiplication by $\psi$ is a bounded map from
$H^{k,p}_{\delta}(M;E)$ to itself for any $k,p,\delta$.
On the compact set $\overline V$,
the $H^{k,p}_{\delta}$ norm is uniformly equivalent 
to the standard $H^{k,p}$
norm, and $P$ is uniformly elliptic.  

First observe that 
\eqref{eq:p-estimate} 
extends to all $u\in H^{k,p}_{\delta}(M;E)$
by continuity.
For any $u\in H^{k,p}_{\delta}(M;E)$, we can write $u = u_\infty + u_0$, where
$u_0 := \psi u$ is supported in $\overline U$ and $u_\infty := (1-\psi)u$
is supported in $M\setminus K$.
We estimate as follows:
\begin{displaymath}
\|u\|_{k,p,\delta}\le \|u_\infty\|_{k,p,\delta} + \|u_0\|_{k,p,\delta}.
\end{displaymath}
For the first term, \eqref{rescaling-sobolev-estimate} gives
\begin{align*}
\|u_\infty\|_{k,p,\delta}
&\le C_1(\|Pu_\infty\|_{k-m,p,\delta}+\|u_\infty\|_{0,p,\delta})\\
&\le C_2(\|Pu_\infty\|_{k-m,p,\delta}+\|Pu_\infty\|_{0,p,\delta})\\
&\le C_2(\|(1-\psi)Pu\|_{k-m,p,\delta} +\|[P,\psi]u\|_{k-m,p,\delta})\\
&\le C_3(\|Pu\|_{k-m,p,\delta} +\|u\|_{k-1,p;\overline V}),
\end{align*}
where in the last line we have used the fact that
$[P,\psi]$ is an operator of order $m-1$
supported on $\overline V$.

For the second term,
standard interior elliptic
estimates yield
\begin{align*}
\|u_0\|_{k,p,\delta}
&= \|u_0\|_{k,p,\delta;\overline U}\\
&\le C_1(\|Pu_0\|_{k-m,p,\delta;\overline V} + \|u_0\|_{0,p,\delta;\overline V})\\
&\le C_1(\|\psi Pu\|_{k-m,p,\delta} +
\|[P,\psi]u\|_{k-m,p,\delta} + \|u\|_{0,p,\delta;\overline V})\\
&\le C_2(\|Pu\|_{k-m,p,\delta} +  \|u\|_{k-1,p;\overline V}).
\end{align*}
Finally, standard interpolation inequalities on the compact set $\overline V$ (see
\cite[Thm.\ 3.70]{Aubin}) allow us to replace $\|u\|_{k-1,p;\overline V}$
by
$C\|u\|_{0,p;\overline V} + \epsilon \|u\|_{k,p,\delta}$ with an arbitrarily small
constant $\epsilon$.  Absorbing the $\epsilon$ term on the left-hand side,
we conclude
\begin{equation}\label{compact-estimate}
\|u\|_{k,p,\delta}\le C(\|Pu\|_{k-m,p,\delta} + \|u\|_{0,p;\overline V} ).
\end{equation}

Now suppose $\{u_i\}$ is a sequence in $\Ker P \intersect
H^{k,p}_{\delta}(M;E)$.  Normalize $u_i$ so that $\|u_i\|_{k,p,\delta} = 1$.  
By the Rellich lemma on the compact set $\overline V$,
some subsequence converges
in $H^{0,p}(\overline V;E)$.  
From \eqref{compact-estimate} we conclude that this
subsequence is Cauchy and hence convergent in $H^{k,p}_{\delta}(M;E)$, and
therefore $\Ker P\intersect H^{k,p}_{\delta}(M;E)$ is finite-dimensional.

Next we need to show that $P$ has closed range.  Since $\Ker P$ is
finite-dimensional, it has a closed complementary subspace $Y\subset
H^{k,p}_{\delta}(M;E)$.  I claim there is a constant $C$ such that
\begin{equation}\label{estimate-on-Y}
\|u\|_{k,p,\delta} \le C\|Pu\|_{k-m,p,\delta}\qquad \text{ for all $u \in Y$.}
\end{equation}
If not, there is a sequence $\{u_i\}\subset Y$ with
\begin{displaymath}
\|u_i\|_{k,p,\delta} = 1 \text{ and } \|Pu_i\|_{k-m,p,\delta} \to 0.
\end{displaymath}
By the Rellich lemma again,
there exists a subsequence
(still denoted by $\{u_i\}$) that converges in $H^{0,p}(\overline V;E)$.  Then
\eqref{compact-estimate} shows that the subsequence also converges in
$H^{k,p}_{\delta}(M;E)$.  However, the limit $u$ then satisfies
$\|u\|_{k,p,\delta}=1$ and $Pu=0$ by continuity, and is therefore a nonzero
element of $Y \intersect \Ker P$, which is a contradiction.

Now if $u_i \in H^{k,p}_{\delta}(M;E)$ with $Pu_i = f_i \to f$ in
$H^{k-m,p}_{\delta}(M;F)$, we can assume without loss of generality that each
$u_i\in Y$, and then \eqref{estimate-on-Y} shows that $\{u_i\}$
converges in $H^{k,p}_{\delta}(M;E)$, which shows that $P$ has closed range.

Conversely, suppose $P\colon H^{k,p}_{\delta}(M;E)\to
H^{k-m,p}_{\delta}(M;E)$ is semi-Fredholm.  
Let $K\subset M$ be a compact subset
chosen as follows:  If $P$ has the unique
continuation property, 
$K$ can be any compact subset of $M$ with nonempty interior;
otherwise, let $K = M\setminus A_\epsilon$, where 
$\epsilon>0$ is chosen small enough that
no element of the finite-dimensional
space $\Ker P\intersect H^{m,p}_{\delta}(M;E)$ vanishes identically on
$K$.  

The key fact is that
there exists $c>0$ such that
\begin{equation}\label{distance-to-kernel}
\|u-v\|_{0,p,\delta} \ge  c\|u\|_{0,p,\delta}
\end{equation}
whenever $u\in H^{0,p}_{\delta}(M;E)$ is supported in $M\setminus K$ and $v\in
\Ker P\intersect H^{m,p}_{\delta}(M;E)$.  
To see this, note that our choice
of $K$ ensures that $\|\cdot\|_{0,p,\delta;K}$ is a norm
on $\Ker P\intersect H^{m,p}_{\delta}(M;E)$.  
Since all norms on a finite-dimensional vector space are
equivalent, there exists a constant $a$ such that
\begin{displaymath}
\|v\|_{0,p,\delta;K} \ge a \|v\|_{0,p,\delta}
\end{displaymath}
for all $v\in \Ker P\intersect H^{m,p}_{\delta}(M;E)$.  
Suppose now that $u$ is supported in
$M\setminus K$, $\|u\|_{0,p,\delta}=1$, and 
$v\in\Ker P\intersect H^{m,p}_{\delta}(M;E)$.  If $\|v\|_{0,p,\delta}\le
1/(1+a)$, then by the reverse triangle inequality
\begin{displaymath}
\|u-v\|_{0,p,\delta}\ge \|u\|_{0,p,\delta} - \|v\|_{0,p,\delta} \ge
1 - \frac{1}{1+a} = \frac{a}{1+a}.
\end{displaymath}
If on the other hand $\|v\|_{0,p,\delta} \ge 1/(1+a)$, then because $u$
vanishes on $K$,
\begin{align*}
\|u-v\|_{0,p,\delta} &\ge \|u-v\|_{0,p,\delta;K} = \|v\|_{0,p,\delta;K}\\
& \ge a\|v\|_{0,p,\delta}\ge \frac{a}{1+a}.
\end{align*}
Inequality \eqref{distance-to-kernel} (with $c=a/(1+a)$)
then follows for general
$u$ by homogeneity.

As above, $\Ker P\cap H^{k,p}_{\delta}(M;E)$ has a closed complementary
subspace $Y$, and both $Y$ and $P(H^{k,p}_{\delta}(M;E))\subset
H^{k-m,p}_{\delta}(M;E)$ are Banach spaces with the induced norms.  Then $P|_Y\colon Y
\to P(H^{k,p}_{\delta}(M;E))$ is bijective, and its inverse $(P|_Y)^{-1}\colon 
P(H^{k,p}_{\delta}(M;E)) \to Y$ is bounded
by the open mapping theorem.  This means there exists a constant $C > 0$
such that
\eqref{estimate-on-Y} holds.

Let $u\in H^{k,p}_{\delta}(M;E)$ be supported in $M\setminus K$, and write $u =
u_0 + u_Y$, with $u_0 \in \Ker P$, $u_Y \in Y$.  It follows from
\eqref{distance-to-kernel} that
\begin{displaymath}
\|u_Y\|_{0,p,\delta} = \|u-u_0\|_{0,p,\delta} \ge c\|u\|_{0,p,\delta},
\end{displaymath}
and therefore, by \eqref{estimate-on-Y} with $k=m$,
\begin{align*}
\|u\|_{0,p,\delta} &\le
c^{-1}\|u_Y\|_{0,p,\delta} \le c^{-1}\|u_Y\|_{m,p,\delta}\\
&\le c^{-1}C\|Pu_Y\|_{0,p,\delta} = c^{-1}C\|Pu\|_{0,p,\delta},
\end{align*}
which is \eqref{eq:p-estimate}.

Finally, suppose that both \eqref{eq:p-estimate} and \eqref{eq:p*-estimate}
hold.
To show that $P$ is actually Fredholm, 
all that remains to be shown is that 
the range of 
$P$ has finite codimension in $H^{k-m,p}_{\delta}(M;F)$.
Recall that $H^{0,p^*}_{-\delta }(M;E)$ is
dual to $H^{0,p}_{\delta}(M;E)$ under the standard $L^2$ pairing.
The argument above, using \eqref{eq:p*-estimate}
instead of \eqref{eq:p-estimate}, shows that
$P^*=P\colon H^{k,p^*}_{-\delta}(M;E)\to H^{k-m,p^*}_{-\delta}(M;E)$
has finite-dimensional kernel.
Any $v\in H^{0,p^*}_{-\delta}(M;E)$ that annihilates
the range of $P$ in $H^{k-m,p}_{\delta}(M;E)$
satisfies in particular $(v,Pu)=0$ for all $u\in C^\infty_c(M;E)$,
so is a distribution solution to $Pv=0$.
By Lemma 
\ref{rescaling-estimates},  $v\in H^{k,p^*}_{-\delta }(M;E)$.
Thus there is at most a finite-dimensional subspace of 
$H^{0,p^*}_{-\delta}(M;E)$ that annihilates
the range of $P$.  Since $P(H^{k,p}_{\delta}(M;E))$ is closed in
$H^{k-m,p}_{\delta}(M;E)$, it has finite codimension.
\end{proof}

%%
%%
%% Fredholm operators and Einstein metrics
%% on conformally compact manifolds
%%
%% by John M. Lee
%% 
%% Chapter 5
%%

\chapter{Analysis on Hyperbolic Space}\label{model-section}

In this chapter, we will analyze the behavior of 
geometric elliptic operators on hyperbolic space,
which serves as a model for the more general case.
Our goal is to show that if $P$ is an operator on
hyperbolic space satisfying the hypotheses of
Theorem \ref{thm:main-fredholm}, then $P$ is
an isomorphism on appropriate weighted Sobolev and H\"older
spaces.  (In the next chapter, we will
use the resulting inverse map to piece together a
parametrix for the analogous operator acting on an
arbitrary asymptotically hyperbolic manifold.)

For the purposes of this chapter,
we will use the Poincar\'e ball model,
identifying hyperbolic space with the unit ball
$\B\subset\R^{n+1}$, with coordinates
$(\xi^1,\dots,\xi^{n+1})$, and with the hyperbolic 
metric $\hyp=4(1-|\xi|)^{-2}\sum_i(d\xi^i)^2$.
The hyperbolic distance function
can be written in terms of the Euclidean norm and dot product as
\begin{displaymath}
d_{\hyp}(\xi,\eta) = \cosh^{-1} \frac{ (1+|\xi|^2)(1+|\eta|^2) - 4 \xi\centerdot \eta}
{ (1-|\xi|^2)(1-|\eta|^2) }.
\end{displaymath}
It will be convenient to use 
\begin{displaymath}
\rho(\xi) =\frac{1}{ \cosh d_{\hyp}(\xi,0)}= \frac{1-|\xi|^2 }{ 1+|\xi|^2 } 
\end{displaymath}
as a defining function for the ball, where
$0=(0,\dots,0)$
denotes the origin in
$\B \subset\R^{n+1}$.

Throughout this chapter, 
$E$ will be a geometric tensor
bundle of weight $r$ over $\B$, 
and $P\colon C^\infty(\B;E)\to C^\infty(\B;E)$
will be a formally self-adjoint geometric elliptic operator of order $m$.
The fact that $P$ is geometric implies that it is isometry
invariant:  
If $\phi$ is any orientation-preserving hyperbolic isometry
and $u$ is any section of $E$, then
\begin{equation}\label{eq:isometry-invariance}
\phi^*(Pu)= P(\phi^*u).
\end{equation}
We will assume that $P$ satisfies
\eqref{eq:asymptotic-L2-estimate}.
Then by Lemma \ref{lemma:semi-fred}, 
$P\colon H^{m,2}(\B ;E)\to 
L^2(\B ;E)$ is Fredholm.  The next
lemma shows that this is equivalent to 
being an isomorphism.

\begin{lemma}\label{fred=iso-on-model}
Suppose $P\colon C^\infty(\B ;E)\to C^\infty(\B ;E)$
is a geometric elliptic 
operator of order $m$ on $\B $.  Then
$P\colon H^{m,2}(\B ;E)\to H^{0,2}(\B ;F)$
is Fredholm if and only if it is an isomorphism.
\end{lemma}

\begin{proof}
If $P$ is an isomorphism, then clearly it is Fredholm.  
Conversely, if $P$ is
Fredholm, then
by Lemma \ref{lemma:semi-fred} $P$ satisfies
\begin{equation}\label{hyperbolic-sobolev-estimate}
\|u\| \le C\|Pu\|
\end{equation}
whenever $u$ is supported in
the complement of some compact set $K$.
Suppose $u$ is {\it any} smooth,
compactly supported section of $E$.  There is a M\"obius transformation
$\phi$ such that $\phi^{-1}(\supp u) \subset \B\setminus K$,
so $\phi^*u$ satisfies 
\eqref{hyperbolic-sobolev-estimate}.
Because $P$ and the $L^2$ norm are preserved by $\phi$,
$u$ itself satisfies the same estimate. 
Therefore, by continuity,
\eqref{hyperbolic-sobolev-estimate} holds for all $u\in H^{m,2}$, so
$\Ker P$ is trivial.  Since $P$ is self-adjoint as an
unbounded operator on $L^2(M;E)$, its index is zero,
which means that it is also surjective.
\end{proof}

Let $K$ be the Green kernel of $P$: That is, $K$ is the
Schwartz kernel of the operator $P^{-1}\colon
L^2(\B ;E)\to H^{m,2}(\B ;E)$.  Invariantly,
$K$ is interpreted as a distributional section of the bundle
$\Hom(\pi_2^*E,\pi_1^*E)$ over $\B \cross \B $, 
where $\pi_j$ is projection
on the $j$th factor.  
For all $f\in L^2(\B;E)$,
\begin{equation}\label{convolution-equation}
P^{-1}f(\xi) = \int_\B K(\xi,\eta) f(\eta) dV_{\hyp}(\eta).
\end{equation}
Equivalently, if we write $K_\eta(\xi) =
K(\xi,\eta)$, $K_\eta$ satisfies
\begin{displaymath}
PK_\eta(\xi) = \delta_\eta(\xi) \Id_{E_\eta}
\end{displaymath}
in the distribution sense and $K_\eta\in L^2$
on the complement of a neighborhood of $\eta$.  Somewhat more explicitly,
for any $\eta\in \B $ and $u_0\in E_\eta$, $K_\eta u_0$ is a (distributional)
section of $E$ satisfying $P(K_\eta u_0) =\delta_\eta u_0$.  
In particular, 
$K_0(\xi) = K(\xi,0)$ can be viewed as a fundamental solution for
$P$ on $\B$ with pole at $0$.

By local elliptic
regularity, $K$ is $C^\infty$ away from the diagonal $\{\xi=\eta\}$.  Since
$P$ is formally self-adjoint, it is easy to check that $K$ satisfies the
symmetry condition $K(\eta,\xi) = K(\xi,\eta)^*$
for $\xi\ne \eta$, where $K(\xi,\eta)^*\colon E_\xi\to
E_\eta$ is the (pointwise) adjoint of $K(\xi,\eta)
\in\Hom(E_\eta,E_\xi)$.

We extend our defining function $\rho$ to 
a function $\rho\colon \overline\B \cross \overline\B \to [0,1]$
of two variables (still denoted by the same symbol) by 
\begin{displaymath}
\rho(\xi,\eta) = \frac{1}{\cosh d_{\hyp}(\xi,\eta)} = 
\frac{ (1-|\xi|^2)(1-|\eta|^2) }{ (1+|\xi|^2)(1+|\eta|^2) - 4 \xi\centerdot \eta}.
\end{displaymath}
Observe that $\rho(\xi,0) = \rho(\xi)$.

Our main technical tool in this chapter is the following 
decay estimate
for $K$.  

\begin{proposition}\label{prop:K-decay}
Let $P\colon C^\infty(\B ;E)\to C^\infty(\B ;E)$ be
a formally self-adjoint 
geometric elliptic operator of order $m$
satisfying 
\eqref{eq:asymptotic-L2-estimate}.
Then $P$ has positive indicial radius $R$, 
and for any $\epsilon>0$
there is a constant $C$ such that
\begin{equation}\label{eq:K-decay}
|K(\xi,\eta)| \le C\rho(\xi,\eta)^{n/2+R-\epsilon}
\end{equation}
whenever $d_{\hyp}(\xi,\eta)\ge 1$.  \(The norm here is the pointwise
operator norm on $\Hom(E_\eta,E_\xi)$ with respect
to the hyperbolic metric.\)
\end{proposition}

\begin{proof}
The isometry invariance of $P$
implies that $K$ has the following equivariance
property for any orientation-preserving 
hyperbolic isometry $\phi$:
\begin{equation}\label{eq:K-equivariance}
K(\phi(\xi),\phi(\eta)) = 
(\phi^*)^{-1}\circ
K(\xi,\eta)\circ\phi^*.
\end{equation}
Also, since $\rho$ is defined purely in terms of the hyperbolic distance
function, it is clearly isometry invariant:
\begin{displaymath}
\rho(\phi(\xi),\phi(\eta)) = \rho(\xi,\eta).
\end{displaymath}
Therefore, to prove \eqref{eq:K-decay},
it suffices to show that 
\begin{displaymath}
|K(\xi,0)| \le C \rho(\xi,0)^{n/2+R-\epsilon} = C \rho(\xi)^{n/2+R-\epsilon}.
\end{displaymath}
Note that $K_0(\xi) = K(\xi,0)$ defines a smooth section of the
bundle $\Hom(E_0,E)$ over $\B\setminus\{0\}$, 
whose fiber at $\xi$ is the vector
space $\Hom(E_0,E_\xi)$.

The group of orientation-preserving
hyperbolic isometries that fix $0$ is
exactly $\SO(n+1)$, acting linearly on the unit ball as isometries
of both the hyperbolic and Euclidean metrics.
Let $L\subset \B$ be the ray
$L=\{(0,\dots,0,t):0\le t <1\}$ along the $\xi^{n+1}$-axis.
The subgroup of $\SO(n+1)$ that fixes 
$L$ pointwise is $\SO(n)\subset \SO(n+1)$, 
realized as the group of linear isometries
acting in the first $n$ variables only.  
Observe that for each
$\xi_0\in L$, $\SO(n)$ 
acts orthogonally (or unitarily if $E$ is a complex
tensor bundle)
on the fiber $E_{\xi_0}$ by pulling back.

Let $E_0 = E_0^{(1)}\oplus\dots\oplus E_0^{(k)}$ be an
orthogonal decomposition
of $E_0$ into irreducible $\SO(n)$-invariant subspaces.  We extend
this to a decomposition of the bundle $E$ over $\B\setminus\{0\}$
as follows.  First,  
for each ${\xi_0}\in L$, let
$E_{\xi_0}^{(i)}$ be the subspace of $E_{\xi_0}$ obtained by parallel translating
$E_0^{(i)}$ along $L$ with respect
to the Euclidean metric $\overline g$; since $\SO(n)$ acts as
Euclidean isometries, it follows that 
\begin{displaymath}
E_{\xi_0} = E_{\xi_0}^{(1)}\oplus
\dots\oplus E_{\xi_0}^{(k)}
\end{displaymath}
is an orthogonal irreducible $\SO(n)$-decomposition of $E_{\xi_0}$.
Then for an arbitrary point $\xi\in \B \setminus\{0\}$, let 
$E_\xi^{(i)} = (\phi^*)^{-1}E_{\xi_0}^{(i)}$,
where $\xi_0$ is the 
unique point of $L$ such that $|\xi_0|=|\xi|$, and $\phi\in \SO(n+1)$
satisfies $\phi(\xi_0) = \xi$.  Since $E_{\xi_0}^{(i)}$ is invariant under 
$\SO(n)$, $E_\xi^{(i)}$ does not depend on the choice of $\phi$.
Since $\phi$ can be chosen locally to depend smoothly on $\xi$
(by means of a smooth local section of the submersion
$\SO(n+1)\to \SO(n+1)/\SO(n) = \S^{n}$), this results
in $k$ smooth subbundles $E^{(1)},\dots,E^{(k)}$ of $E$
over $\B \setminus\{0\}$.

For each pair of indices
$i,j=1,\dots,k$, we choose an  $\SO(n)$-equivariant
linear map $k_0^{(i,j)}\colon 
E_0\to E_{0}$ as follows:  If 
$E_0^{(i)}$ and $E_{0}^{(j)}$ are 
isomorphic as representations of $\SO(n)$, 
let $k_0^{(i,j)}$ be an $\SO(n)$-equivariant Euclidean
isometry
from $E_0^{(i)}$ to $E_0^{(j)}$, extended to be zero on $E_0^{(l)}$
for $l\ne i$;  
and otherwise
let  $k_0^{(i,j)}$ be the zero map.
By Schur's lemma, the nonzero maps $k_0^{(i,j)}$
form a basis for the space of $\SO(n)$-equivariant
endomorphisms of $E_0$.  Let us renumber these
nonzero maps as $k_0^1,\dots,k_0^N$.

Next we extend each map $k_0^{j}$ to a section $k^{j}$ of 
$\Hom(E_0,E)$ over $\B\setminus \{0\}$ in the same way as we extended the
spaces $E_0^{(i)}$: First,
for each point $\xi_0\in L$,
define $k_{\xi_0}^{j}\colon E_0\to E_{\xi_0}$
to be $k_0^{j}$ followed by $\overline g$-parallel translation
along $L$ from $0$ to $\xi_0$; and then for arbitrary $\xi\in \B\setminus \{0\}$,
define $k_\xi^{j} = (\phi^*)^{-1} \circ k_{\xi_0}^{j}\circ \phi^*$,
where $\xi_0\in L$ and 
$\phi\in\SO(n+1)$ satisfies
$\phi(\xi_0)=\xi$.

I claim that there are smooth functions $f_{1},\dots,f_N\colon (0,1)\to \C$
such that 
\begin{equation}\label{eq:K-decomp}
K_0(\xi) = \sum_{j} f_{j}(\rho(\xi)) k_\xi^{j}
\end{equation}
for all $\xi\in \B\setminus\{0\}$.
To see this, first note that for any point $\xi_0\in L$,
the 
equivariance property 
\eqref{eq:K-equivariance} implies that
$K_0(\xi_0)$ is an $\SO(n)$-equivariant linear map
from $E_0$ to $E_{\xi_0}$, and therefore by Schur's lemma
it can
be written as a linear combination of the maps
$k^{j}_{\xi_0}$:
\begin{displaymath}
K_0(\xi_0) = \sum_{j} c_{j}(\xi_0) k_{\xi_0}^{j}.
\end{displaymath}
For any other point $\xi\in \B\setminus\{0\}$, 
let $\xi_0$ be the point of $L$ such that
$|\xi_0|=|\xi|$, and let $\phi\in \SO(n+1)$
satisfy $\phi(\xi_0)=\xi$.
Then \eqref{eq:K-equivariance} yields
\begin{equation}\label{eq:K-decomp-proof}
\begin{aligned}
K_0(\xi) 
&= K(\phi(\xi_0),\phi(0))\\
&= (\phi^*)^{-1}\circ  K(\xi_0,0) \circ \phi^*\\
&= \sum_{j} c_{j}(\xi_0) (\phi^*)^{-1}\circ  k^{j}_{\xi_0} \circ \phi^*\\
&= \sum_{j} c_{j}(\xi_0) k^{j}_{\xi}. 
\end{aligned}
\end{equation}
Since $\rho\colon L\setminus\{0\}\to (0,1)$ is a diffeomorphism,
there are 
functions $f_j\colon(0,1)\to \C$ such that
$f_j(\rho(\xi))= f_j(\rho(\xi_0)) =  c_j(\xi_0) $ whenever $|\xi|=|\xi_0|$,
and then \eqref{eq:K-decomp-proof} is equivalent to \eqref{eq:K-decomp}.
The smoothness of $K_0$ implies that the functions
$f_{j}$ are smooth.

Now the equation $PK_0=0$ reduces to an analytic system of 
ordinary differential equations for the functions
$f_{j}$, and the fact that
$P$ is uniformly degenerate implies that this system
has a regular singular point at $\rho=0$.
Because the sections $k^{j}$ of $\Hom(E_0,E)$
extend smoothly to $\overline\B\setminus\{0\}$, our definition of
the indicial map of $P$ guarantees that the 
characteristic exponents of this system of ODEs are
precisely the characteristic exponents of the
operator $P$.  Therefore, by the standard theory of ODEs
with regular singular points, each coefficient function
$f_j$ satisfies
\begin{equation}\label{eq:fj-estimate}
|f_j(t)|\sim C_j t^{s_j} |\log t|^{k_j}\text{ as $t\to 0$}
\end{equation}
for some characteristic exponent $s_j$ and some nonnegative integer
$k_j$.

If $u_0$ is any tensor in one of the summands
$E_0^{(i)}$, then the images $k_\xi^j(u_0)$
lie in different summands of $E_\xi$ and are therefore
orthogonal, so 
\begin{displaymath}
|K_0(\xi)u_0|_{\overline g}^2 = \sum_j |f_j(\rho(\xi))|^2\, |k_\xi^ju_0|_{\overline g}^2,
\end{displaymath}
where $\overline g$ is the Euclidean metric on $\B$.
Since 
$k_\xi^j$ is  a Euclidean isometry onto its image, 
$|k_\xi^ju_0|_{\overline g}$ is independent
of $\xi$.  Therefore, if 
$k_\xi^ju_0\ne 0$, then on $\B\setminus\{0\}$ we have
\begin{displaymath}
|K_0(\xi)u_0|_{\overline g}\ge C |f_j(\rho(\xi))|
\ge C \rho(\xi)^{\Real s_j} |\log \rho(\xi)|^{k_j}
\end{displaymath}
for some
positive constant $C$.  
Because $K_0(\xi)u_0$ is in $L^2$ away from $0$,
by Lemma \ref{lemma:being-in-Lp}
we must have $\Real s_j>n/2-r$ for each such $j$.
Since for each $j$ there is some $u_0$ such that
$k_\xi^ju_0\ne 0$, the same inequality holds for every $j$.
By definition of the indicial radius $R$,
this implies that in fact $R>0$ and $\Real s_j\ge n/2-r+R$.
Using \eqref{eq:fj-estimate}, we conclude that
for any $\epsilon>0$ there is a 
constant $C$ such that 
\begin{displaymath}
|f_j(t)|\le C t^{n/2-r+R-\epsilon}\text{ for $t$ away from $1$.}
\end{displaymath}
This in turn implies
\begin{align*}
|K_0(\xi)u_0|_{\hyp} &\le C \rho(\xi)^{r}|K_0(\xi)u_0|_{\overline g}\\
&\le C' \rho(\xi)^{r}\rho(\xi)^{n/2-r+R-\epsilon}|u_0|_{\overline g}\\
&= C'' \rho(\xi)^{n/2+R-\epsilon}|u_0|_{\hyp}
\end{align*}
whenever $d(\xi,0)\ge 1$,
which was to be proved.
\end{proof}

The next two lemmas give some estimates that will be
needed to use our decay estimate for proving mapping
properties of $P^{-1}$.

\begin{lemma}\label{lemma:hypergeometric}
For any real numbers $p,q,r$ such that 
$p+1>0$ and $r>q+1>0$,
there exists a constant $C$ depending only on $p,q,r$ 
such that the following estimate holds
for all $u\in [0,1)$:
\begin{displaymath}
\int_0^1 \frac{ t^{p}(1-t)^{q}}{(1-ut)^{r}}\,dt
\le C (1-u)^{q+1-r}.
\end{displaymath}
\end{lemma}

\begin{proof}
We use the following standard integral
representation for hypergeometric functions \cite[p.~59]{Erdelyi}:
\begin{displaymath}
F(\alpha,\beta,\gamma;z) =
\frac{\Gamma(\gamma)}{\Gamma(\beta)\Gamma(\gamma-\beta)}
\int_0^1 \frac{t^{\beta-1}(1-t)^{\gamma-\beta-1}}
{(1-tz)^\alpha}\,dt,
\end{displaymath}
which is valid if $\Real \gamma > \Real \beta > 0$ and $|z|<1$.
The hypergeometric function
$F(\alpha,\beta,\gamma;z)$ is 
analytic for $|z|<1$ and satisfies a second-order ODE that has a
regular singular point at $z=1$ 
with characteristic exponents $0$ and
$\gamma-\alpha-\beta$ \cite[p.\ 246]{BR}.  
As long as $\gamma-\alpha-\beta<0$, therefore, it satisfies
\begin{displaymath}
|F(\alpha,\beta,\gamma;u)|\le C (1-u)^{\gamma-\alpha-\beta} 
\quad \text{if $0\le u<1$}.
\end{displaymath}
Applying this with $\alpha=r$, $\beta=p+1$, and $\gamma = p+q+2$
proves the lemma.
\end{proof}

\begin{lemma}\label{distance-integral-lemma}
Suppose $a$ and $b$ are real numbers such that $a+b>n$ and $a>b$.  There
exists a constant $C$ depending only on $n,a,b$ such that the following
estimate holds for all $\xi,\zeta\in \B$:
\begin{displaymath}
\int_{\B} \rho(\xi,\eta)^a \rho(\eta,\zeta)^b \, dV_{\hyp}(\eta) \le C \rho(\xi,\zeta)^b.
\end{displaymath}
\end{lemma}

\begin{proof}
By an isometry, we can arrange that $\zeta=0$ and
$\xi=(0,\dots,0,r)$ is on the positive $\xi^{n+1}$-axis.  
Substituting $\zeta=0$ into the integral, we must estimate
\begin{equation}\label{integral-to-estimate}
I:=\int_{\B } \left(\frac {(1-|\xi|^2)(1-|\eta|^2) }{(1 + |\xi|^2)(1+|\eta|^2) -
4 \xi\centerdot \eta} \right)^a
\left(\frac {1-|\eta|^2}{1 + |\eta|^2} \right)^b\, dV_{\hyp}(\eta).
\end{equation}

Parametrize the ball by the map
$\Phi\colon(0,1)\cross (0,\pi) \cross \S^{n-1}\to \B $ given by
\begin{displaymath}
\Phi(s,\theta,\omega) = (s\omega^1\sin\theta,\dots,
s\omega^n\sin\theta,s\cos\theta),
\end{displaymath}
so that $s$ is the Euclidean 
distance from $0$ and $\theta$ is the angle from the
positive $\xi^{n+1}$-axis.
In these coordinates, the hyperbolic metric is
\begin{displaymath}
\hyp = \frac{4}{(1-s^2)^2}(ds^2 + s^2 d\theta^2 + s^2\sin^2\theta \,\gcirc),
\end{displaymath}
where $\gcirc$ represents the standard metric on $\S^{n-1}$.  The
hyperbolic volume element is therefore
\begin{displaymath}
dV_{\hyp} = \frac{2^{n+1} s^n \sin^{n-1}\theta }{(1-s^2)^{n+1}} \,ds\,d\theta\, dV_{\gcirc},
\end{displaymath}
where $dV_{\gcirc}$ is the volume element on $\S^{n-1}$.

In these coordinates, we have $|\xi|=r$, $|\eta|=s$ and $\xi\centerdot \eta =
rs\cos\theta$.  The integrand in
\eqref{integral-to-estimate}
is constant on each sphere $\S^{n-1}$, so
we can immediately   integrate over $\S^{n-1}$ and write $I$
as a constant multiple of
\begin{displaymath}
\int_0^1\int_0^\pi
\left(\frac {(1-r^2)(1-s^2) }{(1 + r^2)(1+s^2) -
4rs\cos\theta} \right)^a
\left(\frac {1-s^2}{1 + s^2} \right)^b\frac{ s^n\sin^{n-1}\theta }{(1-s^2)^{n+1}} d\theta\, ds.
\end{displaymath}

Since we are only interested in estimates up to a constant multiple, we
will write $f\sim g$ to mean that $f/g$ is bounded above and below by 
positive constants
depending only on $a$, $b$, and $n$.  Thus, for example, $1-s^2
= (1-s)(1+s) \sim 1-s$, $1-r^2\sim 1-r$, $1+r^2\sim 1 + s^2\sim 1$,
and 
\begin{equation}\label{new-integral-to-estimate}
I\sim \int_0^1\int_0^\pi
\left(\frac{(1-r)(1-s) }{(1 + r^2)(1+s^2) -
4rs\cos\theta} \right)^a
(1-s)^b\frac{ s^n\sin^{n-1}\theta }{(1-s)^{n+1}} d\theta\, ds.
\end{equation}

The $\theta$ integral can be simplified by the substitution
$\cos\theta=2t-1$
to obtain
\begin{equation}\label{theta-integral}
\begin{aligned}
\int_0^\pi
&\frac { \sin^{n-1}\theta}{((1 + r^2)(1+s^2) -
4rs\cos\theta)^a}\,d\theta\\
&\qquad =
2^{n-1}B(r,s)^{a}
\int_{0}^{1} \frac{t^{n/2-1} (1-t)^{n/2-1} }
{(1-8 B(r,s) r s t)^a}\, dt,
\end{aligned}
\end{equation}
where
\begin{displaymath}
B(r,s) = \frac{1}{(1+r^2)(1+s^2)+4rs}\sim 1.
\end{displaymath}
Because our hypothesis guarantees that $n-a < b < a$
and therefore $a>n/2$, 
Lemma \ref{lemma:hypergeometric} shows that
the right-hand integral in 
\eqref{theta-integral}
is bounded by a constant multiple of $(1-8B(r,s) r s)^{n/2-a}$.
Substituting this into \eqref{new-integral-to-estimate}
yields
\begin{equation}\label{newer-integral-to-estimate}
I\le C (1-r)^a\int_0^1
\frac {s^n(1-s)^{a+b-n-1}}{(1 - 8 B(r,s) r s)^{a-n/2}}\,ds.
\end{equation}

A computation shows that
\begin{displaymath}
1 - 8 B(r,s) r s = \frac{(1-rs)^2 + (r-s)^2}{(1+rs)^2 + (r+s)^2}
\sim (1-rs)^2 + (r-s)^2\ge (1-rs)^2.
\end{displaymath}
Inserting this into \eqref{newer-integral-to-estimate},
we conclude that 
\begin{displaymath}
I\le C(1-r)^a\int_0^1
\frac {s^n(1-s)^{a+b-n-1}}{(1-rs)^{2a-n}}\,ds.
\end{displaymath}
Lemma \ref{lemma:hypergeometric}
then shows that this is bounded by
a multiple of $(1-r)^b\sim \rho(\xi,\zeta)^b$.
\end{proof}

The following estimate is the key to proving sharp mapping
properties of $P^{-1}$.

\begin{lemma}\label{lemma:K-integral}
Let $P$ satisfy the hypotheses of Proposition \ref{prop:K-decay}.
Then for any real number $b$ satisfying
$n/2-R<b<n/2+R$, there exists a constant $C$ such that
\begin{align*}
\int_\B |K(\xi,\eta)|\, \rho(\eta)^b\, dV_{\hyp}(\eta)&\le C\rho(\xi)^b,\\
\int_\B |K(\xi,\eta)|\, \rho(\xi)^b\, dV_{\hyp}(\xi)&\le C\rho(\eta)^b.
\end{align*}
\end{lemma}

\begin{proof}
Since $|K(\xi,\eta)| = |K(\eta,\xi)^*| = |K(\eta,\xi)|$ by self-adjointness,
the two inequalities are equivalent, so it suffices to prove the second one.
We will write
\begin{multline*}
\int_\B |K(\xi,\eta)| \, \rho(\xi)^b\, dV_{\hyp}(\xi)\\
= \int_{d_{\hyp}(\xi,\eta)\le 1}|K(\xi,\eta)| \, \rho(\xi)^b\, dV_{\hyp}(\xi)
+ \int_{d_{\hyp}(\xi,\eta)\ge 1}|K(\xi,\eta)| \, \rho(\xi)^b\, dV_{\hyp}(\xi)
\end{multline*}
and estimate each term separately.

For the first term, we observe
that $K$ is uniformly
locally integrable near $\xi=\eta$:
Since $K_0(\xi)=K(\xi,0)$ satisfies $PK_0 = \delta_0 \Id_{E_0}$,
and the Dirac delta function is in the Sobolev space
$H^{-1,q}$ for $1<q<1+1/n$ (defined as the dual space
to $H^{1,q^*}$, $q^* = q/(q-1)$),
local elliptic regularity implies that $K_0\in L^q_{\text{loc}}\subset
L^1_{\text{loc}}$.  Therefore, if $\eta\in \B$ is arbitrary and $\phi$ is any
M\"obius transformation sending $\eta$ to $0$, the 
change of variables $\xi'=\phi(\xi)$ yields 
\begin{equation}\label{eq:K-integrable}
\int_{d_{\hyp}(\xi,\eta)\le 1} |K(\xi,\eta)|\,dV_{\hyp}(\xi) 
= \int_{d_{\hyp}(\xi',0)\le 1} |K(\xi',0)|\,dV_{\hyp}(\xi') \le C.
\end{equation}

Using the triangle inequality together with the
elementary fact that $\cosh(A+B)\le 2 \cosh A \cosh B$ for 
$A,B\ge 0$, we estimate
\begin{displaymath}
\cosh d_{\hyp}(\xi,0)
\le \cosh (d_{\hyp}(\xi,\eta)+d_{\hyp}(\eta,0))
\le 2\cosh d_{\hyp}(\xi,\eta)\cosh d_{\hyp}(\eta,0).
\end{displaymath}
It follows that $\rho(\eta)\le 2\rho(\xi)$
on the set where $d_{\hyp}(\xi,\eta)\le 1$.
By symmetry, the same inequality holds with $\xi$ and $\eta$ reversed.
Thus
\begin{align*}
\int_{d_{\hyp}(\xi,\eta)\le 1}&|K(\xi,\eta)| \, \rho(\xi)^b\, dV_{\hyp}(\xi)\\
&\le \bigg(\sup_{d_{\hyp}(\xi,\eta)\le 1}\rho(\xi)^b\bigg)\int_{d_{\hyp}(\xi,\eta)\le 1}|K(\xi,\eta)| \, dV_{\hyp}(\xi)\\
&\le C\rho(\eta)^b.
\end{align*}

For the second term,
choose 
$\epsilon>0$ small enough that
\begin{equation}\label{eq:s-eps-est}
\frac n 2-R+\epsilon < b < \frac n 2+R-\epsilon.
\end{equation}
Then with $a = n/2 +R-\epsilon$, we have
$a+b>n$ and $a>b$, so we can use 
Proposition \ref{prop:K-decay}
and 
Lemma \ref{distance-integral-lemma}
to conclude
\begin{align*}
\int_{d_{\hyp}(\xi,\eta)\ge 1}|K(\xi,\eta)| \, \rho(\xi)^b\, dV_{\hyp}(\xi)
&\le \int_{d_{\hyp}(\xi,\eta)\ge 1}\rho(\xi,\eta)^a \, \rho(\xi)^b\, dV_{\hyp}(\xi)\\
&\le C\rho(\xi)^b.
\end{align*}
\end{proof}

\begin{proposition}\label{prop:K-maps-to-Lp}
If $1<p<\infty$, $k\ge m$, and $|\delta +n/p-n/2|<R$, then 
there exists a constant $C$ such that 
\begin{equation}\label{eq:hyperbolic-Lp-est}
\|u\|_{k,p,\delta} \le C \|Pu\|_{k-m,p,\delta}
\end{equation}
for all $u\in H^{k,p}_\delta(\B;E)$.
\end{proposition}

\begin{proof}
Using Lemma \ref{rescaling-estimates}, it suffices
to prove that
\begin{displaymath}
\|u\|_{0,p,\delta}\le C \|Pu\|_{0,p,\delta}
\end{displaymath}
for all $u\in H^{k,p}_{\delta}(\B;E)$.
Because $C^\infty_c(\B;E)$ is dense in 
$H^{k,p}_{\delta}(\B;E)$, it suffices to prove this
inequality for $u\in C^\infty_c(\B;E)$.
Since $u = P^{-1}(Pu)$ in that case, it suffices to prove
the estimate
\begin{equation}\label{eq:5.x}
\|P^{-1}f\|_{0,p,\delta}\le C \|f\|_{0,p,\delta}
\text{ for all $f\in C^\infty_c(\B;E)$}.
\end{equation}

Put 
\begin{align*}
p^* &= \frac{p}{p-1},\\
a &= \frac{1}{p^*} \left(\delta +\frac n p\right),
\end{align*}
so that 
\begin{equation}\label{eq:m-estimates}
\begin{aligned}
\frac n 2 - R &< ap^* < \frac n 2 + R,\\
\frac n 2 - R &< ap-\delta p < \frac n 2 + R.
\end{aligned}
\end{equation}

By H\"older's inequality and Lemma 
\ref{lemma:K-integral}, we estimate
\begin{align*}
|P^{-1}f(\xi)|_{\hyp} 
&\le \int_\B |K(\xi,\eta)| \, |f(\eta)|_{\hyp} \, dV_{\hyp}(\eta)\\
&= \int_\B \left(|K(\xi,\eta)|^{1/p}\rho(\eta)^{-a}|f(\eta)|_{\hyp}\right)
\left(|K(\xi,\eta)|^{1/p^*}\rho(\eta)^{a} \right)  dV_{\hyp}(\eta)\\
&\le \left(\int_\B |K(\xi,\eta)|\,\rho(\eta)^{-ap}|f(\eta)|_{\hyp}^p \,dV_{\hyp}(\eta)\right)^{1/p}
\cross\\
&\qquad \qquad
\left(\int_\B |K(\xi,\eta)|\,\rho(\eta)^{ap^*}\,dV_{\hyp}(\eta)\right)^{1/p^*}\\
&\le C\rho(\xi)^a
\left(\int_\B |K(\xi,\eta)|\,\rho(\eta)^{-ap}|f(\eta)|_{\hyp}^p\,dV_{\hyp}(\eta)\right)^{1/p}.
\end{align*}
Therefore,
\begin{align*}
\|P^{-1}f\|_{0,p,\delta}^p
&= \int_\B \rho(\xi)^{-\delta p} |P^{-1}f(\xi)|^p_{\hyp}\,dV_{\hyp}(\xi)\\
&\le C^p\int_\B\int_\B \rho(\xi)^{ap- \delta p} |K(\xi,\eta)|\, \rho(\eta)^{-ap} |f(\eta)|_{\hyp}^p\,
dV_{\hyp}(\eta)\, dV_{\hyp}(\xi).
\end{align*}
By  Lemma 
\ref{lemma:K-integral} again,
we can evaluate the $\xi$ integral first to obtain 
\begin{align*}
\|P^{-1}f\|_{0,p,\delta}^p
&\le C'\int_\B \rho(\eta)^{ap-\delta p} \rho(\eta)^{-ap} |f(\eta)|_{\hyp}^p\,
dV_{\hyp}(\eta)\\
&= C'\|f\|_{0,p,\delta}^p.
\end{align*}
\end{proof}

\begin{theorem}\label{thm:hyperbolic-sobolev-iso}
Let $P\colon C^\infty(\B ;E)\to C^\infty(\B ;E)$ be
a formally self-adjoint 
geometric elliptic operator of order $m$
satisfying 
\eqref{eq:asymptotic-L2-estimate}.
If $k\ge m$, $1<p<\infty$, and $|\delta +n/p-n/2|<R$, then
the natural extension
$P\colon H^{k,p}_{\delta}(\B;E)\to H^{k-m,p}_{\delta}(\B;E)$
is an isomorphism.
\end{theorem}

\begin{proof}
Injectivity is an immediate consequence of \eqref{eq:hyperbolic-Lp-est}.
To prove surjectivity, let $f\in H^{k-m,p}_{\delta}(\B;E)$
be arbitrary, and let $f_i\in C^\infty_c(\B;E)$ be a sequence 
such that $f_i\to f$ in $H^{k-m,p}_{\delta}(\B;E)$.
Set $u_i = P^{-1}f_i\in H^{m,2}(\B;E)$, so that $Pu_i = f_i$.
Then each $u_i$ is in $H^{0,p}_{\delta}(\B;E)$ by \eqref{eq:5.x}, 
and in $H^{k,p}_{\delta}(\B;E)$
by Lemma \ref{rescaling-estimates}, and 
\eqref{eq:hyperbolic-Lp-est} shows that
$\{u_i\}$ is Cauchy in $H^{k,p}_{\delta}(\B;E)$.
It follows that $u = \lim u_i\in H^{k,p}_{\delta}(\B;E)$
satisfies $Pu=f$ as desired, so $P$ is surjective.
The continuity of the inverse map then follows from 
\eqref{eq:hyperbolic-Lp-est}.
\end{proof}

Now we turn our attention to the H\"older case.
First we prove an estimate analogous to 
\eqref{eq:5.x}.

\begin{proposition}\label{prop:K-maps-to-Linfty}
If $|\delta -n/2|<R$, 
there exists a constant $C$ such that
\begin{equation}\label{eq:hyperbolic-sup-est}
\|P^{-1}f\|_{0,0,\delta}\le C \|f\|_{0,0,\delta}
\end{equation}
for all $f\in C^{0,0}_{\delta}(\B;E)$.
\end{proposition}

\begin{proof}
By Lemma \ref{lemma:K-integral},
\begin{align*}
|P^{-1}f(\xi)|_{\hyp}
&\le \int_\B |K(\xi,\eta)|\, |f(\eta)|_{\hyp} \, dV_{\hyp}(\eta)\\
&\le C\int_\B |K(\xi,\eta)|\,\rho(\eta)^\delta\, \|f\|_{0,0,\delta} \, dV_{\hyp}(\eta)\\
&\le C' \rho(\xi)^\delta  \|f\|_{0,0,\delta},
\end{align*}
which implies
\begin{align*}
\|P^{-1}f\|_{0,0,\delta} = \sup_{\xi\in \B}\big( \rho(\xi)^{-\delta}|P^{-1}f(\xi)|_{\hyp}\big)
\le C' \|f\|_{0,0,\delta}.
\end{align*}
\end{proof}

\begin{theorem}
Let $P\colon C^\infty(\B;E)\to C^\infty(\B;E)$ be
a formally self-adjoint 
geometric elliptic operator of order $m$
satisfying 
\eqref{eq:asymptotic-L2-estimate}.
If $0<\alpha<1$, $k\ge m$, and $|\delta -n/2|<R$, then
the natural extension
$P\colon C^{k,\alpha}_{\delta}(\B;E)\to C^{k-m,\alpha}_{\delta}(\B;E)$
is an isomorphism.
\end{theorem}

\begin{proof}
To prove surjectivity, let $f\in C^{k-m,\alpha}_{\delta}(\B;E)$
be arbitrary and set $u = P^{-1}f$, so that
$u\in C^{0,0}_{\delta}(\B;E)$ by 
Proposition 
\ref{prop:K-maps-to-Linfty}.
An easy computation shows that $Pu=f$ in the
distribution sense, so $u\in C^{k,\alpha}_{\delta}(\B;E)$
by Lemma \ref{rescaling-estimates}.

To prove injectivity,
choose 
$\delta'$ close to $\delta$ and $p$ large such that
$\delta>\delta'+n/p$ and $|\delta'+n/p-n/2|<R$.
Then $C^{k,\alpha}_{\delta}(\B;E)\subset
H^{k,p}_{\delta'}(\B;E)$ by 
by Lemma \ref{properties-of-spaces}.
Since $P$ is injective on $H^{k,p}_{\delta'}(\B;E)$ 
by Theorem \ref{thm:hyperbolic-sobolev-iso}, it
is injective on the smaller space
$C^{k,\alpha}_{\delta}(\B;E)$.
\end{proof}

%%
%%
%% Fredholm operators and Einstein metrics
%% on conformally compact manifolds
%%
%% by John M. Lee
%% 
%% Chapter 6
%%

\chapter{Fredholm Theorems}\label{section:fredholm}

In this chapter, we return to the general case of 
a connected $(n+1)$-manifold
$(M,g)$, assumed to be 
asymptotically hyperbolic 
of class $C^{l,\beta}$, with $l\ge 2$ and $0\le\beta<1$.
Let $E$ be a geometric tensor bundle over $M$, and let 
$P\colon C^\infty(M;E)\to C^\infty(M;E)$ be a
formally self-adjoint geometric elliptic operator of order $m\ge 1$.
We will prove Fredholm properties of $P$ by using 
M\"obius coordinates
near the boundary to piece together a parametrix modeled
on the inverse operator on hyperbolic space.

For this purpose, we will need
a slightly modified version of 
M\"obius coordinates.
Whereas the original M\"obius coordinates defined in Chapter 
\ref{mobius-section} were valid 
in a neighborhood of an
interior point,
for our parametrix construction we will need coordinates
that are 
defined all the way up to the boundary
and adjusted to make the background metric $\gbar$ 
close to  the
Euclidean metric on a neighborhood of a boundary point.
These coordinates will be used to transfer the operator
$\breve P^{-1}$ to
$M$ with an error that decays to one higher order along
the boundary.  

For each point $\hat p\in \del M$, choose some
neighborhood $\Omega$ on which background coordinates
$(\theta,\rho)$ 
are defined on a set of the form
\eqref{eq:bkg-chart-set}.
Let
$\omega^1,\dots,\omega^n\in C^{l,\beta}_{(0)}(\Omega,T^*\overline M)$
be $1$-forms 
chosen so that 
$(\omega^1,\dots,\omega^n,d\rho)$
is an
orthonormal coframe for  $\overline g$
at each point of $\del M\cap \Omega$
(recall that $|d\rho|_{\overline g}\equiv 1$ along $\del M$).
Let $A^\alpha_\beta$, $B^\alpha$ be the coefficients
of $\omega^\alpha$ at $\hat p$, defined by
\begin{displaymath}
\omega^\alpha_{\hat p} = A^\alpha_\beta d\theta^\beta_{\hat p} + 
B^\alpha d\rho_{\hat p},
\end{displaymath}
and let $(\tw\theta^1,\dots,\tw\theta^n)$ be the functions
defined on  $\Omega$ by
\begin{displaymath}
\tw\theta^\alpha = A^\alpha_\beta \theta^\beta + 
B^\alpha \rho.
\end{displaymath}
Then $(\tw\theta^1,\dots,\tw\theta^n,\rho)$ form 
coordinates on $\Omega$, and 
in these new coordinates $\overline g$ has the
matrix $\delta_{ij}$ at $\hat p$.  
For  $0<a$ and $0<r<c$, define open subsets
$Y_a\subset \H$ and $Z_r(\hat p)\subset\Omega\subset M$ by
\begin{align*}
Y_a &= \{(x,y)\in \H: |x|<a, 0<y<a\},\\
Z_r(\hat p) &= \{(\tw\theta,\rho)\in \Omega: |\tw\theta|<r, 0<\rho<r\}.
\end{align*}
For $0<r<c$, define a chart
$\Psi_{\hat p,r}\colon Y_1\to Z_r(\hat p)$
by 
\begin{equation*}
(\tw\theta,\rho) = \Psi_{\hat p,r}(x,y) = (rx,ry).
\end{equation*}
We will call $\Psi_{\hat p,r}$ a 
{\it boundary M\"obius chart} of radius $r$ centered
at $\hat p$.
Recall that $\hyp$ denotes the hyperboolic metric
on the upper half-space.

\begin{lemma}\label{lemma:bdry-mobius-charts}
There is a constant $C>0$ such that 
for any $\hat p\in \del M$ and any sufficiently small
$r>0$,
\begin{equation}\label{eq:bdry-g-estimate}
\|\Psi_{\hat p,r}^{*}g - \hyp\|_{l,\beta;Y_1} \le 
 rC.
\end{equation}
\end{lemma}

\begin{proof}
Because $Y_1$ is not precompact in $\H$, we have to interpret
the $C^{l,\beta}$ norm on the left-hand side
of \eqref{eq:bdry-g-estimate} as an intrinsic H\"older norm,
defined by \eqref{eq:def-holder-norm}.
For each point $(x_0,y_0)\in Y_1$, we have a M\"obius
chart $\Phi_{(x_0,y_0)}\colon B_2\to V_2(x_0,y_0)\subset  \H$
defined by
\begin{displaymath}
\Phi_{(x_0,y_0)}(x,y) = (x_0 + y_0 x, y_0 y).
\end{displaymath}
Then we need to get an upper bound for
\begin{displaymath}
\sup _{(x_0,y_0)\in Y_1} 
\|\Phi_{(x_0,y_0)}^* ( \Psi_{\hat p,r} ^* g - \hyp)
\|_{C^{l,\beta}(B_2)}.
\end{displaymath}
Since $\Phi_{(x_0,y_0)}$ is a hyperbolic isometry, the norm
above is the same as 
\begin{equation}\label{eq:phi-psi-g}
\| ( \Psi_{\hat p,r}\circ \Phi_{(x_0,y_0)}) ^* g - \hyp
\|_{C^{l,\beta}(B_2)}.
\end{equation}
In $(\tw\theta,\rho)$ coordinates,
we have
\begin{displaymath}
\Psi_{\hat p,r}\circ \Phi_{(x_0,y_0)}(x,y) = 
(r x_0 + r y_0 x, r y_0 y).
\end{displaymath}
Let us 
abbreviate this composite map as 
$\zeta(x,y) = (r x_0 + r y_0 x, r y_0 y)$, so that
\begin{displaymath}
( \Psi_{\hat p,r}\circ \Phi_{(x_0,y_0)}) ^* g - \hyp
 = y^{-2}\zeta^*(\overline g_{ij}-\delta_{ij}) dx^i \, dx^j.
\end{displaymath}
Since  the $C^{l,\beta}(B_2)$ norm in \eqref{eq:phi-psi-g} is
just the norm of the components in $(x,y)$-coordinates,
and $y^{-2}$ is uniformly bounded on $B_2$ together with all its
derivatives,
it suffices to show that
\begin{displaymath}
\|\zeta^*f\|_{C^{l,\beta}(B_2)} \le C r \|f\|_{C^{l,\beta}_{(0)}}
\end{displaymath}
for any function $f\in
C^{l,\beta}_{(0)}(\Omega)$ that 
vanishes at $\hat p$.
Moreover, the $(\tw\theta,\rho)$ coordinates 
are uniformly $C^{l+1,\beta}$-equivalent to the
original background coordinates $(\theta,\rho)$, because the coefficients
$B^\alpha$, the matrix
$(A^\alpha_\beta)$, and its inverse are uniformly bounded.
Thus the 
$C^{l,\beta}_{(0)}$ norm of $f$ in $(\tw\theta,\rho)$ coordinates
is uniformly bounded by the global norm
$\|f\|_{C^{l,\beta}_{(0)}}$.

To bound the sup norm of $\zeta^*f$, we use the mean value theorem
and the fact that $f(0,0) = f(\hat p) = 0$  to estimate
\begin{align*}
|\zeta^*f(x,y)|
&= |f(r x_0 + r y_0 x, r y_0 y) - f(0,0)|\\
&= 
\left| df _{(a_0,b_0)}
(r x_0 + r y_0 x, r y_0 y)\right|\\
&\le 
Cr\|f\|_{C^{1,0}_{(0)}}
\end{align*}
where $(a_0,b_0)$ is some point
on the line between $(0,0)$ and $(r x_0 + r y_0 x, r y_0 y)$.
For any coordinate $x^k$ $(k=1,\dots,n+1)$, we have
\begin{align*}
|\del_{x^k}(\zeta^{*}f)(x,y)| &= |r y_0 \del_{\theta^k}f(r x_0 + r y_0 x, r y_0 y)|\\
&\le r\|f\|_{C^{1,0}_{(0)}}.
\end{align*}
The H\"older norm of the first derivatives 
is estimated as follows:
\begin{align*}
&\frac{|\del_{x^k}(\zeta^{*}f)(x,y)-\del_{x^k}(\zeta^{*}f)(x',y')|}
{|(x,y) - (x',y')|^{\alpha}}\\
&\qquad= \frac{|ry_0 \del_{\theta^k}f(r x_0 + r y_0 x, r y_0 y)-
r y_0 \del_{\theta^k}f(r x_0 + r y_0 x', r y_0 y')|}
{|(x,y) - (x',y')|^{\alpha}}\\
&\qquad\le \frac{r \|f\|_{C^{1,\alpha}_{(0)}} |(r x_0 + r y_0 x, r y_0 y)-
 (r x_0 + r y_0 x', r y_0 y')|^\alpha}
{|(x,y) - (x',y')|^{\alpha}}\\
&\qquad\le Cr^{1+\alpha}\|f\|_{C^{1,\alpha}_{(0)}}.
\end{align*}
The general case now follows by induction on $l$, using
the fact that $\del_{x^k}(\zeta^{*}f)
= r y_0 \Phi_{\hat p}^* (\del_{\theta^k}f)$.
\end{proof}

We need to explore how the weighted Sobolev and H\"older norms
behave under boundary M\"obius charts.
The function $y$ is not a defining
function for hyperbolic space because it blows up at infinity;
however, by patching together via a partition of unity,
it is easy to construct a smooth defining function $\rho_0$ for
hyperbolic space that is equal to $y$ on $Y_1$.
Then for any boundary M\"obius chart $\Psi_{\hat p,r}$,
it follows that $\Psi_{\hat p,r}^* \rho = 
y =\rho_0$ on $Y_1$ for $r$ sufficiently
small, so 
by the same reasoning that led to 
\eqref{eq:mobius-uniform-holder} and \eqref{eq:mobius-uniform-sobolev},
the weighted norms have the following
scaling behavior under boundary M\"obius charts:
\begin{align}
C^{-1}r^{-\delta}\|\Psi_{\hat p,r}^*u\|_{k,\alpha;Y_1}
&\le \|u\|_{k,\alpha,\delta;Z_r(\hat p)}
\le Cr^{-\delta}\|\Psi_{\hat p,r}^*u\|_{k,\alpha;Y_1},
\label{eq:scaled-holder}\\
C^{-1}r^{-\delta}\|\Psi_{\hat p,r}^*u\|_{k,p;Y_1}
&\le \|u\|_{k,p,\delta;Z_r(\hat p)}
\le C r^{-\delta}\|\Psi_{\hat p,r}^*u\|_{k,p;Y_1}.
\label{eq:scaled-sobolev}
\end{align}

Choose a specific smooth bump function $\psi
\colon \H\to [0,1]$
that is equal to $1$ on $A_{1/2}$ and supported in $A_{1}$.
 For any $\hat p\in \del M$
and any $r>0$, let $(\tw\theta,\rho)$ be the coordinates
on a neighborhood $\Omega$ of $\hat p$ constructed above,
and define 
$\psi_{\hat p,r}\in C^{l,\beta}(\Omega)$
by 
\begin{displaymath}
\psi_{\hat p,r}(\tw\theta,\rho) = (\Psi_{\hat p,r}^{-1})^*\psi
=\psi(\tw\theta/r,\rho/r).
\end{displaymath}
Because the different choices of $(\tw\theta,\rho)$
coordinates are all uniformly
bounded in $C^{l+1,\beta}_{(0)}(\Omega)$ 
with respect to each other, and $|d\rho|_g$, $|d\tw\theta^\alpha|_g$
are both in $C^{l,\beta}_1(M)$
with norms independent of $\hat p$, it follows that the functions
$\psi_{\hat p,r}$ are uniformly bounded in $C^{l,\beta}(\Omega)$,
independently of $\hat p$ and $r$.

By the same argument as in Lemma 
\ref{covering-by-mobius-coordinates}, there
is a number $N$ such that for any $r>0$ we can choose
(necessarily finitely many) points $\{\hat p_1,\dots,\hat p_m\}\subset\del M$
such that the sets $\{Z_{r/2}(\hat p_i)\}$ cover 
$A_{r/2} = \{p\in M: \rho(p)<r/2\}$
and no more than $N$ of the sets $\{Z_{r}(\hat p_i)\}$
intersect nontrivially at any point.
For any such covering, let 
$\Psi_i = \Psi_{\hat p_i,r}$ and
$\psi_i = \psi_{\hat p_i,r}$.  Let 
$\psi_0\in C^\infty_c(M)$ be a smooth bump function that is
supported in $M\setminus A_{r/4}$ and 
equal to $1$ on $M\setminus A_{r/2}$,
and define
\begin{displaymath}
\phi_i = \frac{\psi_i}{\big(\sum_{j=0}^m \psi_i^2\big)^{1/2}}.
\end{displaymath}
It follows that $\{\phi_i^2\}$ is a partition of unity for $\overline M$ 
subordinate to the cover $\{M\setminus A_{r/4}, Z_{r}(\hat p_i)\}$.
Moreover,  at  each point of $\overline M$, at
least one of the functions $\psi_i$ is equal to $1$
at and at most $N$ of them
are nonzero, so the functions $\phi_i$ are still
uniformly bounded in $C^{l,\beta}_{(0)}(\overline M)$.

Let $\breve E$ be the tensor bundle over hyperbolic space
associated with the same $\Ortho(n+1)$ or 
$\SO(n+1)$ representation as $E$,
and let $\Phyp$ be the operator on hyperbolic space with the
same local coordinate expression as $P$.
For each boundary M\"obius chart $\Psi_i$, 
let $g_i$ be the metric $\Psi_i^* g$ defined
on $Y_1\subset \H$, and let
$P_i\colon C^\infty(Y_1;E)\to C^\infty(Y_1;E)$ be the operator
defined by
\begin{displaymath}
P_i u := \Psi_i^* P (\Psi_i^{-1*}u).
\end{displaymath}
Then 
Lemma \ref{lemma:bdry-mobius-charts} implies
that $P_i$ is close to $P$ 
in the following sense: 
For each $\delta\in \R$, $0<\alpha<1$,
$1<p<\infty$, and  $k$ such that
$m\le k\le l$ and $m<k+\alpha\le l+\beta$,
there is a constant
$C$ (independent of $r$ or $i$) such 
that for all $u\in C^{k,\alpha,\delta}(M;E)$,
\begin{equation}\label{eq:6.y}
\|P_iu - Pu\|_{k-m,\alpha,\delta} \le Cr\|u\|_{k,\alpha,\delta},
\end{equation}
and for all $u\in H^{k,p,\delta}(M;E)$,
\begin{equation}\label{eq:Pi-close-sobolev}
\|P_iu - Pu\|_{k-m,p,\delta} \le Cr\|u\|_{k,p,\delta}.
\end{equation}

Now assume that $P$ satisfies the hypotheses of
Theorem \ref{thm:main-fredholm}.  In particular, there is some
constant $C$ such that the $L^2$ estimate
\eqref{eq:asymptotic-L2-estimate} holds on the
complement of some compact set.
Choosing $\hat p\in \del M$ arbitrarily 
and $r$ sufficiently small, \eqref{eq:Pi-close-sobolev}
and \eqref{eq:scaled-sobolev} together
imply that $\Phyp$ satisfies an analogous
estimate (perhaps with a larger constant)
for all smooth sections $u$ of $\breve E$ compactly supported in 
$Y_1$.  But 
if $u\in C^\infty_c(\H;E)$ is 
arbitrary, 
there is a M\"obius transformation that takes 
$\supp u$ into $Y_1$, so the same estimate holds
globally on $\H$.  Therefore, by the results of
Chapter \ref{model-section}, $\Phyp$ is invertible
on $C^{k,\alpha}_\delta(\H;E)$ for $|\delta-n/2|<R$,
and on $H^{k,p}_\delta(\H;E)$ for $|\delta+n/p-n/2|<R$.

For any sufficiently small $r>0$,
define operators $Q_r,S_r,T_r\colon C^\infty_c(M;E)\to C^\infty_c(M;E)$ by
\begin{align*}
Q_r u &= \sum_i \phi_i (\Psi_i^{-1})^* \Phyp^{-1}\Psi_i^* (\phi_i u),\\
S_r u &= \sum_i \phi_i (\Psi_i^{-1})^* \Phyp^{-1}(P_i - \Phyp)\Psi_i^* (\phi_i u),\\
T_r u &= \sum_i \phi_i (\Psi_i^{-1})^* \Phyp^{-1}\Psi_i^* ([\phi_i,P]u).
\end{align*}

\begin{proposition}\label{prop:parametrix}
Let $P\colon C^\infty(M;E)\to C^\infty(M;E)$
satisfy the hypotheses of Theorem \ref{thm:main-fredholm}.
\begin{enumerate}\letters
\item\label{part:parametrix-sobolev}
If $|\delta +n/p-n/2|<R$ and $1<p<\infty$, 
then $Q_r$, $S_r$, and $T_r$  extend to
bounded maps as follows:
\begin{align*}
Q_r&\colon H^{0,p}_{\delta}(M;E) \to H^{m,p}_{\delta}(M;E),\\
S_r&\colon H^{m,p}_{\delta}(M;E) \to H^{m,p}_{\delta}(M;E),\\
T_r&\colon H^{m-1,p}_{\delta}(M;E)\to H^{m,p}_{\delta_1}(M;E),
\end{align*}
for any $\delta_1$ such that $\delta\le\delta_1\le \delta+1$ and
$|\delta_1+n/p-n/2|<R$.
Moreover, there exists $r_0>0$ such that 
if  $u\in H^{m,p}_{\delta}(M;E)$ is supported in
$A_{r}$ for $0<r<r_0$, then 
\begin{equation}\label{eq:parametrix-eqn}
Q_r Pu = u + S_r u +T_r u
\end{equation}
 and 
\begin{equation}\label{eq:sobolev-parametrix-estimate}
\|S_r u\|_{m,p,\delta}\le Cr \|u\|_{m,p,\delta}
\end{equation}
for some constant $C$ independent of $r$ and $u$.
\item\label{part:parametrix-holder}
If $|\delta -n/2|<R$, $0<\alpha<1$, 
and $m+\alpha\le l+\beta$,  
then $Q_r$, $S_r$, and $T_r$  extend to
bounded maps as follows:
\begin{align*}
Q_r&\colon C^{0,\alpha}_{\delta}(M;E) \to C^{m,\alpha}_{\delta}(M;E),\\
S_r&\colon C^{m,\alpha}_{\delta}(M;E) \to C^{m,\alpha}_{\delta}(M;E),\\
T_r&\colon C^{m-1,\alpha}_{\delta}(M;E)\to C^{m,\alpha}_{\delta_1}(M;E),
\end{align*}
for any $\delta_1$ such that $\delta\le\delta_1\le \delta+1$ and
$|\delta_1-n/2|<R$.
Moreover, there exists $r_0>0$ such that 
if $u\in C^{m,\alpha}_{\delta}(M;E)$ is supported in
$A_r$  for $0<r<r_0$, then 
\begin{equation}
Q_r Pu = u + S_r u +T_r u
\end{equation}
and
\begin{equation}\label{eq:holder-parametrix-estimate}
\|S_ru\|_{m,\alpha,\delta}\le Cr \|u\|_{m,\alpha,\delta}
\end{equation}
for some constant $C$ independent of $r$ and $u$.
\end{enumerate}
\end{proposition}

\begin{proof}
The fact that $Q_rPu = u + S_ru +T_ru$ in $A_r$
is just a computation:
\begin{align*}
Q_rPu
&= \sum_i \phi_i (\Psi_i^{-1})^* \Phyp^{-1}\Psi_i^* (\phi_i Pu)\\
&= \sum_i \phi_i (\Psi_i^{-1})^* \Phyp^{-1}\Psi_i^* (P(\phi_i u))
+ \sum_i \phi_i (\Psi_i^{-1})^* \Phyp^{-1}\Psi_i^* ([\phi_i, P]u)\\
&= \sum_i \phi_i (\Psi_i^{-1})^* \Phyp^{-1}P_i\Psi_i^* (\phi_i u)
+ T_ru\\
&= \sum_i \phi_i (\Psi_i^{-1})^* \Phyp^{-1}\Phyp\Psi_i^* (\phi_i u)\\
&\qquad + \sum_i \phi_i (\Psi_i^{-1})^* \Phyp^{-1}(P_i-\Phyp)\Psi_i^* (\phi_i u)
+ T_ru\\
&= u + S_ru + T_ru.
\end{align*}

To check the mapping properties of $S_r$, we begin by observing that
the fact that the functions 
$\phi_i$ are uniformly bounded in $C^{l,\beta}_{(0)}(\overline M)
\subset C^{l,\beta}(M)$
implies by Lemma 
\ref{properties-of-spaces}\eqref{part:mult-cts}
that 
multiplication by $\phi_i$ is a bounded map from 
$H^{j,p}_{\delta}(Z_r(\hat p_i);E)$
to itself for each $i$ and all $0\le j\le l$, with norm bounded
independently of $i$ and $r$.
The fact that $S_r$ maps $H^{m,p}_{\delta}(M;E)$
to itself 
then follows from 
Proposition \ref{prop:K-maps-to-Linfty} and
\eqref{eq:scaled-sobolev}, because the factors of
$r^\delta$ and $r^{-\delta}$ introduced by $\Psi_i^*$ and its
inverse cancel each other.
Moreover, \eqref{eq:6.y} implies \eqref{eq:sobolev-parametrix-estimate}
whenever $u$ is supported in $A_r$, for some constant
$C$ independent of $r$ and $u$. 

The mapping properties of $T_r$ will
follow from a similar argument once we show that the
commutator $[\phi_i,P]$
maps $H^{m-1,p}_{\delta}(M;E)$
to $H^{m,p}_{\delta_1}(M;E)$.
Observe that each term in the coordinate expression for $[\phi_i,P]u$
is a product of four factors: a constant, 
a $p$th covariant derivative of $u$, 
a $q$th covariant derivative of $\phi_i$, 
and a polynomial in the components of 
$g$, $(\det g)^{-1/2}$, and their  derivatives up through
order $r$, with $p+q+r\le m$ and $q\ge 1$.
Since $\phi_i$ is uniformly bounded 
in $C^{l,\beta}(M)$,
the result follows.
The argument for the H\"older case is identical.
\end{proof}

\begin{corollary}\label{cor:parametrix}
Let $P\colon C^\infty(M;E)\to C^\infty(M;E)$
satisfy the hypotheses of Theorem \ref{thm:main-fredholm}.
\begin{enumerate}\letters
\item
If $|\delta +n/p-n/2|<R$, $|\delta_1+n/p-n/2|<R$, 
$\delta\le\delta_1\le \delta+1$, 
and $1<p<\infty$, 
then there exist $r>0$ and bounded operators
\begin{align*}
\tw Q&\colon H^{0,p}_{\delta}(M;E) \to H^{m,p}_{\delta}(M;E),\\
\tw T&\colon H^{m-1,p}_{\delta}(M;E)\to H^{m,p}_{\delta_1}(M;E)
\end{align*}
such that 
\begin{displaymath}
\tw Q Pu = u + \tw T u
\end{displaymath} 
whenever 
$u\in H^{m,p}_{\delta}(M;E)$ is supported in 
$A_{r}$.
\item
If $|\delta -n/2|<R$, $|\delta_1-n/2|<R$,
$\delta\le \delta_1\le\delta+1$, $0<\alpha<1$, 
and $m+\alpha\le l+\beta$,  
then there exist $r>0$ and bounded operators
\begin{align*}
\tw Q&\colon C^{0,\alpha}_{\delta}(M;E) \to C^{m,\alpha}_{\delta}(M;E),\\
\tw T&\colon C^{m-1,\alpha}_{\delta}(M;E)\to C^{m,\alpha}_{\delta_1}(M;E)
\end{align*}
such that 
\begin{equation}\label{eq:strong-parametrix-eqn}
\tw Q Pu = u + \tw T u
\end{equation} 
whenever 
$u\in C^{m,\alpha}_{\delta}(M;E)$ is supported in 
$A_{r}$.
\end{enumerate}
\end{corollary}

\begin{proof}
Just choose $r$ small enough that 
\eqref{eq:sobolev-parametrix-estimate}
holds with $Cr<1/2$.  
It then follows that
$\Id +S_r\colon H^{m,p}_\delta(M;E) \to H^{m,p}_\delta(M;E) $
has a bounded inverse.  We just set 
\begin{align*}
\tw Q = (\Id + S_r)^{-1}\circ Q_r,\\
\tw T = (\Id + S_r)^{-1}\circ T_r,
\end{align*}
and then \eqref{eq:strong-parametrix-eqn} follows immediately from 
\eqref{eq:parametrix-eqn}.
Once again, the argument for the H\"older case is identical.
\end{proof}

Our first application of this parametrix construction
is a significant strengthening
of Lemma \ref{rescaling-estimates}, giving improved
decay for solutions to $Pu=f$ when
$u$ and $f$ are in appropriate spaces.
We begin with a special case.

\begin{lemma}\label{lemma:optimal-decay}
Assume $P$ satisfies the hypotheses 
of Theorem \ref{thm:main-fredholm}.
\begin{enumerate}\letters
\item\label{part:optimal-sobolev-decay}
Suppose that 
$1<p<\infty$, $m\le k \le l$,
$|\delta +n/p-n/2|<R$,
and $|\delta '+n/p-n/2|<R$.
If $u\in H^{0,p}_{\delta}(M;E)$ 
and $Pu\in H^{k-m,p}_{\delta'}(M;E)$,
then $u\in H^{k,p}_{\delta'}(M;E)$.
\item\label{part:optimal-holder-decay}
Suppose that 
$0<\alpha<1$, $m< k+\alpha\le l+\beta$,
$|\delta -n/2|<R$,
and $|\delta '-n/2|<R$.
If $u\in C^{0,0}_{\delta}(M;E)$ 
and $Pu\in C^{k-m,\alpha}_{\delta'}(M;E)$,
then $u\in C^{k,\alpha}_{\delta'}(M;E)$.
\end{enumerate}
\end{lemma}

\begin{proof}
If $\delta'\le \delta$, the result is a trivial
consequence of Lemmas 
\ref{rescaling-estimates} and 
\ref{properties-of-spaces}\eqref{inclusions},
so assume $\delta'>\delta$.
Consider part \eqref{part:optimal-sobolev-decay}.
By Lemma \ref{rescaling-estimates}, it suffices to show
that $u\in H^{0,p}_{\delta'}(M;E)$.
For any small $r>0$,
by means of a bump function we can write $u=u_0+u_\infty$,
where $\supp u_0$ is compact and $\supp u_\infty\subset A_r$.
Local elliptic regularity gives 
$u_0\in H^{k,p}_{\delta'}(M;E)$.
Since $Pu_\infty$ agrees with $Pu$ off of a compact
set, $Pu_\infty\in H^{0,p}_{\delta'}(M;E)$, so 
by Corollary \ref{cor:parametrix},
if $r$ is small enough, 
\begin{align*}
u_\infty &= \tw QPu_\infty - \tw Tu_\infty\\
&\in H^{m,p}_{\delta'}(M;E) + H^{1,p}_{\delta_1}(M;E)\\
&\subset H^{0,p}_{\delta_1}(M;E),
\end{align*}
where $\delta_1 = \min(\delta',\delta+1)$.
Iterating this argument finitely 
many times, we conclude that $u_\infty\in H^{0,p}_{\delta'}(M;E)$.
By Lemma \ref{properties-of-spaces},
this implies that 
$u\in H^{k,p}_{\delta'}(M;E)$ as claimed.
The argument for the H\"older case is the same.
\end{proof}

\begin{proposition}\label{prop:optimal-regularity}
Suppose $P$ satisfies the hypotheses 
of Theorem \ref{thm:main-fredholm}, and $u$
is either in $H^{0,p_0}_{\delta_0}(M;E)$ for some 
$|\delta _0+n/p_0-n/2|<R$ and $1<p_0<\infty$, or in
$C^{0,0}_{\delta_0}(M;E)$ for some
$|\delta _0-n/2|<R$.
\begin{enumerate}\letters
\item\label{part:optimal-reg-sobolev}
If $Pu \in 
H^{k-m,p}_{\delta}(M;E)$ 
for $|\delta +n/p-n/2|<R$, $1<p<\infty$, and $m\le k \le l$,
then $u\in H^{k,p}_{\delta}(M;E)$.
\item\label{part:optimal-reg-holder}
If $Pu \in 
C^{k-m,\alpha}_{\delta}(M;E)$ 
for $|\delta -n/2|<R$, $0<\alpha<1$, and $m< k+\alpha \le l+\beta$,
then $u\in C^{k,\alpha}_{\delta}(M;E)$.
\end{enumerate}
\end{proposition}

\begin{proof}
If $u\in C^{0,0}_{\delta_0}(M;E)$ with 
$|\delta _0-n/2|<R$, then $u\in H^{0,p}_{\delta}(M;E)$
whenever $\delta + n/p < \delta_0$ by Lemma \ref{properties-of-spaces}.
Since such $\delta$ and $p$ can be chosen that also satisfy
$|\delta  + n/p - n/2|<R$, 
it suffices to prove the proposition under
the hypothesis that 
$u\in H^{0,p_0}_{\delta_0}(M;E)$ with 
$|\delta _0+n/p_0-n/2|<R$.  
For the rest of the proof, 
we assume this.

First we treat case \eqref{part:optimal-reg-sobolev}.
Assume that $Pu\in H^{k-m,p}_{\delta}(M;E)$
with $|\delta +n/p-n/2|<R$, and 
let $\scr P$ be the following set: 
\begin{displaymath}
\scr P = 
\{ p'\in (1,\infty):
u\in H^{0,p'}_{\delta'}(M;E)
\text{ for some $\delta'$ with $|\delta '+n/p'-n/2|<R$}\}.
\end{displaymath}
Clearly $p_0\in \scr P$ by hypothesis.
We will show that $p\in \scr P$.
It will then follow from Lemma \ref{lemma:optimal-decay} 
that $u\in H^{k,p}_{\delta}(M;E)$, which will prove case 
\eqref{part:optimal-reg-sobolev}.

{\sc Claim 1:}
If $p_1\in \scr P$, then $(1,p_1]\subset \scr P$.
To prove this, assume $p_1\in \scr P$ and $1<p'<p_1$.
The fact that $p_1\in \scr P$ means that there is
some $\delta_1$ with $|\delta _1+n/p_1-n/2|<R$ such that $u\in H^{0,p_1}_{\delta_1}(M;E)$.
By Lemma \ref{properties-of-spaces}, $u\in H^{0,p'}_{\delta'}(M;E)$
for any $\delta'$ such that $\delta_1+n/p_1 > \delta'+n/p'$.  Choosing $\delta'$ 
so that $\delta'+n/p'$ is sufficiently close to $\delta_1+n/p_1$,
we can ensure that $|\delta '+n/p'-n/2|<R$.
This implies that $p'\in\scr P$ as desired.

{\sc Claim 2:}
If $p_1\in \scr P$ and $p_2$ satisfies $p_1<p_2\le p$ and
\begin{equation}\label{eq:choice-of-p2}
\frac{p_2}{p_1} \le \min\left(\frac{n+2}{n+1},1+\frac{\epsilon}{2n}\right),
\end{equation}
where 
\begin{displaymath}
\epsilon = \delta + \frac n p - \frac n 2 + R > 0,
\end{displaymath}
then $p_2\in \scr P$.
The assumption that $p_1\in \scr P$ means that $u\in H^{0,p_1}_{\delta_1}(M;E)$
for some $\delta_1$ with $|\delta _1+n/p_1-n/2|<R$.  
Choose $\delta'$ satisfying 
\begin{equation}\label{eq:choice-of-s'}
\delta+\frac n p >
\delta'+\frac n {p_1} > \delta + \frac n p - \frac{\epsilon}{2} .
\end{equation}
By virtue of the first inequality above, 
Lemma \ref{properties-of-spaces} implies that
\begin{equation}\label{eq:Pu-in-spaces}
Pu \in H^{k-m,p}_{\delta}(M;E)\subset
H^{k-m,p_1}_{\delta'}(M;E).
\end{equation}
Since the two inequalities of
\eqref{eq:choice-of-s'} guarantee that $|\delta '+n/p_1-n/2|<R$, 
Lemma \ref{lemma:optimal-decay} shows that 
$u \in  H^{k,p_1}_{\delta'}(M;E)$.
Our restriction on $p_2$ guarantees that
\begin{align*}
\frac{n+1}{p_1}
&\le \frac{n+2}{p_2}\\
&= \frac{n+1}{p_2} + \frac{1}{p_2}\\
&\le \frac{n+1}{p_2} + k,
\end{align*}
so Lemma \ref{properties-of-spaces}\eqref{prop:Sobolev}
implies that $u\in H^{0,p_2}_{\delta'}(M;E)$.
Now \eqref{eq:choice-of-p2}
implies that 
\begin{align*}
\frac{n}{p_1} - \frac{n}{p_2} 
&= \frac{n}{p_2}\left(\frac{p_2}{p_1} - 1\right)\\
&\le \frac{n}{p_2}\left( \frac{\epsilon}{2n}\right)\\
&< \frac{\epsilon}{2}.
\end{align*}
Therefore, using \eqref{eq:choice-of-s'}, we obtain
\begin{align*}
\delta' + \frac{n}{p_2}-\frac{n}{2}
&= \left(\delta' + \frac{n}{p_1}-\frac{n}{2}\right) -
\left(\frac{n}{p_1} - \frac{n}{p_2}\right)\\
&> \left(\delta + \frac n p -\frac{n}{2}- \frac{\epsilon}{2}\right) -
\left(\frac{\epsilon}{2}\right)\\
&=-R,\\
\delta' + \frac{n}{p_2}-\frac{n}{2}
&<\delta + \frac n p - \frac n 2\\
&<R,\\
\end{align*}
which proves that $p_2\in\scr P$ as claimed.

{\sc Claim 3:} $p\in \scr P$.
If $p\le p_0$, this follows immediately
from Claim 1 together with the obvious fact that $p_0\in\scr P$.
Otherwise, just iterate Claim 2, starting with $p_0\in \scr P$.
After finitely many iterations, we can conclude that $p\in \scr P$. 

Finally we turn to case
\eqref{part:optimal-reg-holder}.
Suppose $Pu\in C^{k-m,\alpha}_{\delta}(M;E)$ with $|\delta -n/2|<R$,
and choose $p'$ large and $\delta'$ close to $\delta$
satisfying
\begin{equation}\label{eq:large-p'}
\delta >\delta'+\frac n{p'} ,
\end{equation}
Then by Lemma \ref{properties-of-spaces}, 
$Pu\in H^{k-m,p'}_{\delta'}(M;E)$.
If we choose $\delta'+n/p'$ sufficiently close to $\delta$, 
we have $|\delta '+n/p'-n/2|<R$, and thus
$u\in H^{k,p'}_{\delta'}(M;E)$
by part \eqref{part:optimal-reg-sobolev} above.
If $p$ is also chosen large enough that $(n+1)/p' \le k-\alpha$,
the Sobolev embedding theorem 
(Theorem \ref{properties-of-spaces}\eqref{prop:Sobolev}) implies
$u\in C^{0,\alpha}_{\delta'}(M;E)$,
and if $\delta'$ is sufficiently close to $\delta$ we will 
have $|\delta '-n/2|<R$.
Then Lemma \ref{lemma:optimal-decay} implies $u\in C^{k,\alpha}_{\delta}(M;E)$.
\end{proof}

Now suppose $P\colon C^\infty(M;E)\to C^\infty(M;E)$ satisfies
the hypotheses of Theorem \ref{thm:main-fredholm}.  
By Lemma \ref{lemma:fred},
estimate \eqref{eq:asymptotic-L2-estimate}
implies that $P\colon H^{m,2}(M;E)\to H^{0,2}(M;E)$
is Fredholm. 
Let $Z = \Ker P\cap L^2(M;E)$, which is equal to 
$\Ker P \cap H^{m,2}(M;E)$ by Lemma \ref{rescaling-estimates}.
Then $Z$ is finite-dimensional, and 
Proposition \ref{prop:optimal-regularity}
shows that 
$Z\subset H^{k,p}_{\delta}(M;E)$ whenever 
$1<p<\infty$, $m\le k \le l$, and $|\delta +n/p-n/2|<R$.
In fact,
\begin{equation}\label{eq:Z=kerP}
Z= \Ker P\colon H^{k,p}_{\delta}(M;E)
\to H^{k-m,p}_\delta(M;E),
\end{equation}
because any $u\in \Ker P\cap H^{k,p}_{\delta}(M;E)$
is also in $L^2(M;E)$
by Proposition \ref{prop:optimal-regularity}.
Similarly, if 
$0<\alpha<1$, $m< k+\alpha \le l+\beta$, and $|\delta -n/2|<R$,
then $Z= \Ker P\colon C^{k,\alpha}_{\delta}(M;E)\to C^{k-m,\alpha}_\delta(M;E)$.

Because of these observations, whenever 
$1<p<\infty$, $m\le k \le l$, and $|\delta +n/p-n/2|<R$,
we can define a subspace $Y^{k,p}_\delta\subset H^{k,p}_\delta(M;E)$
by 
\begin{displaymath}
Y^{k,p}_\delta=\{u\in  H^{k,p}_\delta(M;E): (u,v) = 0 \text{ for all $v\in Z$}\},
\end{displaymath}
where $(u,v)$ represents the standard $L^2$ pairing.
Since $Z\subset H^{0,p^*}_{-\delta}(M;E) = 
(H^{0,p}_\delta(M;E))^*
\subset (H^{k,p}_\delta(M;E))^*$,
it follows that $Y^{k,p}_\delta$ is a well-defined
closed subspace of $H^{k,p}_\delta(M;E)$.
Similarly, if 
$0<\alpha<1$, $m< k+\alpha \le l+\beta$, and $|\delta -n/2|<R$,
we can define
\begin{displaymath}
Y^{k,\alpha}_\delta=\{u\in  C^{k,\alpha}_\delta(M;E): (u,v) = 0 \text{ for all $v\in Z$}\},
\end{displaymath}
since $C^{k,\alpha}_{\delta}(M;E) \subset H^{0,p}_{\delta'}(M;E)$
for $\delta>\delta'+n/p>n/2-R$ implies
$Z\subset H^{0,p^*}_{-\delta'}(M;E) =
(H^{0,p}_{\delta'}(M;E))^*
\subset (C^{k,\alpha}_{\delta}(M;E))^*$.

The next result
is the main structure theorem
for operators satisfying the 
hypotheses of Theorem \ref{thm:main-fredholm}.

\begin{theorem}\label{thm:structure}
Suppose $P\colon C^\infty(M;E)\to C^\infty(M;E)$
satisfies the hypotheses of Theorem \ref{thm:main-fredholm}.
\begin{enumerate}\letters
\item\label{part:sobolev-structure}
If $1<p<\infty$, $0\le k\le l$, and $|\delta +n/p-n/2|<R$,
there exist bounded operators $G,H \colon H^{k,p}_\delta(M;E)\to 
H^{k,p}_\delta(M;E)$ such that 
$G(H^{k-m,p}_\delta(M;E))\subset H^{k,p}_\delta(M;E)$
for $k\ge m$, and 
\begin{align}
Y^{k,p}_\delta &= \Ker H ,\label{eq:kerH}\\
Z &= \Image H ,\label{eq:imH}\\
u &= GPu + H u \text{ for $u\in H^{k,p}_\delta(M;E)$, $m\le k\le l$},\label{eq:GPu}\\
u &= PGu + H u \text{ for $u\in H^{k,p}_\delta(M;E)$, $0\le k\le l$}.\label{eq:PGu}
\end{align}
\item\label{part:holder-structure}
If $0<\alpha<1$, $0< k+\alpha\le l+\beta$, and $|\delta-n/2|<R$,
there exist bounded operators $G,H \colon C^{k,\alpha}_\delta(M;E)\to 
C^{k,\alpha}_\delta(M;E)$ such that 
$G(C^{k-m,\alpha}_\delta(M;E))\subset C^{k,\alpha}_\delta(M;E)$
for $k\ge m$, and 
\begin{align}
Y^{k,\alpha}_\delta &= \Ker H ,\label{eq:kerH-holder}\\
Z &= \Image H ,\\
u &= GPu + H u \text{ for $u\in C^{k,\alpha}_\delta(M;E)$, $m\le k\le l$},\\
u &= PGu + H u \text{ for $u\in C^{k,\alpha}_\delta(M;E)$, 
$0\le k\le l$}.\label{eq:PGu-holder}
\end{align}
\end{enumerate}
\end{theorem}

\begin{proof}
We begin with the Sobolev
case, part \eqref{part:sobolev-structure}.
First consider the special case $p=2$, $k=0$, and $\delta=0$
(in which case this is basically the standard construction
of a partial inverse for a Fredholm operator on $L^2$).
As noted above, 
the assumption of an $L^2$ estimate \eqref{eq:asymptotic-L2-estimate}
near the boundary implies that 
$P\colon H^{m,2}_0(M;E)\to L^2(M;E)$ is Fredholm
by Lemma \ref{lemma:fred}.  

By definition, $Y^{0,2}_0\subset L^2(M;E)$
is precisely the orthogonal complement of $Z=\Ker P$
in $L^2(M;E)$, so we have an orthogonal 
direct sum decomposition
$L^2(M;E) = Z \oplus Y^{0,2}_0$.
Since $P$ is formally self-adjoint, any 
$u\in H^{m,2}(M;E)$
satisfies $(Pu,v) = (u,Pv) = 0$
for all $v\in Z$, so 
$P(H^{m,2}(M;E))\subset Y^{0,2}_0$.
On the other hand, if $v\in L^2(M;E)$ is
orthogonal to $P(H^{m,2}(M;E))$,
then for any smooth, compactly supported section
$u$ of $E$ we have $(v,Pu) = 0$, so
$v$ is a distributional solution to $Pv=0$, which means
$v\in Z$.  This shows that $(P(H^{m,2}(M;E)))^\perp = Z$,
and since $P(H^{m,2}(M;E))$ is closed in $L^2(M;E)$
we have $P(H^{m,2}(M;E)) = Y^{0,2}_0$.

Now $P\colon Y^{m,2}_0\to Y^{0,2}_0$ is bijective
and bounded, 
so by the open mapping theorem it has a
bounded inverse $(P|_{Y^{m,2}_0})^{-1}\colon 
Y^{0,2}_0\to Y^{m,2}_0$.  Define $G\colon L^2(M;E)
\to L^2(M;E)$ by 
\begin{displaymath}
Gu = 
\begin{cases}
(P|_{Y^{m,2}_0})^{-1}u & u\in Y^{0,2}_0,\\
0 & u\in Z,
\end{cases}
\end{displaymath}
and define $H \colon L^2(M;E)\to L^2(M;E)$ 
to be the orthogonal projection onto $Z$.
Then 
\eqref{eq:kerH} and \eqref{eq:imH} are immediate
from the definition of $H $, and 
\eqref{eq:GPu} (for $u\in H^{m,2}(M;E)$)
and \eqref{eq:PGu} (for all $u$) follow 
by considering $u\in Y^{0,2}_0$ and $u\in Z$
separately.

Next consider the case of arbitrary 
$\delta$ satisfying $|\delta|<R$, 
still with $p=2$ and $k=0$.
If $\delta>0$
and $u\in H^{0,2}_\delta(M;E)\subset L^2(M;E)$,
then $PGu = u - H u \in H^{0,2}_\delta(M;E)$,
so $Gu\in H^{m,2}_\delta(M;E)$ by 
Lemma \ref{lemma:optimal-decay}.
Thus the restriction of
$G$ to $H^{0,2}_\delta(M;E)$
takes its values in $H^{m,2}_\delta(M;E)$,
as does
$H $ by \eqref{eq:Z=kerP}.
In this case, \eqref{eq:kerH}--\eqref{eq:PGu} 
are satisfied because they are
already satisfied on the bigger space $L^2(M;E)$ (or $H^{m,2}(M;E)$
in case of \eqref{eq:GPu}),
and $Y^{0,2}_\delta = Y^{0,2}_0\cap H^{0,2}_\delta(M;E)$.

On the other hand, if $\delta<0$, we can use the
fact that $H^{0,2}_{\delta}(M;E)=(H^{0,2}_{-\delta}(M;E))^*$
to extend
the definition of $G$ and $H $ to $H^{0,2}_{\delta}(M;E)$
by duality:  For any 
$u\in H^{0,2}_\delta(M;E)$,
let $Gu$ and $H u$ 
be the elements of 
$H^{0,2}_\delta(M;E)$ defined uniquely by 
\begin{align}
(Gu,v)&= (u,Gv),\label{eq:G-adjoint}\\
(H u,v)&= (u,H v)\label{eq:H-adjoint}
\end{align}
for all $v\in H^{0,2}_{-\delta}(M;E)$.
In other words, $G,H\colon H^{0,2}_{\delta}(M;E)\to
H^{0,2}_\delta(M;E)$
are defined to be the dual maps of 
$G,H\colon H^{0,2}_{-\delta}(M;E)\to
H^{0,2}_{-\delta}(M;E)$.  Since $H$ and $G$ are
self-adjoint on $L^2(M;E)$ ($H$ because it is an orthogonal projection,
and $G$ because $P$ is self-adjoint as an unbounded
operator), these are indeed extensions of the original
maps $G$ and $H$.

To see that 
these extended operators 
satisfy \eqref{eq:kerH}--\eqref{eq:PGu},
we observe that 
\eqref{eq:H-adjoint}
implies that $Hu=0$ for $u\in H^{0,2}_\delta(M;E)$
exactly when $(u,v)=0$ 
for all $v$ in the image of 
$H\colon H^{0,2}_{-\delta}(M;E)
\to H^{0,2}_{-\delta}(M;E)$; since this image is 
exactly $Z$, it follows that $\Ker H = Y^{0,2}_\delta$,
which is \eqref{eq:kerH}.
Since the restriction of $H$ to 
$Z\subset L^2(M;E)\subset H^{0,2}_{\delta}(M;E)$ is
the identity, it follows that $Z\subset\Image H$.
On the other hand, for any
$u\in H^{0,2}_\delta(M;E)$,
we have $(Hu,Pv) = (u,H Pv) = 0$ for all $v\in C^\infty_c(M;E)
\subset H^{0,2}_{-\delta}(M;E)$,
which means that $H u$ is a weak solution to $P(H u)=0$.  Thus
$\Image H  \subset \Ker P = Z$,
which proves
\eqref{eq:imH}.
Equations \eqref{eq:GPu}
and \eqref{eq:PGu} then follow easily from our definitions
by duality.

Next we generalize to $1<p<\infty$ and $|\delta+n/p-n/2|<R$,
still with $k=0$.  If $p>2$, we can choose $\delta'$ 
such that $|\delta'|<R$ and $\delta+n/p>\delta'+n/2$, so that
$H^{0,p}_\delta(M;E)\subset H^{0,2}_{\delta'}(M;E)$.
Arguing as above, we see that the restrictions of $G$ and $H $ map
$H^{0,p}_\delta(M;E)$ to $H^{m,p}_\delta(M;E)$.
On the other hand, for $p<2$, we can extend $G$ and $H $ 
to maps from $H^{0,p}_\delta(M;E) = (H^{0,p^*}_{-\delta}(M;E))^*$ to itself by
duality as above.  In both cases \eqref{eq:kerH}--\eqref{eq:PGu}
are satisfied, by restriction or duality as appropriate.

Now consider the general case of 
$H^{k,p}_{\delta}(M;E)$ with 
$0\le k\le l$, $1<p<\infty$, and $|\delta+n/p-n/2|<R$.
Since $Z\subset H^{k,p}_{\delta}(M;E)$, 
it is clear that $H $ restricts to a map of 
$H^{k,p}_{\delta}(M;E)$ to itself.
If $u\in 
H^{k,p}_\delta(M;E)\subset H^{0,p}_\delta(M;E)$ for 
$0\le k \le l-m$, 
observe as above that
$PGu = u-H u\in H^{k,p}_\delta(M;E)$, so 
$Gu\in H^{k+m,p}_\delta(M;E)$ by Lemma \ref{rescaling-estimates}.
For $l-m\le k \le l$, we have
$G(H^{k,p}_\delta(M;E))\subset G(H^{l-m,p}_\delta(M;E))
\subset H^{l,p}_\delta(M;E)\subset H^{k,p}_\delta(M;E)$.
Thus in each case $G$ and $H $ restrict to maps from
$H^{k,p}_\delta(M;E)$ to itself, and 
properties \eqref{eq:kerH}--\eqref{eq:PGu}
are satisfied by restriction.

Finally, consider case \eqref{part:holder-structure},
and assume that $0<\alpha<1$, $0< k+\alpha \le l+\beta$,
and $|\delta-n/2|<R$.
We can choose $\delta'\in\R$ satisfying $|\delta'|<R$
and $\delta >\delta'+n/2$, so that
$C^{k,\alpha}_\delta(M;E)\subset H^{k,2}_{\delta'}(M;E)$
and the results of part \eqref{part:sobolev-structure}
apply to $H^{k,2}_{\delta'}(M;E)$.
Then the restrictions of $G$ and $H $ map 
$C^{k,\alpha}_{\delta}(M;E)$ to itself by the same
argument as above, and properties 
\eqref{eq:kerH-holder}--\eqref{eq:PGu-holder}
are automatically satisfied by restriction.
\end{proof}

The next construction will be useful in proving that
$P$ is not Fredholm outside the expected  range of weights.

\begin{lemma}\label{lemma:critical-asymp-solutions}
Suppose $P\colon C^\infty(M;E)\to C^\infty(M;E)$
satisfies the hypotheses of Theorem \ref{thm:main-fredholm},
where $E$ is a bundle of tensors of weight $r$.  
Let $s_0$ be a characteristic exponent of $P$,
and let $\delta_0 = \Real s_0+r$.
Given any compact subset $K\subset M$,
there is an infinite-dimensional subspace
$W\subset C^{l,\beta}_{\delta_0}(M;E)$ such that
every nonzero $w\in W$ has the following properties.
\begin{enumerate}\letters
\item\label{part:supp-w}
$\supp w\subset M\setminus K$.
\item
$Pw\in C^{0,0}_{\delta_0+1}(M;E)$.
\item\label{part:notin-Hp}
If $1<p<\infty$ and $\delta\ge \delta_0-n/p$, then
$w \notin H^{0,p}_{\delta}(M;E)$.
\end{enumerate}
\end{lemma}

\begin{proof}
Let $\hat p\in \del M$ be arbitrary, and let
$V$ be any neighborhood of $\hat p$ in $\overline M$.
Since the characteristic
exponents are constant on 
$\del M$ by Lemma \ref{lemma:const-char-exps},
there is a tensor $\overline w_{\hat p}\in E_{\hat p}$
such that $I_{s_0}(P)\overline w_{\hat p} = 0$.
We can extend $\overline w_{\hat p}$ to a $C^{l,\beta}$
tensor field $\overline w$ 
on a neighborhood of ${\hat p}$ in $\del M$,
still satisfying $I_{s_0}(P)\overline w=0$, as follows.
Shrinking $V$
if necessary, we can choose background coordinates on $V$,
and for
each $\hat q\in V\cap \del M$, let $A(\hat q)$
be the matrix of $I_{s_0}(P)\colon E_{\hat q}\to E_{\hat q}$.
By Lemma 
\ref{lemma:Clbeta-indicial-map},
the matrix entries of
$A(\hat q)$ are $C^{l,\beta}$ functions of $\hat q$.
If $\gamma$ is any smooth, positively-oriented 
closed curve in $\C$
whose interior contains $0$ but no other eigenvalues of 
$A(\hat p)$,
then the projection onto the kernel of $A(\hat q)$
can be written as $-1/(2\pi i)\int_\gamma (A(\hat q)-z \Id)^{-1} dz$.
(Here we use the fact that the eigenvalues of $A(\hat q)$
and their multiplicities are independent of $\hat q$.)
We define $\overline w$ to be the tensor field on $V\cap\del M$
whose coordinate expression is 
\begin{displaymath}
\overline w_{\hat q} = 
\frac{-1}{2\pi i}
\int_\gamma (A(\hat q)-z \Id)^{-1}\overline w_{\hat p} dz,
\end{displaymath}
which is a $C^{l,\beta}$ tensor field along $\del M$
satisfying $I_{s_0}(P)\overline w=0$.
If we extend $\overline w$ arbitrarily to a $C^{l,\beta}$
tensor field on $V$, and let
$w = \rho^{s_0}\phi \overline w$ where 
$\phi$ is any smooth cutoff function
that is positive at $\hat p$ and supported in $V$,
then $w\in C^{l,\beta}_{\delta}(M;E)$,
and 
\begin{displaymath}
|Pw|_{\overline g} = O(\rho^{\Real s_0+1}),
\end{displaymath}
which implies
\begin{equation}\label{eq:Pw-estimate}
|Pw|_{g} = O(\rho^{\delta_0+1}),
\end{equation}
so $Pw\in C^{0,0}_{\delta_0+1}(M;E)$.
On the other hand, 
$w \notin H^{0,p}_{\delta}(M;E)$
for $\delta\ge \delta_0-n/p$ by Lemma \ref{lemma:being-in-Lp}.

Now
choose countably many points $\hat p_i\in \del M$ and 
disjoint neighborhoods $V_i$ of $\hat p_i$.
For each $i$ we can construct $w_i$ as above with
support in $V_i$, so the space spanned by $\{w_i\}$
is clearly infinite-dimensional.
\end{proof}

\begin{proof}[Proof of Proposition \ref{prop:L2-Fredholm}]
It follows immediately from Lemma \ref{lemma:fred} that
$P$ is Fredholm as a map from $H^{m,2}(M;E)$ to $H^{0,2}(M;E)$ 
if and only if $P$ satisfies 
an estimate of the form \eqref{eq:asymptotic-L2-estimate}.
By Lemma \ref{rescaling-estimates}, the kernel and range
of $P$ as an unbounded operator on $L^2(M;E)$ are
the same as those of $P\colon H^{m,2}(M;E)\to H^{0,2}(M;E)$.
 The proposition follows.
\end{proof}

Finally, we are in a position to prove our main
Fredholm theorem, Theorem
\ref{thm:main-fredholm} from the Introduction.

\begin{proof}[Proof of Theorem \ref{thm:main-fredholm}]
We will prove parts \eqref{sobolev-fredholm-result}
and \eqref{holder-fredholm-result} together.
The proof of sufficiency is identical for the Sobolev
and H\"older cases, so we do only the Sobolev case.

Suppose $1<p<\infty$, $m\le k\le l$, and $|\delta +n/p-n/2|<R$,
and let $G,H \colon H^{k,p}_{\delta}(M;E)\to H^{k,p}_{\delta}(M;E)$
be as in Theorem \ref{thm:structure}.
We have already remarked that 
the kernel of $P\colon H^{k,p}_{\delta}(M;E)\to H^{k-m,p}_{\delta}(M;E)$
is equal to $Z$, which is finite-dimensional, and in fact is the
same as the $L^2$ kernel.
We will show that the range of $P$ is closed by showing that
it is equal to $Y^{k-m,p}_\delta$.
If $f=Pu\in P(H^{k,p}_{\delta}(M;E))$, then clearly $(f,v) =
(Pu,v) = (u,Pv) = 0$ for all $v\in Z$, so $f\in Y^{k-m,p}_\delta$.
On the other hand, if $f\in Y^{k-m,p}_\delta$, then 
$f = PGf + H f = PGf$ by \eqref{eq:PGu} and 
\eqref{eq:kerH}, so $f\in P(H^{k,p}_{\delta}(M;E))$.
Since every 
$f\in H^{k-m,p}_\delta(M;E)$ can be written $f=PGf+H f$,
where $PGf\in P(H^{k,p}_\delta(M;E)) = Y^{k-m,p}_\delta$
and $H f\in Z$, it follows that 
$H^{k-m,p}_\delta(M;E)  = Y^{k-m,p}_\delta \oplus Z$.
Therefore,
\begin{displaymath}
\frac{ H^{k-m,p}_\delta(M;E)} {P(H^{k,p}_{\delta}(M;E))}
= 
\frac{  Y^{k-m,p}_\delta \oplus Z} {Y^{k-m,p}_\delta }
\cong
Z,
\end{displaymath}
which is finite-dimensional.  This also shows that
the cokernel and kernel of $P$ have the same dimension,
so $P$ has index zero.

Next we will prove the necessity of the 
stated conditions on $\delta$.  
In fact, we will show that $P$ has infinite-dimensional kernel when $\delta$
is strictly below the Fredholm range, and infinite-dimensional cokernel
when $\delta$ is strictly above;
in the borderline case $C^{k,\alpha}_{n/2-R}$, we will also show that
$P$ has infinite-dimensional kernel, and in all other borderline
cases, we will show that it fails to have closed range.

First 
we address 
the H\"older case below the Fredholm range.
Assume $0<\alpha<1$, $m<k+\alpha\le l+\beta$, and 
$\delta\le n/2 -R$,
and consider $P\colon C^{k,\alpha}_{\delta}(M;E)\to C^{k-m,\alpha}_{\delta}(M;E)$.
We will prove that $P$ is not Fredholm in this case by
showing that it has an infinite-dimensional kernel.

The definition of the indicial radius $R$ and the symmetry
of the characteristic exponents about $\Real s = n/2-r$
imply that 
$P$ has a characteristic exponent $s_0$ with
$\Real s_0 = n/2-r-R$.
Let $W$ be the subspace of $C^{l,\beta}_{n/2-R}(M;E)$ 
given by Lemma \ref{lemma:critical-asymp-solutions} for this characteristic exponent.
(The compact set $K$ is irrelevant in this case.)
For any $p>n$, 
$P(W)\subset C^{0,0}_{n/2-R+1}(M;E) \subset H^{0,p}_{n/2-R}(M;E)$
by Lemma \ref{properties-of-spaces}.
If we choose $p$ large enough, then $\delta = n/2-R$
will satisfy $|\delta +n/p-n/2| = |n/p-R|<R$, so there exist
operators $G,H \colon H^{0,p}_{n/2-R}(M;E)\to H^{0,p}_{n/2-R}(M;E)$ satisfying
\eqref{eq:kerH}--\eqref{eq:PGu}.
Let $W_0\subset W$ be the linear subspace defined
by 
\begin{equation}
W_0 = \{w\in W: H Pw=0\}.
\end{equation}
Because $H $ takes its values in the finite-dimensional
space $Z$, the space $W_0$ is also infinite-dimensional.
Note that for $w\in W_0$, 
$Pw\in H^{0,p}_{n/2-R}(M;E)$ implies
$GPw \in H^{m,p}_{n/2-R}(M;E)\subset C^{0,0}_{n/2-R}(M;E)$
by Lemma \ref{properties-of-spaces}\eqref{prop:Sobolev}.
Define
$X\colon W_0\to C^{0,0}_{n/2-R}(M;E)$ by
$Xw = w - GPw$.  It follows from \eqref{eq:PGu}
that 
\begin{displaymath}
PXw = Pw - PGPw = H Pw = 0 \quad \text{for all $w\in W_0$}.
\end{displaymath}
Therefore,
\begin{displaymath}
X(W_0)\subset \Ker P\cap C^{0,0}_{n/2-R}(M;E)
\subset C^{l,\beta}_{n/2-R}(M;E)
\end{displaymath}
 by Lemma
\ref{rescaling-estimates}.
Moreover, $X$ is injective because $Xw=0$ implies
$w=GPw\in H^{0,p}_{n/2-R}(M;E)$, which implies that $w=0$ by
assertion \eqref{part:notin-Hp} of Lemma
\ref{lemma:critical-asymp-solutions}.
Thus we have shown that $X(W_0)$ is an infinite-dimensional
subspace of $\Ker P\cap C^{l,\beta}_{n/2-R}(M;E)$.
Since $C^{l,\beta}_{n/2-R}(M;E)\subset C^{k,\alpha}_{\delta}(M;E)$
whenever $\delta\le n/2 -R$ and 
$m<k+\alpha\le l+\beta$, it follows that 
$P$ has infinite-dimensional kernel on
$C^{k,\alpha}_{\delta}(M;E)$ in all such cases.

Next consider the Sobolev case 
below the Fredholm
range.  When $1<p<\infty$, $m\le k\le l$, and 
$\delta<n/2-n/p-R$, we have
$C^{l,\beta}_{n/2-R}(M;E) \subset H^{k,p}_{\delta}(M;E)$
by Lemma \ref{properties-of-spaces}, so
$P$ has infinite-dimensional kernel in $H^{k,p}_{\delta}(M;E)$ as well.

Now we consider the exponents strictly above the Fredholm range,
beginning with the 
Sobolev case.
Suppose $1<p<\infty$, $m<k+\alpha\le l+\beta$, and 
$\delta> n/2-n/p+R$.
Recall that the dual space to $H^{0,p}_{\delta}$ is 
$H^{0,p^*}_{-\delta}$ (where $p^*$ is the
conjugate exponent, $1/p+1/p^*=1$), acting by way of the standard $L^2$ pairing.  
Since $-\delta< n/2-n/p^*-R$, the argument above shows that
$P^*=P\colon H^{m,p^*}_{-\delta }(M;E)\to H^{0,p^*}_{-\delta }(M;E)$
has infinite-dimensional kernel.  
Each element $v$ of the infinite-dimensional
space $\Ker P\cap H^{m,p^*}_{-\delta}(M;E)$
thus defines a continuous linear functional
on $H^{k-m,p}_{\delta }(M;E)$ by $u\mapsto (u,v)$,
and each such linear functional 
annihilates 
$P(H^{k,p}_{\delta }(M;E))$
by Lemma 
\ref{lemma:P-symmetric}.
It follows that the range of $P$ has infinite
codimension in $H^{k-m,p}_{\delta }(M;E)$.

For the H\"older case, suppose 
$0<\alpha<1$, $m<k+\alpha\le l+\beta$, and 
$\delta>n/2+R$.
Choose $\delta'$ close to $\delta$
and $p$ sufficiently large 
that $n/2+R<\delta'+n/p<\delta$.  It follows that 
$-\delta'<n/2-n/p^*-R$, which implies 
as above that
$P$ has an infinite-dimensional kernel in 
$H^{m,p^*}_{-\delta'}(M;E)$.
The fact that
$C^{k-m,\alpha}_{\delta}(M;E)\subset H^{0,p}_{\delta'}(M;E)$
implies that 
$H^{0,p^*}_{-\delta'}(M;E) = 
(H^{0,p}_{\delta'}(M;E))^* 
\subset (C^{k-m,\alpha}_{\delta}(M;E))^*$.
As above, each
linear functional on $C^{k-m,\alpha}_{\delta}(M;E)$ defined by
an element of 
$\Ker P\cap H^{m,p^*}_{-\delta'}(M;E)$ annihilates the range of $P$,
so once again we conclude that 
$P\colon C^{k,\alpha}_{\delta}(M;E)
\to C^{k-m,\alpha}_{\delta}(M;E)$ has 
infinite-dimensional cokernel.

Next we consider the borderline cases.
The lower borderline H\"older case $\delta=n/2-R$
was already treated above, when we showed that $P$ has infinite-dimensional kernel
in $C^{k,\alpha}_{\delta}(M;E)$ whenever
$\delta\le n/2-R$.  
In all remaining cases, we will show that $P$ does not have
closed range.

We begin with the upper borderline Sobolev case, $H^{k,p}_{\delta_p}(M;E)$ with
$\delta_p = n/2-n/p+R$.
Let us assume that $P\colon H^{k,p}_{\delta_p}(M;E)\to H^{k-m,p}_{\delta_p}(M;E)$
has closed range and derive a contradiction.
If $\delta<\delta_p$ is chosen sufficiently close to $\delta_p$ that
$|\delta-n/2+n/p|<R$, the argument at the beginning of this proof showed
that $\Ker P\cap H^{k,p}_{\delta}(M;E)$ is finite-dimensional.
Since $H^{k,p}_{\delta_p}(M;E)\subset H^{k,p}_{\delta}(M;E)$,
we see that $\Ker P\cap H^{k,p}_{\delta_p}(M;E)$ is finite-dimensional as well,
so $P\colon H^{k,p}_{\delta_p}(M;E)\to H^{k-m,p}_{\delta_p}(M;E)$ is semi-Fredholm.
It follows from Lemma \ref{lemma:semi-fred} that there is a
compact set 
$K\subset M$ and a constant $C$ such that
\begin{equation}\label{eq:6.z}
\|u\|_{0,p,\delta_p}\le C\|Pu\|_{0,p,\delta_p}
\end{equation}
when $u\in H^{m,p}_{\delta_p}(M;E)$
is supported in $M\setminus K$.

By definition of $R$,
$P$ has a characteristic exponent $s_0$ whose real part is equal to
$n/2-r+R$. 
Let $w$ be any element of the space $W$ defined
in Lemma \ref{lemma:critical-asymp-solutions} corresponding to this
characteristic exponent, with 
$\supp w \subset M\setminus K$, 
$w\notin H^{0,p}_{\delta_p}(M;E)$, 
but $Pw\in C^{0,0}_{n/2+R+1}(M;E)\subset H^{0,p}_{\delta_p}(M;E)$.
Let $\{\psi_\epsilon\}$ be a family of cutoff functions
as in Lemma \ref{psi-epsilon}, and define $w_\epsilon = (1-\psi_\epsilon) w$.
Note that for any $\epsilon>0$, $w_\epsilon$
is in $C^{l,\beta}_{\text{loc}}(M;E)$ and 
compactly supported, so it is in $H^{k,p}_{\delta_p}(M;E)$
for all $k\le l$.
Because $w\notin H^{0,p}_{\delta_p}(M;E)$ 
and $w_\epsilon \to w$ uniformly on compact sets
as $\epsilon\to 0$,
we have $\|w_\epsilon\|_{0,p,\delta_p}\to \infty$.
On the other hand, 
\begin{equation}\label{eq:P-psi-e}
Pw_\epsilon = (1-\psi_\epsilon) Pw
- [P,\psi_\epsilon] w\in H^{0,p}_{\delta_p}(M;E).
\end{equation}
If we can show that $\|Pw_\epsilon\|_{0,p,\delta_p}$
remains bounded as $\epsilon\to 0$, we will have a
contradiction to \eqref{eq:6.z}.

The fact that  $Pw\in H^{0,p}_{\delta_p}(M;E)$ implies
that $(1-\psi_\epsilon) Pw \to Pw$ in the $H^{0,p}_{\delta_p}$ norm,
so the first term in \eqref{eq:P-psi-e}
is clearly bounded in $H^{0,p}_{\delta_p}$.
For the second term, observe that the
commutator $[\Del,\psi_\epsilon]w = w\otimes d\psi_\epsilon $ 
is an operator of order
zero with coefficients that are uniformly bounded in 
$C^{l,\beta}(M)$ and supported on the set 
where $\epsilon/2\le \rho \le \epsilon$.
It follows by induction that 
$[P,\psi_\epsilon]$ 
is an operator of order
$m-1$ with bounded coefficients supported in the same set,
and therefore
\begin{align*}
\|[P,\psi_\epsilon]w\|_{0,p,\delta_p}^p
&\le C \sum_{0\le j\le m-1} 
\int_{\epsilon/2\le \rho \le \epsilon}
|\Del^j w|^p\,dV_g.
\end{align*}
Since $|\Del^j w|^p$ is integrable for $0\le j\le m$,
each integral above goes to zero as $\epsilon\to 0$ by the dominated
convergence theorem.
This contradicts \eqref{eq:6.z} and completes the proof that $P$ does not
have closed range in this case.

Next consider the upper borderline H\"older case, $\delta = n/2+R$.
The argument is almost the same as in the Sobolev case,
except in this case we have to set $w_\epsilon = \psi_\epsilon w$ and show that
$\|Pw_\epsilon\|_{0,\alpha,n/2+R}\to 0$
while $\|w_\epsilon\|_{0,\alpha,n/2+R}$ remains bounded below by a positive
constant.  The details are left to the reader.

The only case left is the lower borderline Sobolev case, 
$H^{k,p}_{\delta}(M;E)$ with 
$\delta= n/2-n/p-R$.
Since $-\delta= n/2-n/p^*+R$, we showed above that
$P\colon H^{m,p^*}_{-\delta }(M;E)\to H^{0,p^*}_{-\delta }(M;E)$
does not have closed range, and thus neither does 
$P\colon H^{0,p^*}_{-\delta }(M;E)\to H^{0,p^*}_{-\delta }(M;E)$
considered as an unbounded operator (since its domain is exactly 
$H^{m,p^*}_{-\delta }(M;E)$).
Since a closed, densely defined operator has closed range if and only
if its adjoint does (cf.\ \cite[Theorem IV.5.13]{Kato}), 
this implies that 
$P\colon H^{0,p}_{\delta}(M;E)\to H^{0,p}_{\delta}(M;E)$
does not have closed range, and then it follows
from the regularity results of 
Proposition \ref{prop:optimal-regularity} that 
$P\colon H^{k,p}_{\delta}(M;E)\to H^{k-m,p}_{\delta}(M;E)$
does not either.

Finally, to prove part \eqref{part:pos-indicial}, 
just observe that $P$ being 
Fredholm on $L^2(M;E)$
implies that it is Fredholm from $H^{m,2}_0(M;E)$ to $H^{0,2}_0(M;E)$.
By part \eqref{sobolev-fredholm-result}, this 
in turn implies that $|0|<R$. 
\end{proof}

%%
%%
%% Fredholm operators and Einstein metrics
%% on conformally compact manifolds
%%
%% by John M. Lee
%% 
%% Chapter 7
%%

\chapter{Laplace Operators}\label{L2-section}

In this chapter we specialize to Laplace operators.
Throughout this chapter, $(M,g)$ will be a
connected asymptotically hyperbolic
$(n+1)$-manifold of class $C^{l,\beta}$ for some $l\ge 2$ and $0\le\beta< 1$.

Let $E$ be a geometric tensor bundle of weight $r$
over $M$.  
A \defname{Laplace operator} is a 
second-order geometric
operator $P\colon C^\infty(M;E)\to C^\infty(M;E)$
that can be written in the form
$P=\Del^*\Del+\scr K$, where $\Del^*\Del$ is the covariant Laplacian and
$\scr K\colon E\to E$ is a bundle endomorphism (i.e., a differential
operator of order zero).  Note that our definition of geometric
operators guarantees that the coefficients of $\scr K$ in any local
frame are contractions of tensor products of $g$,
$g^{-1}$, $dV_g$, and the curvature tensor.  

To get sharp Fredholm results for a specific Laplace operator,
we need to compute its indicial radius.
In general, this is
just a straightforward
computation in coordinates near the boundary.
Here are some of the results.

\begin{lemma}\label{lemma:cov-laplacian-ind-radius}
The covariant Laplacian $\Del^*\Del$ 
on trace-free symmetric $r$-tensors has indicial radius
\begin{displaymath}
R =  \sqrt{\frac{n^2}4+r}.
\end{displaymath}
\end{lemma}

We will postpone the proof of this lemma until
after Proposition 
\ref{prop:Is-vs-I0} below.
For the record, we also note the following,
which is proved in 
\cite{Mazzeo-thesis,Mazzeo-Hodge}.

\begin{lemma}[Mazzeo]\ 
\label{lemma:laplace-beltrami-ind-radius}
The Laplace-Beltrami operator $\Delta = 
dd^* + d^*d$ on $q$-forms
has the following indicial radius:
\begin{displaymath}
R = \begin{cases}
\dfrac n 2 - q, & 0\le q \le \dfrac n 2,\\
\dfrac 1 2\vphantom{\dfrac {\tw{\tw n}} 2}, & q = \dfrac{n+1}2,\\
q - \dfrac {n+2} 2, \vphantom{\dfrac {\tw{\tw n}} 2}& \dfrac {n+2} 2 \le q \le n+1.
\end{cases}
\end{displaymath}
\end{lemma}

The following lemma greatly simplifies the computation of
the indicial radius of a Laplace operator.

\begin{proposition}\label{prop:Is-vs-I0}
Let $P = \Del^*\Del + \scr K$ be a Laplace operator
acting on a geometric tensor bundle of weight $r$.
For any $s\in\C$, 
\begin{displaymath}
I_s(P) = I_0(P) + s(n-s-2r).
\end{displaymath}
\end{proposition}

\begin{proof}
Since $I_s(\scr K) = \scr K|_{\del M}$ is independent
of $s$, we need only consider the case $P=\Del^*\Del$.
As in the
proof of Proposition 2.7 of
\cite{Graham-Lee}, we compute (in background coordinates)
\begin{equation}\label{eq:asymptotics-of-delrho}
\begin{aligned}
\rho_{;ij} &= \partial_i\partial_j \rho -
  \Gamma_{ij}^k\partial_k\rho\\
&= \rho^{-1}(2\del_i \rho\del_j \rho -  \gbar_{ij}) + o(\rho^{-1});\\
\Delta(\rho^s)
&= -s(s-1)\rho^{s-2}g^{jk}\rho_{;j}\rho_{;k} - s\rho^{s-1}\rho_{;k}{}^k\\
&= s(n-s)\rho^s + o(\rho^{s});\\
|\Del\rho|_g^2 
&= |d\rho|_g^2 \\
&= \rho^2|d\rho|_{\overline g}^2 =\rho^2+o(\rho^2).
\end{aligned}
\end{equation}
Therefore,
\begin{equation}\label{eq:rho-s-P-rhos}
\begin{aligned}
\rho^{-s} \Del^*\Del (\rho^s\overline u)
&= - \rho^{-s} \Tr_g \Del^2 (\rho^s\overline u)\\
&= - \rho^{-s} \Tr_g \Del ( s\rho^{s-1} \overline u \otimes \Del\rho
+ \rho^s \Del\overline u)\\
&= - \rho^{-s} \Tr_g \big( s(s-1)\rho^{s-2} \overline u \otimes \Del\rho \otimes \Del\rho
+2s\rho^{s-1} \Del\overline u \otimes \Del\rho\\
&\qquad +s\rho^{s-1} \overline u \otimes \Del^2\rho
+ \rho^s \Del^2\overline u\big)\\
&= -s(s-1)\rho^{-2}|\Del\rho|_g^2 \overline u
-2s\rho^{-1} \Del_{\grad \rho}\overline u\\
&\qquad+s\rho^{-1}\Delta\rho \overline u
+ \Del^*\Del\overline u\\
&= -s(s-1) \overline u
-2s\rho^{-1} \Del_{\grad \rho}\overline u
+s(n-1)\overline u
+ \Del^*\Del\overline u+o(1).
\end{aligned}
\end{equation}
To compute the second term above,
assume $\overline u$ is a tensor of type $\binom q p$
with $q-p=r$, and
let $D$ be the difference tensor $D=\Del-\overline\Del$.  The 
components of $\Del\overline u$ are 
\begin{align*}
\overline u^{i_1\dots i_p}
_{j_1\dots j_q k}
&= \del_k \overline u^{i_1\dots i_p}
_{j_1\dots j_q}
+ \sum_{s=1}^p \Gamma _{mk}^{i_s}
\overline u^{i_1\dots m \dots i_p}
_{j_1\dots j_q}
- \sum_{s=1}^q
\Gamma _{j_sk}^{m}
\overline u^{i_1\dots i_p}
_{j_1\dots m \dots j_q}\\
&=O(1)+ \sum_{s=1}^p  D_{mk}^{i_s}
\overline u^{i_1\dots m \dots i_p}
_{j_1\dots j_q}
- \sum_{s=1}^q  D_{j_sk}^{m}
\overline u^{i_1\dots i_p}
_{j_1\dots m \dots j_q}.
\end{align*}
Let us introduce the shorthand notations $\rho_i = \del_i\rho$
and $\overline \rho^i = \overline g^{ij}\del_j\rho$.
Using 
formula \eqref{eq:difference-tensor-components}
for the components of $D$,
together with the fact that $\rho_i\overline\rho^i=1+O(\rho)$,
we obtain
\begin{align*}
\rho^{-1}g^{kl}\rho_l D_{jk}^{i}
&= - \overline \rho^{k}(\delta^i_j \rho_k + \delta^i_k \rho_j
- \overline g_{jk} \overline \rho^{i} )
= -\delta^i_j + O(\rho).
\end{align*}
Therefore, 
\begin{align*}
\rho^{-1}(\Del_{\grad\rho}\overline u)^{i_1\dots i_p}
_{j_1\dots j_q}
&= \rho^{-1} g^{kl}\rho_l
\overline u^{i_1\dots i_p}
_{j_1\dots j_q;k}\\
&= 
-\sum_{s=1}^p  \delta^{i_s}_m
\overline u^{i_1\dots m \dots i_p}
_{j_1\dots j_q}
+ \sum_{s=1}^q  \delta^{m}_{j_s}
\overline u^{i_1\dots i_p}
_{j_1\dots m \dots j_q} + O(\rho)
\\
&= r\overline u^{i_1\dots i_p}
_{j_1\dots j_q} + O(\rho).
\end{align*}
Inserting this back into 
\eqref{eq:rho-s-P-rhos},
we obtain
\begin{align*}
\rho^{-s} \Del^*\Del (\rho^s\overline u)
= s(n-s-2r)\overline u
+ \Del^*\Del\overline u+o(1),
\end{align*}
which implies the result.
\end{proof}

\begin{proof}[Proof of Lemma \ref{lemma:cov-laplacian-ind-radius}]
By Prop.\ \ref{prop:Is-vs-I0}, to determine $I_s(\Del^*\Del)$ it suffices
to compute $I_0(\Del^*\Del)$.  
The cases $r=0$ and $r=2$ follow from Corollary 2.8 and 
Lemma 2.9 of \cite{Graham-Lee}.  For the general case, 
suppose $r\ge 1$.
Assuming $\overline u$  
is trace-free, symmetric, and
smooth up to the boundary, and using the notation of the preceding
proof, we compute
\begin{align*}
(\Del^*\Del \overline u)_{i_1\dots i_r}
&= - \overline u_{i_1\dots i_r;l}{}^l\\
&= - g^{lm} \left( \del _m \overline u_{i_1\dots i_r;l} - D^j_{lm} \overline u_{i_1\dots i_r;j}
- \sum_{s=1}^r D^j_{i_s m} \overline u_{i_1\dots j \dots i_r;l}
\right) + O(\rho).
\end{align*}
As in the preceding proof, we have 
\begin{displaymath}
\overline u_{i_1 \dots i_r;l} = 
- \sum_{t=1}^r D^k_{i_t l} \overline u_{i_1\dots k \dots i_r} + O(1),
\end{displaymath}
and therefore,
\begin{align*}
(\Del^*\Del \overline u)_{i_1\dots i_r}
&=  
\sum_{t=1}^r g^{lm} (\del _m D^k_{i_t l}) \overline u_{i_1\dots k \dots i_r}
- \sum_{t=1}^r g^{lm} D^j_{lm} D^k_{i_t j} \overline u_{i_1\dots k \dots i_r}\\
&\quad 
- \sum_{\substack{
1\le s,t\le r\\
s\ne t}}
g^{lm} D^j_{i_s m} 
D^k_{i_t l} \overline u_{i_1\dots j\dots k \dots i_r}
- \sum_{s=1}^r 
g^{lm} D^j_{i_s m} D^k_{jl} \overline u_{i_1\dots k \dots i_r}
\\
&\quad + O(\rho).
\end{align*}
Using 
\eqref{eq:difference-tensor-components}
again, we compute
\begin{align*}
g^{lm} \del_m D^k_{i_tl}
&= g^{lm} \rho^{-2} \rho_m\left( 
\delta^k_{i_t}\rho_l + \delta^k_l\rho_{i_t} - \overline g_{i_tl}\overline \rho^k
\right)+O(\rho)= \delta^k_{i_t}+O(\rho);\\
g^{lm} D^j_{lm} D^k_{i_t j}
&= -(n-1) \delta^k_{i_t}+O(\rho);\\
g^{lm} D^j_{i_sm} D^k_{{i_t}l}
&= \delta^j_{i_s}\delta^k_{i_t} +
\overline g^{jk} \rho_{i_s} \rho_{i_t}
 - \delta_{i_s}^k \rho_{i_t}\overline \rho^j
 - \delta^j_{i_t}\rho_{i_s} \overline \rho^k + \overline g_{{i_s}{i_t}}\overline \rho^j\overline \rho^k
  +O(\rho);\\
g^{lm} D^j_{i_sm} D^k_{jl}
&= -(n-1) \rho_{i_s} \overline \rho^k+O(\rho).
\end{align*}
Inserting these above and using the fact that $\overline u$ is symmetric and trace-free, 
we obtain
\begin{align*}
(\Del^*\Del \overline u)_{i_1\dots i_r}
&=  
r \overline u_{i_1\dots i_r}
+ r(n-1) \overline u_{i_1\dots i_r}
- r(r-1)\overline u_{i_1\dots i_r}-0\\
&\quad 
+ (r-1) \sum_{t=1}^r \rho_{i_t} \overline\rho^j \overline u_{i_1\dots j \dots i_r}
+ (r-1) \sum_{s=1}^r \rho_{i_s} \overline\rho^k \overline u_{i_1\dots k \dots i_r}\\
&\quad
- \sum_{\substack{
1\le s,t\le r\\
s\ne t}}
\overline g_{i_si_t} \overline \rho^j\overline \rho^k \overline u_{i_1\dots j \dots k \dots i_r}
+ (n-1)\sum_{s=1}^r \rho_{i_s}\overline\rho^k \overline u_{i_1\dots k\dots i_r}
+O(\rho)\\
&= r(n+1-r)\overline u_{i_1\dots i_r} 
+ (2r+n-3) \sum_{t=1}^r \rho_{i_t} \overline\rho^j \overline u_{i_1\dots j \dots i_r}\\
&\quad 
- \sum_{\substack{
1\le s,t\le r\\
s\ne t}}
\overline g_{i_si_t} \overline \rho^j\overline \rho^k \overline u_{i_1\dots j \dots k \dots i_r}
+O(\rho).
\end{align*}

Suppose 
$s$ is any indicial root of $\Del^*\Del $ and $\overline u$ is a corresponding 
unit eigentensor.  Using the formula above for $\Del^*\Del \overline u$
together with Proposition
\ref{prop:Is-vs-I0}, and observing that $\overline u$ is trace-free, we compute
\begin{align*}
0 &= \left\< \overline u, I_{s}(\Del^*\Del) \overline u\right\>_{\overline g}\\
&= r(n+1-r)\left|\overline u\right|_{\overline g}^2 
+ (2r+n-3) \left|\grad_{\overline g}\rho\into\overline u\right|_{\overline g}^2
+ s(n-s-2r)\left|\overline u\right|_{\overline g}^2\\
&\ge r(n+1-r) + s(n-s-2r),
\end{align*}
which implies that each indicial root satisfies 
\begin{displaymath}
\left|s - \left(\frac n 2 -r\right)\right|^2\ge \frac{n^2}{4}+r.
\end{displaymath}
It follows that the indicial radius of $\Del^*\Del $ is at least $\sqrt{ n^2/4+r }$.

On the other hand, if $\overline u$ 
is chosen so that $\grad_{\overline g}\rho\into\overline u=0$ along $\del M$ (i.e., $\overline u$ is
purely tangential),
we find that 
$I_s(\Del^*\Del )\overline u = \bigl(r(n+1-r) + s(n-s-2r)\bigr)\overline u$.
Solving $I_s(\Del^*\Del )\overline u=0$ for $s$, we find that 
two of the indicial roots of $\Del^*\Del $ are 
\begin{displaymath}
s = \frac n 2 - r \pm \sqrt{ \frac {n^2}{4}+r}.
\end{displaymath}
This proves that the indicial radius is exactly $\sqrt{ n^2/4+r }$ as claimed.
\end{proof}

\begin{corollary}\label{cor:P+lambda}
Let $P = \Del^*\Del + \scr K$ be a Laplace operator
acting on a geometric tensor bundle of weight $r$,
and suppose $P$ has indicial radius $R>0$.  
For $c\in \R$, the indicial radius $R'$
of $P+c$ is positive if and only if 
$c+R^2>0$, in which case 
\begin{displaymath}
R' = \sqrt{c+R^2}.
\end{displaymath}
\end{corollary}

\begin{proof}
Comparing the formulas given by 
Proposition \ref{prop:Is-vs-I0} for $I_s(P)$
and $I_{n/2-r}(P)$, we find that
\begin{displaymath}
I_s(P) = I_{n/2-r}(P) - (s-n/2+r)^2.
\end{displaymath}
Observe that $I_{n/2-r}(P)$ is self-adjoint by 
Proposition \ref{prop:Is*}, so it has
real eigenvalues.
Now $s$ is a characteristic root of $P$ precisely
when   
$s$ is a solution to the quadratic equation
$(s-n/2+r)^2 =\mu$
for some eigenvalue
$\mu$ of $I_{n/2-r}(P)$.
If some eigenvalue were nonpositive,
this equation would have a root with real part
equal to $n/2-r$, which would imply $R=0$.
Therefore, the assumption $R>0$ means that
all of the eigenvalues $\{\mu_i\}$ of $I_{n/2-r}(P)$ are
strictly positive, and therefore the characteristic roots of $P$
are $s=n/2-r\pm\sqrt{\mu_i}$, with $R^2= \min \{\mu_i\}$.

Since $I_s(P+c) = I_s(P)+c$,
the characteristic roots of $P+c$ are
$s=n/2-r\pm\sqrt{c+\mu_i}$,
and the one with smallest real part greater than
$n/2-r$ is $s=n/2-r\pm\sqrt{c+R^2}$.
Thus the indicial radius of $P+c$ is
$\sqrt{c+R^2}$ as claimed.
\end{proof}

\begin{lemma} 
Let $\DeltaL$
be the Lichnerowicz Laplacian, and let
$c$ be a real constant.  If $n^2/4-2n+c> 0$,
then the indicial radius of 
$\DeltaL + c$
acting on symmetric $2$-tensors is 
\begin{equation}\label{eq:ind-radius-deltaL}
R =  \sqrt{ \frac{n^2}4 - 2 n+c}.
\end{equation}
\end{lemma}

\begin{proof}
Observe first that
$\DeltaL$ preserves the splitting
of symmetric $2$-tensors into trace and trace-free parts:
\begin{displaymath}
\Sigma^2M = \R g \oplus\Sigma^2_0M ,
\end{displaymath}
where $\Sigma^2M$ is the bundle of symmetric covariant
$2$-tensors,
$\Sigma^2_0M$ is the subbundle of tensors that are
trace-free with respect to $g$, and 
$\R g\subset \Sigma^2M$ is the 
real line bundle of multiples
of $g$.
On $\R g$, 
$\mathring{Rc}(u g) = \mathring{Rm}(u g)$, so 
$\DeltaL$ acts as the ordinary
Laplacian:
\begin{displaymath}
\DeltaL(u g) = (\Del^*\Del u)g.
\end{displaymath}
It follows from 
Lemma \ref{lemma:cov-laplacian-ind-radius}
(or Lemma \ref{lemma:laplace-beltrami-ind-radius})
that the indicial radius of $\Del^*\Del$ on functions
is $R = n/2$, so the indicial radius of $\Del^*\Del + c$
is $\sqrt{n^2/4 + c}$, which is greater than
\eqref{eq:ind-radius-deltaL}.  Therefore, it suffices to show
that $\DeltaL+c$ acting on trace-free symmetric
$2$-tensors has indicial radius given by 
\eqref{eq:ind-radius-deltaL}.

The asymptotic formula
\eqref{eq:Rijkl} for the Riemann curvature tensor
implies that the action of 
$\mathring{Rm}$ and $\mathring{Rc}$
on $\Sigma^2_0M$
near the boundary is given by 
\begin{equation}\label{eq:est-for-Rm-Rc}
\begin{aligned}
\mathring{Rc}(u) &= -nu + O(\rho|u|);\\
\mathring{Rm}(u) &= u + O(\rho|u|).
\end{aligned}
\end{equation}
Thus the indicial radius
of $\DeltaL+c$ on $\Sigma^2_0M$
is the same as that of $\Del^*\Del -2n-2+c$,
which is 
$\sqrt{n^2/4 -2n+c}$ by
Lemma \ref{lemma:cov-laplacian-ind-radius}
and Corollary \ref{cor:P+lambda}.
\end{proof}

The main thing that needs to be checked in order
to apply Theorem \ref{thm:main-fredholm} 
is the $L^2$ estimate 
\eqref{eq:asymptotic-L2-estimate}.
For some operators, an appropriate asymptotic estimate 
follows from 
an obvious integration by parts, 
such as $\Del^*\Del+c$ when $c$ is a positive constant:
\begin{align*}
(u,\Del^*\Del u+c u)
= \|\Del u\|^2 + c \|u\|^2
\ge c \|u\|^2,
\end{align*}
from which $\|u\| \le c^{-1}\|(\Del^*\Del+c)u\|$
follows by the Cauchy-Schwartz inequality.
However, when the zero-order term is not strictly 
positive, we need
to work a bit harder.

As a warmup for the general $L^2$ estimates we will
prove below, consider first the ordinary
Laplacian $\Delta=d^*d$ on functions.  The use of positive eigenfunctions, and
more generally positive functions
satisfying differential inequalities, is a common tool for estimating the
lower bound of the spectrum of elliptic operators; see for example
\cite{CY,Lee,McKean} and especially \cite{Sullivan}, where such functions
play a central role.  The following
lemma was proved originally by Cheng and Yau
\cite[p.\ 345]{CY}.

\begin{lemma}[Cheng-Yau]
\label{Cheng-Yau}
Let $M$ be any Riemannian manifold.  If there exists a positive, locally
$C^2$ function $\phi$ on $M$ such that $\Delta\phi/\phi\ge\lambda$, then
\begin{equation}\label{CY-estimate}
(u,\Delta u)\ge \lambda\|u\|^2
\end{equation}
for all smooth compactly supported functions $u$.
\end{lemma}

The proof of this lemma in \cite{CY} uses the maximum principle.  Here
is a simple proof based on integration by parts, which serves to
motivate the somewhat more delicate estimates below.

\begin{proof}
Let $u \in C^\infty_c(M)$.
The divergence theorem gives
\begin{align*}
0 &= \int_M d^*(u^2\phi^{-1}d\phi)\,dV_g  \\
&= \int_M
\left( -2u \phi^{-1}\left<du,d\phi\right>_g+u^2\phi^{-2}|d\phi|^2
+ u^2\phi^{-1}\Delta\phi \right) \, dV_g.
\end{align*}
Thus
\begin{align*}
0 &\le \int_M
\left|\phi d(\phi^{-1}{u})\right|^2 \,dV_g \\
&= \int_M \left| du - u \phi^{-1}d\phi\right|^2\,dV_g\\
&= \int_M \left(| du|^2 - 2u \phi^{-1}\left<du,d\phi\right>_g+
u^2 \phi^{-2}|d\phi|^2\right)\,dV_g\\
&= \int_M \left( u\Delta u - u^2 \phi^{-1}{\Delta \phi}
\right)\, dV_g  \le \int_M \left( u\Delta u - \lambda u^2
\right)\, dV_g ,
\end{align*}
which is equivalent to \eqref{CY-estimate}
\end{proof}

Lemma \ref{Cheng-Yau} is a global result; but 
for proving asymptotic estimates, we need only 
find a function
$\phi$ that satisfies $\Delta\phi/\phi\ge \lambda$ 
on the complement of a compact set.
In particular, on an asymptotically hyperbolic $(n+1)$-manifold,
asymptotic computations show that $\Delta(\rho^{n/2})/\rho^{n/2}$ can be
made arbitrarily close to $n^2/4$ near $\dm$.  
From this it will follow, for example, that 
$\Delta-\lambda$ on functions
satisfies the hypotheses of Theorem \ref{thm:main-fredholm}
for any constant
$\lambda<n^2/4$.

It is interesting to note that this simple estimate immediately yields a
(rather crude) estimate for the covariant Laplacian on tensor fields.

\begin{lemma}\label{lemma:crude-tensor-est}
If \eqref{CY-estimate} holds for smooth functions compactly supported in
some open set $U\subset M$, then for any smooth tensor field $w$ compactly
supported in $U$, we have
\begin{displaymath}
(w,\Del^*\Del w) \ge \lambda\|w\|^2
\end{displaymath}
with the same constant $\lambda$.
\end{lemma}

\begin{proof}
Inequality
\eqref{CY-estimate} implies
\begin{displaymath}
\lambda\|u\|^2 \le \|\Del u\|^2,
\end{displaymath}
which extends continuously to compactly supported functions in
$H^{1,2}(U)$.  If $w$ is a smooth, compactly supported tensor field,
the function $|w|$ is Lipschitz, hence in $H^{1,2}(U)$.  Kato's
inequality says that 
$|\Del|w|\,|\le |\Del w|$ almost everywhere (see \cite[Prop.\
3.49]{Aubin}, where this is proved for scalar functions; the proof
extends easily to tensors).  Therefore
\begin{displaymath}
\lambda\|w\| ^2
\le \|\Del|w|\,\|^2\le \|\Del w\| ^2= (w,\Del^*\Del w).
\end{displaymath}
\end{proof}

A version of this argument (using the maximum principle instead of
integration by parts) was used implicitly in our proof that the
linearized Einstein operator is invertible on weighted H\"older spaces
over hyperbolic space \cite{Graham-Lee}.  
Unfortunately, as we noted there, this estimate was not sharp,
and led to less-than-optimal Fredholm results for tensors.  
For the application to Einstein metrics in this monograph,
we no longer need a sharp asymptotic estimate, because of the
sharp Fredholm theorems of the preceding chapter.  However,
with other applications in mind, it is useful to see how far the
asymptotic estimates can be pushed.

The key to finding 
improved asymptotic estimates on tensors 
turns out to be to consider $r$-tensor fields as
$(r-1)$-tensor-valued 1-forms.  Therefore we must make a short
digression to discuss the properties of tensor-valued differential
forms.

Let $E$ be any geometric 
tensor bundle over $M$, and let $\Lambda^qE:=E \tprod
\Lambda^qM$ denote the bundle of $E$-valued
$q$-forms on $M$.  Let $D\colon C^\infty(M;\Lambda^qE)\to C^\infty(M;
\Lambda^{q+1}E)$ denote the exterior covariant differential on
$E$-valued forms, defined by
\begin{displaymath}
D(\sigma\tprod\alpha) = \Del\sigma \wedge \alpha + \sigma \tprod
d\alpha,
\end{displaymath}
for $\alpha\in C^\infty(M;\Lambda^qM)$ and $\sigma\in C^\infty(M;E)$.
(The wedge product above is computed by wedging the 1-form component of
$\Del\sigma$ with $\alpha$ to yield a section of $\Lambda^{q+1}E$.)  We
will study the covariant Laplace-Beltrami operator on $E$-valued forms,
defined by $\Delta = DD^* + D^*D$, where $D^*$ is the formal adjoint of
$D$.

For a scalar 1-form $\alpha$, we let $\alpha\vee\colon \Lambda^{q+1}E \to
\Lambda^{q}E $ denote the (pointwise) adjoint of the operator
$\alpha\wedge\colon \Lambda^qE \to \Lambda^{q+1}E$
with respect to $g$, so that $\left<
\alpha\wedge \omega, \eta \right>_g = \left< \omega,\alpha \vee
\eta \right>_g$.  In particular, if $\beta$ is also a scalar
$1$-form, then $\alpha\vee\beta=\<\alpha,\beta\>_g$.

For any function $u\in C^2(M)$, let $H(u)$
denote the covariant Hessian of $u$ acting as 
bundle
endomorphism $H(u)\colon \Lambda^1M\to \Lambda^1M$, and 
extended to $\Lambda^qM$ as a derivation.  In terms
of any orthonormal basis,
\begin{equation}\label{formula-for-hessian}
H(u)\omega = u_{;ij}
e^i\wedge(e^j\vee\omega)
\end{equation}
(where $u_{;ij}$ are the components of $\Del^2u$), since both sides are
derivations that agree on $\Lambda^1M$.  
We extend this endomorphism to $\Lambda^qE$ by letting it act on the
differential form component alone.

\begin{lemma}\label{properties-of-D}
Suppose $\alpha$ and $\beta$ are scalar $1$-forms, $\omega$ is an
$E$-valued $q$-form, and $u$ is a function.  For any local orthonormal
frame $\{e_j\}$ for $TM$ and dual coframe $\{e^j\}$, we have the
following facts:
\begin{enumerate}\letters
\item \label{formula-for-D}
$D\omega = e^j \wedge \Del_{e_j} \omega$.
\item $D^*\omega = - e^j\vee \Del_{e_j} \omega$.
\item $D(u\omega) = u D\omega + du \wedge \omega$.
\item $D^*(u\omega) = u D^*\omega - du\vee \omega$.
\item \label{distribute-wedge}
$\alpha\wedge(\beta\vee\omega) +
\beta\vee(\alpha\wedge\omega) = \left<\alpha,\beta\right>_g\omega$.
\item \label{distribute-D}
$D(du\vee\omega) = - du\vee D\omega + H(u)\omega +
\Del_{\grad u} \omega$.
\item \label{distribute-D*}
$D^*(du\wedge\omega) = -du\wedge D^*\omega + H(u) \omega -
\Del_{\grad u}\omega + (\Delta u)\omega$.
\end{enumerate}
\end{lemma}

\begin{proof}
Parts \eqref{formula-for-D} through \eqref{distribute-wedge} are
standard, and can be found, for example, in \cite[Ch.\ 2 and Section 6.1]{Wu}.
For \eqref{distribute-D},
choose a point $p\in M$ and a
frame $\{e_j\}$ such that $\Del e_j =
\Del e^j = 0$ at $p$.  Then, computing at $p$ and using
\eqref{formula-for-D}, \eqref{distribute-wedge},
and \eqref{formula-for-hessian}, we have
\begin{align*}
D(du\vee\omega) &= e^j\wedge\Del_j(u_{;k}e^k\vee\omega)\\
 &= e^j\wedge(u_{;jk} e^k \vee\omega) +
    u_{;k} e^j\wedge(e^k\vee\Del_j\omega)\\
 &= H(u)\omega + u_{;k} g^{jk}\Del_j\omega
    - u_{;k} e^k\vee(e^j\wedge\Del_j\omega)\\
 &= H(u)\omega + \Del_{\grad u}\omega - du\vee D\omega.
\end{align*}
The computation for \eqref{distribute-D*} is similar.
\end{proof}

The integral formula in the next lemma, a tensor analogue of
\eqref{CY-estimate}, is the key to proving sharp asymptotic $L^2$ 
estimates.  It unifies and generalizes Bochner-type formulas that have
been used in various settings, such as the weighted Bochner formula
introduced by Witten \cite{Witten} to prove the Morse inequalities, and
a closely related formula for scalar differential forms used by Donnelly
and Xavier \cite{DX} to analyze the spectrum of the scalar
Laplace-Beltrami operator on negatively-curved manifolds, and by
Lars Andersson \cite{Andersson} for the same operator on asymptotically
hyperbolic manifolds.  Later Andersson and Chru{\'s}ciel \cite{AC} used the
formula presented here (based on an early draft of the present monograph) to
obtain Fredholm results for the ``vector Laplacian'' $L^*L$
that arises in the
constraint equations of general relativity (see the Introduction). 

\begin{lemma}\label{basic-L2-estimate}
For any smooth, compactly supported section $\omega$ of $\Lambda^q E$,
and any positive $C^2$ function $\phi$ on $M$, the following integral
formula holds:
\begin{equation}\label{integral-formula}
\begin{aligned}
( \omega, \Delta \omega) &= \int_M
\phi^{-1}\Delta\phi |\omega|_g^2 +
2 \left< \omega, H({\log \phi})\omega\right>_g\\
&\qquad+ |\phi D(\phi^{-1}\omega)|_g^2 + |\phi^{-1}
D^*(\phi\omega)|_g^2\,dV_g\\
& \ge \int_M\left<\omega,(\phi^{-1}\Delta\phi  +
2 H({\log \phi}))\omega\right>_g\,dV_g.
\end{aligned}
\end{equation}
\end{lemma}

\begin{proof}
Let $u=\log \phi$.
Using Lemma \ref{properties-of-D}, we compute
\begin{align*}
\int_M &|e^u D(e^{-u}\omega)|_g^2 + |e^{-u}D^*(e^u\omega)|_g^2
    \,dV_g  \\
 &= \int_M |D\omega - du\wedge\omega|_g^2
    + |D^*\omega - du\vee\omega|_g^2\,dV_g \displaybreak[0]\\
 &= \int_M |D\omega|_g^2 - 2\left<D\omega,du\wedge\omega\right>_g +
    |du\wedge\omega|_g^2 \displaybreak[0]\\
 &\qquad + |D^*\omega|_g^2 - 2\left<D^*\omega,du\vee\omega\right>_g +
    |du\vee\omega|_g^2 \,dV_g \displaybreak[0]\\
 &= \int_M \left<\omega,\Delta\omega\right>_g - 2\left<du \vee D\omega,\omega\right>_g
    + \left<du\vee(du\wedge\omega),\omega\right>_g\\
 &\qquad - 2\left<\omega,D(du\vee\omega)\right>_g
    + \left<du\wedge(du\vee\omega),\omega\right>_g\,dV_g \\
 &= \int_M \left<\omega,\Delta\omega\right>_g + |du|_g^2 |\omega|_g^2
    - 2\left<\omega,H(u)\omega\right>_g
    - 2\left<\omega,\Del_{\grad u}\omega\right>_g\,dV_g \\
 &= \int_M \left<\omega,\Delta\omega\right>_g + |du|_g^2 |\omega|_g^2
    - 2\left<\omega,H(u)\omega\right>_g
    - (\Delta u)|\omega|_g^2\,dV_g \\
 &= \int_M \left<\omega,\Delta\omega\right>_g - e^{-u}\Delta(e^u)|\omega|_g^2
    - 2\left<\omega,H(u)\omega\right>_g\,dV_g .
\end{align*}
\end{proof}

To make use of this formula, we will use a power of the defining
function $\rho$ as our weight function $\phi$.  
It is convenient
to introduce the following notation: If 
$P\colon C^\infty(M;E)\to C^\infty(M;E)$ 
is a differential operator on $M$ and $\lambda$ is a real
number, 
we write
\begin{equation}\label{eq:asymptotic-est}
(u,Pu)\gtrsim \lambda\|u\|^2
\end{equation}
to mean that  for every $\epsilon>0$, there exists 
a compact set $K_\epsilon$ such that 
\begin{displaymath}
(u,Pu)\ge (\lambda-\epsilon)\|u\|^2
\end{displaymath}
whenever $u$ is smooth and compactly supported in 
$M\setminus K_\epsilon$.

\begin{lemma}\label{asymptotic-estimate-lemma}
Let $M$ be an  asymptotically hyperbolic $(n+1)$-manifold of class
$C^{l,\beta}$, with $l\ge 2$ and $0\le \beta<1$,
and suppose either $0 \le q < n/2$ or $(n+2)/2 < q \le n$.  The covariant
Laplace-Beltrami operator satisfies the following asymptotic estimate on
$E$-valued $q$-forms:
\begin{displaymath}
(\omega,\Delta \omega)\gtrsim \lambda\|\omega\|^2,
\end{displaymath}
where $\lambda = (n-2q)^2/4$ if $0 \le q < n/2$ and $\lambda =
(n+2-2q)^2/4$ if
$(n+2)/2 < q \le n$.
\end{lemma}

\begin{proof}
We will use \eqref{integral-formula} for forms $\omega$ supported near
the boundary.  Let $\rho$ be a smooth defining function.  
Using \eqref{eq:asymptotics-of-delrho}, we compute
\begin{align*}
(\log\rho)_{;ij}
=\rho^{-2}\rho_{;i}\rho_{;j} - g_{ij} + O(\rho).
\end{align*}
The tensor $g$ 
acts as the identity on
$1$-forms, and since we extend it to act as a derivation, it 
acts on $E$-valued $q$-forms as $q$ times the identity.
Using \eqref{formula-for-hessian}, therefore, the action of
$H(\log\rho)$ on an $E$-valued $q$-form $\omega$ can be written
\begin{align*}
\left<\omega,
H(\log\rho)\omega\right>_g
&= \left<\omega,(\log\rho)_{;ij}e^i\wedge(e^j\vee \omega)\right>_g\\
&= \left<\omega,\frac{d\rho}{\rho}\wedge\left( \frac{d\rho}{\rho}\vee
\omega\right)\right>_g
- q|\omega|_g^2 +
O(\rho |\omega|_g^2)\\
&= \left|\frac{d\rho}{\rho}\vee\omega \right|_g^2 - q|\omega|_g^2 +
O(\rho |\omega|_g^2).
\end{align*}
Thus with $\phi=\rho^s$, where $s$ is a constant to be determined later,
the integrand on the right-hand side of \eqref{integral-formula} can be
estimated as follows:
\begin{equation}\label{eq:Hlogrho-est}
\begin{aligned}
\langle\omega,&(\rho^{-s} \Delta(\rho^s)  +
2s H({\log \rho}))\omega\rangle \\
&\ge s(n-s-2q)|\omega|_g^2 + 2s \left|\frac{d\rho}{\rho} \vee \omega\right|_g^2
- O(\rho|\omega|_g^2).
\end{aligned}
\end{equation}

Given $\delta>0$, we can choose $\epsilon$ small enough that the
absolute value of the $O(\rho|\omega|_g^2)$ factor above
is bounded by $\delta|\omega|_g^2$ on
$A_\epsilon$.  If $s\ge 0$, we then have for $\omega$ compactly
supported in $A_\epsilon$
\begin{displaymath}
(\omega,\Delta\omega)
\ge (s(n-s-2q)-\delta) \|\omega\|_g^2.
\end{displaymath}
This estimate is optimal when $s=(n-2q)/2$, so as long as $(n-2q)/2>0$,
which is to say $q<n/2$, we obtain the conclusion of the lemma
with $\lambda= (n-2q)^2/4$.

If on the other hand $s\le0$, we use the fact that $|d\rho/\rho|_g$
approaches $1$ uniformly at $\dm$.  We choose $\epsilon$ small enough
that the $O(\rho|\omega|_g^2)$ factor in 
\eqref{eq:Hlogrho-est}
is bounded by $(\delta/2)|\omega|_g^2$ and
$2s|d\rho/\rho|_g^2 \ge 2s-\delta/2$ on $A_\epsilon$, and conclude that
\begin{align*}
(\omega,\Delta\omega)
&\ge (s(n-s-2q)-\delta/2)\|\omega\|^2
+ (2s-\delta/2) \|\omega\|^2 \\
&= (s(n-s-2q+2)-\delta)\|\omega\|^2.
\end{align*}
This in turn is optimal when $s=(n+2-2q)/2$.  Thus as long as
$q>(n+2)/2$, so that $s<0$, we obtain the conclusion of the lemma
with $\lambda=(n+2-2q)^2/4$.
\end{proof}

When applied to scalar-valued forms, this result can be used to obtain
an elementary proof of the sharp
$L^2$ Fredholm theorems and spectral bounds
for the Laplace-Beltrami operator originally obtained by Mazzeo
\cite{Mazzeo-thesis,Mazzeo-Hodge}.  (See \cite{Andersson}, where this is
carried out in detail.)  Our main interest in this monograph, however, is in
the covariant Laplacian acting on symmetric tensors.  In this case, we
consider the bundle $\Sigma^rM$ of symmetric $r$-tensors as a subbundle
of the bundle $\Lambda^1 T^{r-1}M$ of $(r-1)$-tensor-valued $1$-forms.

The following lemma gives a 
Weitzenb\"ock formula relating the covariant Laplacian on such
tensors to $\Delta$.
For this purpose, we define 
zero-order operators 
$\tw{Rc}(u)$ and $\tw{Rm}(u)$ 
acting on $r$-tensors by
\begin{align*}
\tw{Rc}(u)_{i_1\dots i_{r-1}j} &= R_j{}^k u_{i_1\dots i_{r-1}k};\\
\tw{Rm}(u)_{i_1\dots i_{r-1}j} &= \sum_{p=1}^{r-1} \tensor
R{\down{i_p}\up {l} \down{j} \up {k}}
u_{i_1\dots l\dots i_{r-1}k}
.
\end{align*}
Note that on symmetric $2$-tensors,
$\tw{Rm}$ agrees with the operator
$\mathring{Rm}$ defined by 
\eqref{eq:def-Rc-Rm}, but $\tw{Rc}$ is not the
same as $\mathring{Rc}$ in general.

\begin{lemma}\label{weitzenbock-lemma}
For a
section $u$ of $\Lambda^1 T^{r-1}M$,
\begin{equation}\label{weitzenbock}
\Delta u = \Del^*\Del u + \tw{Rc}(u) - \tw{Rm}(u).
\end{equation}
\end{lemma}

\begin{proof}
This is easiest to see in components, noting that the {\it last} index
of $u$ is considered to be the $1$-form index.
\begin{align*}
(D^*Du)_{i_1\dots i_{r-1}j}
&= - u_{i_1\dots i_{r-1}j;k}{}^k +
u_{i_1\dots i_{r-1}k;j}{}^k;\\
(DD^*u)_{i_1\dots i_{r-1}j}
&= - u_{i_1\dots i_{r-1}k;}{}^{k}{}_{j}.
\end{align*}
Applying the Ricci identity to the commutator
\begin{displaymath}
u_{i_1\dots i_{r-1}k;j}{}^k- u_{i_1\dots i_{r-1}k;}{}^{k}{}_{j}
\end{displaymath}
yields the result.
\end{proof}

Specializing Lemma \ref{asymptotic-estimate-lemma} to $\Sigma^r_0M$, we
obtain the following sharp asymptotic estimate.  

\begin{lemma}\label{symmetric-tensor-estimate}
The
following asymptotic estimate holds for any
smooth, compactly supported,
trace-free symmetric $r$-tensor $u$:
\begin{displaymath}
(u,\Del^*\Del u)\gtrsim
\left(\frac{n^2}{4}+r\right)\|u\|^2.
\end{displaymath}
\end{lemma}

\begin{proof}
Using \eqref{eq:Rijkl}, we compute that the curvature operators
$\tw{Rc}$ and $\tw{Rm}$ have the following asymptotic
behavior on trace-free 
symmetric $r$-tensors near $\del M$:
\begin{align*}
\tw{Rc}(u) &= -nu + O(\rho|u|_g);\\
\tw{Rm}(u) &= (r-1)u + O(\rho|u|_g).
\end{align*}
To prove the lemma, just use
Lemma \ref{asymptotic-estimate-lemma} with $q=1$ together with
Lemma \ref{weitzenbock-lemma} to obtain
\begin{align*}
(u,\Del^*\Del u)
&= (u,\Delta u) - (u,\tw{Rc}(u)) + (u,\tw{Rm}(u))\\
&\gtrsim \frac{(n-2)^2}{4} \|u\|^2 + n\|u\|^2 + (r-1)\|u\|^2\\
&=\left( \frac{n^2}{4}+r\right) \|u\|^2.
\end{align*}
\end{proof}

\begin{remark}
The same method
yields similar estimates  for 
 Laplace operators acting on
any other tensor bundle, but in cases other than fully symmetric or
fully antisymmetric tensors, the estimates so obtained appear not to be
sharp.  It would be interesting to know if a modified version of
\eqref{integral-formula} could be used to sharpen the estimate in those
cases.
\end{remark}

\begin{lemma}\label{lemma:lichnerowicz-estimate}
The
following asymptotic estimate holds for any
smooth, compactly supported,
trace-free symmetric $2$-tensor $u$:
\begin{displaymath}
(u,\Delta_{\text{L}} u)\gtrsim
\left(\frac{n^2}{4}-2n\right)\|u\|^2.
\end{displaymath}
\end{lemma}

\begin{proof}
This follows immediately from the preceding lemma
together with the 
asymptotic formulas \eqref{eq:est-for-Rm-Rc}
for $\mathring{Rm}$ and $\mathring{Rc}$.
\end{proof}

\begin{proof}[Proof of Propositions \ref{prop:lichnerowicz},
\ref{prop:cov-Laplacian}, 
\ref{prop:hodge-laplacian}, and
\ref{prop:vector-laplacian}.]
The operators $\Delta_{\text{L}}+c$ and $\Del^*\Del+c$
are Fredholm on $L^2$ 
for the claimed values of $c$ by virtue of 
Lemmas \ref{lemma:lichnerowicz-estimate} and
\ref{symmetric-tensor-estimate}.
On the other hand, the indicial radius
computations at the beginning of this chapter show 
that these are precisely the values of $c$ for
which these operators have positive indicial radius,
so these are the only values of $c$ for which the
operators are Fredholm.
The result for the Hodge Laplacian follows similarly
from Lemma
\ref{asymptotic-estimate-lemma} in the case of scalar forms.
The claims about the essential spectrum follow immediately
from the Fredholm results: Since each of these
operators is self-adjoint on $L^2$, 
its spectrum is contained in $\R$, 
and a real number $\lambda$ is in the essential spectrum
if and only if $P-\lambda\colon L^2(M;E)\to L^2(M;E)$
is not Fredholm.  

For the vector Laplacian $L^*L$, 
L. Andersson proved in \cite[Lemma 3.15]{Andersson}
that the following asymptotic estimate holds:
\begin{displaymath}
(V,L^*LV) \gtrsim \frac{n^2}{8}\|V\|^2.
\end{displaymath}
It follows from Proposition \ref{prop:L2-Fredholm}
that $L^*L$ is Fredholm on $L^2$, and then the
rest of the claims follow from Theorem \ref{thm:main-fredholm}.
\end{proof}

We conclude this chapter by observing that the 
sharp a priori $L^2$ estimates developed here for
Laplace operators actually lead directly to sharp
$L^2$ Fredholm theorems for such operators, without any
need for the parametrix construction of the preceding
chapter.  This follows from the 
simple device of replacing $u$ by 
$\rho^{-\delta}u$ to convert
unweighted $L^2$ estimates 
to weighted ones.  (Cf.~ also \cite[Lemma 3.8]{Andersson}.)

\begin{lemma}\label{weighted-from-unweighted}
Let $E$ be any geometric
tensor bundle over $M$, and let $P=\Del^*\Del+\scr K$ be a 
Laplace operator acting on sections of $E$.  Suppose that $P$
satisfies the asymptotic estimate 
\begin{equation}\label{asymptotic-estimate}
(u,Pu)\gtrsim \lambda\|u\|^2
\end{equation}
for some  $\lambda>0$.  If $|\delta|^2<\lambda$, 
there is a compact set $K\subset M$
such that the following weighted
estimate holds for all 
$u\in C^\infty_c(M\setminus K;E)$:
\begin{equation}\label{weighted-asymptotic-estimate}
\|u\|_{0,2,\delta}
\le C\|Pu\|_{0,2,\delta}.
\end{equation}
\end{lemma}

\begin{proof}
If $u\in C^\infty_c(M;E)$, then
\begin{equation}\label{del*del-phiu}
\begin{aligned}
\big(\rho^{-\delta} u,\Del^*\Del(\rho^{-\delta} u)\big)
&=\|\Del (\rho^{-\delta} u)\|_g^2\\
&=\int_M \left(\rho^{-2\delta} |\Del u|_g^2
-2\delta \rho^{-2\delta}
\left< \Del u, u\tprod \frac{d\rho}{\rho}\right>_{\kern -4pt g}\right.\\
&\qquad\qquad 
+\left.\delta^2\rho^{-2\delta}|u|_g^2\left|\frac{d\rho}{\rho}\right|_g^2\right)\,dV_g.
\end{aligned}
\end{equation}
By the divergence theorem and the fact that
$\Delta|u|_g^2=2\left<u,\Del^*\Del u\right>_g -2|\Del u|_g^2$,
\begin{displaymath}
\begin{aligned}
0 &= \frac 1 2 \int_M d^*\left(  \rho^{-2\delta} d|u|_g^2\right)\,dV_g\\
&= \int_M\rho^{-2\delta} \left<u,\Del^*\Del u\right>_g 
-\rho^{-2\delta} |\Del u|_g^2 
+ 2\delta \rho^{-2\delta} \left< \Del u, u\tprod \frac{d\rho}{\rho}\right>_{\kern -4pt g}
\, dV_g.
\end{aligned}
\end{displaymath}
Substituting this into \eqref{del*del-phiu}, we obtain
\begin{equation}\label{phiu}
\big(\rho^{-\delta} u,\Del^*\Del(\rho^{-\delta} u)\big)
= (\rho^{-\delta} u, \rho^{-\delta} \Del^*\Del u)
+ \delta^2\int_M\rho^{-2\delta}|u|_g^2\left|\frac{d\rho}{\rho}\right|_g^2\,dV_g.
\end{equation}

Given $\epsilon>0$, choose $r>0$ small enough  that
\begin{equation}\label{unweighted-L2-estimate}
(u,Pu)\ge (\lambda-\epsilon)\|u\|^2
\end{equation}
whenever $u$ is smooth and compactly supported in $A_r$, 
and such that  $|d\rho/\rho|_g^2\le 1+\epsilon/(\delta^2)$
on $A_r$.
Applying \eqref{unweighted-L2-estimate} 
to $\rho^{-\delta} u$ and using \eqref{phiu}, we obtain
\begin{displaymath}
\begin{aligned}
(\lambda-\epsilon)\|\rho^{-\delta}u\|^2
&\le \big(\rho^{-\delta}u,(\Del^*\Del+\scr K)(\rho^{-\delta} u)\big)\\
&\le \big(\rho^{-\delta}u,\rho^{-\delta} (\Del^*\Del+\scr K)u\big)
+(\delta^2+\epsilon)\|u\|_{0,2,\delta}^2\\
&\le \|u\|_{0,2,\delta}\|Pu\|_{0,2,\delta}
+(\delta^2+\epsilon)\|u\|_{0,2,\delta}^2,
\end{aligned}
\end{displaymath}
which implies \eqref{weighted-asymptotic-estimate}
as long as $2\epsilon < \lambda - \delta^2$.
\end{proof}

With this estimate in hand, it follows
immediately from Lemma 
\ref{lemma:fred} that $P$ is 
Fredholm on $H^{k,2}_{\delta}(M;E)$
provided that $|\delta|^2<\lambda$.
If only $L^2$ results are
needed, this provides an exceedingly
elementary approach that can be used in many
cases.

%%
%%
%% Fredholm operators and Einstein metrics
%% on conformally compact manifolds
%%
%% by John M. Lee
%% 
%% Chapter 8
%%

\chapter{Einstein Metrics}\label{section:einstein}

In this chapter, we will apply the linear theory we have developed
so far to prove the existence of Einstein metrics.  To do so, we first need
to describe a systematic construction of asymptotic solutions to the
Einstein equation, which will
then be corrected by appealing to 
an appropriate linear isomorphism theorem 
and the inverse function theorem.
The construction of asymptotic solutions is 
very similar to 
the one we described on hyperbolic space
in \cite{Graham-Lee}, but more delicate because
we start with less regularity at the boundary.

Suppose $(M,h)$ is an asymptotically hyperbolic Einstein
$(n+1)$-manifold of class $C^{l,\beta}$, with $l\ge 2$ and
$0< \beta<1$.
Let $\rho$ be a smooth defining function for $h$, and put
$\overline h = \rho^2 h$, $\hat h = \overline h|_{\del M}$.

Recall the spaces $C^{k,\alpha}_{(s)}(\overline M;E)$
defined in Chapter \ref{spaces-section}.
By definition of the indicial map, if $P$ is
a uniformly degenerate operator of order $m$ and 
$\overline u\in C^{m,\alpha}_{(0)}(\overline M;E)$, then
the behavior of 
$P(\rho^s\overline u)$ at the boundary is dominated by
$\rho^sI_s(P)\overline u$.
This in turn follows from the fact that
$\rho^{-s}\del_{i_1}\dotsm \del_{i_j}(\rho^s \overline u)
=
\rho^{-s}\del_{i_1}\dotsm \del_{i_j}(\rho^s) \overline u + o(\rho^{-j})$
in background coordinates whenever $\overline u\in C^{j,\alpha}_{(0)}(\overline M)$.
In the construction of our asymptotic solutions, 
we will need to make similar estimates when $\overline u$ has
somewhat less regularity.  The key is the following lemma.

\begin{proposition}\label{prop:rho-derivs}
Suppose $0<\alpha<1$, $0<k+\alpha\le l+\beta$, $\delta\in\R$,
and $s,j$ are integers satisfying $1\le j\le s
< k+\alpha\le l+\beta$.  
If $u\in C^{k,\alpha}_{(s)}(\overline M)$,
then for any indices $1\le i_1,\dots,i_j\le n+1$,
\begin{displaymath}
\rho^{-\delta}\del_{i_1}\dotsm \del_{i_j}(\rho^\delta u)
- \rho^{-\delta-s} \del_{i_1}\dotsm \del_{i_j}(\rho^{\delta+s})u 
\in C^{k-j,\alpha}_{(s-j+\alpha)}(\overline M).
\end{displaymath}
\end{proposition}

\begin{proof} 
The proof is by induction on $j$.
For $j=1$, 
consider first the case $\delta=0$, in which case
the claim is 
\begin{equation}\label{eq:j=0-case}
\del_i
u = s\rho^{-1}(\del_i\rho)u
\mod C^{k-1,\alpha}_{(s-1+\alpha)}(\overline M).
\end{equation}
We will prove this claim by induction on $k$.  
Observe that 
$\del_i u$ and $s\rho^{-1} (\del_i \rho)u$ are in 
$C^{k-1,\alpha}_{(s-1)}(\overline M)$ by 
parts \eqref{part:derivs-map-vanishing}
and \eqref{part:vanishing-mult-neg}
of Lemma 
\ref{lemma:properties-of-vanishing-spaces}.  Thus to
prove \eqref{eq:j=0-case}, it suffices to show that the
difference between the two terms is $O(\rho^{s-1+\alpha})$.

For $k=1$, the only value of $s$ that satisfies the
hypotheses is $s=1$.
Since $\rho=\theta^{n+1}$ is a coordinate function,
we consider separately the cases 
$i<n+1$ and 
$i=n+1$.
When $i< n+1$ the right-hand side of \eqref{eq:j=0-case} is zero.
On the other hand, the left-hand side is in $C^{0,\alpha}_{(0)}(\overline M)$
and vanishes on $\del M$ (since $u$ vanishes on $\del M$ and $\del_i$ is
tangent to $\del M$), so it is $O(\rho^\alpha)$.
When $i=n+1$, $\rho^{-1}u\in C^{0,\alpha}_{(0)}(\overline M)$
by Lemma 
\ref{lemma:properties-of-vanishing-spaces}\eqref{part:vanishing-mult-neg},
and by definition of the derivative,
\begin{align*}
\del_{n+1}u(\theta,0) 
&= \lim _{\rho\to 0}\frac{u(\theta,\rho) - u(\theta,0)}{\rho}\\
&= \lim _{\rho\to 0}\rho^{-1} u(\theta,\rho).
\end{align*}
It follows that $\del_{n+1}
u - \rho^{-1}(\del_{n+1}\rho)u =
\del_{n+1}
u - \rho^{-1}u \in C^{0,\alpha}_{(0)}(\overline M)$
and vanishes on $\del M$, so it too is 
$O(\rho^\alpha)$.
This proves \eqref{eq:j=0-case} in the $k=s=1$ case.

Now let $k>1$, and suppose
\eqref{eq:j=0-case} holds for all smaller values of $k$.
Observe that
$v=\rho^{1-s}u\in C^{k-s+1,\alpha}_{(1)}(\overline M)$
by Lemma 
\ref{lemma:properties-of-vanishing-spaces}\eqref{part:vanishing-mult-neg}.
Therefore, by the chain rule and the inductive hypothesis,
\begin{align*}
\del_i u
&= \del_i (\rho^{s-1}v)\\
&= (s-1) \rho^{s-2}(\del_i \rho)v
+ \rho^{s-1} \del_i v\\
&= (s-1) \rho^{-1} (\del_i \rho)u 
+ \rho^{s-1}(\rho^{-1}(\del_i\rho)v+ O(\rho^\alpha))\\
&= s\rho^{-1}(\del_i \rho) u  + O(\rho^{s-1+\alpha}).
\end{align*}
Finally, for 
$\delta\ne 0$, the product rule gives 
\begin{align*}
\rho^{-\delta}\del_i(\rho^\delta u) 
&= \rho^{-\delta}(\delta\rho^{\delta-1} (\del_i\rho) u + \rho^\delta \del_i u)\\
&= \delta\rho^{-1} (\del_i\rho) u + s\rho^{-1}(\del_i\rho)u
+ O(\rho^{s-1+\alpha})\\
&= \rho^{-\delta-s}\del_i(\rho^{\delta+s})u + O(\rho^{s-1+\alpha}).
\end{align*}
This completes the proof of the $j=1$ step.

Now suppose the proposition is true for some $j\ge 1$.  For
any $(j+1)$-tuple $(i_1,\dots,i_{j+1})$, beginning with the 
induction hypothesis in the form
\begin{align*}
&\del_{i_2}\dotsm\del_{i_{j+1}}(\rho^{\delta}u)
- \rho^{-s}\del_{i_2}\dotsm\del_{i_{j+1}}(\rho^{\delta+s})u\\
&\qquad\in \rho^\delta C^{k-j,\alpha}_{(s-j+\alpha)}(\overline M)
\subset C^{k-j,\alpha}_{(\delta+s-j+\alpha)}(\overline M),
\end{align*}
we apply the chain rule to obtain 
\begin{align*}
\rho^{-\delta}\del_{i_1}\dotsm\del_{i_{l+1}}(\rho^\delta u) 
&=
\rho^{-\delta}\big(-s\rho^{-s-1}(\del_{i_1}\rho)\del_{i_2}\dotsm
\del_{i_{l+1}}(\rho^{\delta+s})u\\ 
&\qquad +\rho^{-s} \del_{i_1}\dotsm\del_{i_{l+1}}(\rho^{\delta+s})u\\
&\qquad +\rho^{-s} \del_{i_2}\dotsm\del_{i_{l+1}}(\rho^{\delta+s}) \del_{i_1}u\big)
+O(\rho^{s-j-1+\alpha})\\
&= \rho^{-\delta-s} \del_{i_1}\dotsm\del_{i_{l+1}}(\rho^{\delta+s})u
+O(\rho^{s-j-1+\alpha}),
\end{align*}
where in the last line we have used the induction hypothesis again to 
evaluate $\del_{i_1}u$.
\end{proof}

\begin{corollary}\label{cor:indicial-op-irregular}
Let $P\colon C^\infty(M;E)\to C^\infty(M;E)$ be a self-adjoint,
elliptic, geometric partial differential operator of order $m\le l$.
Suppose $0<\alpha<1$, $0<k+\alpha\le l+\beta$, $\delta\in\R$,
and $u\in C^{k,\alpha}_{(s)}(\overline M;E)$, where 
$s$ is an integer satisfying $1\le s
< k+\alpha\le l+\beta$.  
Then 
\begin{equation}
\rho^{-\delta-s}P(\rho^\delta u)|_{\del M}
 = I_{\delta+s}(P)\hat u,
\end{equation}
where $\hat u = \rho^{-s}u|_{\del M}$.
\end{corollary}

\begin{proof}
Write $Pu$ in background coordinates as
\begin{displaymath}
Pu =  \sum_{0\le j\le m}\sum_{i_1,\dots,i_j}
\rho^j a^{i_1,\dots,i_j} \del_{i_1}\dotsm\del_{i_j} u,
\end{displaymath}
where $u$ is vector-valued and the coefficient functions
$a^{i_1,\dots,i_j}$ are matrix-valued.
If $\overline u\in C^{m,\alpha}_{(0)}(M;E)$, then an easy computation
shows that 
\begin{displaymath}
\rho^{-\delta-s}P(\rho^{\delta+s}\overline u) 
= \rho^{-\delta-s} \sum_{0\le j\le m}\sum_{i_1,\dots,i_j}
\rho^j a^{i_1,\dots,i_j}
\del_{i_1}\dotsm\del_{i_j}(\rho^{\delta+s})\overline  u
+ o(1),
\end{displaymath}
which implies that 
\begin{displaymath}
I_{\delta+s}(P)\hat u
 = \sum_{0\le j\le m}\sum_{i_1,\dots,i_j}(\rho^{-\delta-s} 
\rho^j a^{i_1,\dots,i_j} \del_{i_1}\dotsm\del_{i_j}(\rho^{\delta+s}))|_{\del M}\hat u.
\end{displaymath}
On the other hand, Proposition  
\ref{prop:rho-derivs} shows that 
\begin{displaymath}
\rho^{-\delta-s}P(\rho^{s}u)
= 
\sum_{0\le j\le m}\sum_{i_1,\dots,i_j}
\rho^{-\delta-s}
\rho^j a^{i_1,\dots,i_j}
\del_{i_1}\dotsm\del_{i_j}(\rho^{\delta+s}) (\rho^{-s}
u)+O(\rho^\alpha),
\end{displaymath}
which proves the result.
\end{proof}

We will need an extension operator from tensor fields
on $\del M$ to $C^{l,\beta}$
tensor fields on $\overline M$.  
In \cite{Graham-Lee}, we did this by 
parallel translating along $\overline h$-geodesics
normal to $\del M$.  However, in our present 
circumstances, because $\overline h$ may not be smooth 
up to the boundary, this parallel translation 
would lose regularity.  Instead, we define
our extension operator as follows.

Let $\pi\colon T\overline M|_{\del M}\to T\del M$
denote the $\overline h$-orthogonal projection, which is 
clearly a $C^{l,\beta}$ bundle map over $\del M$.
Given any $C^{l,\beta}$
$2$-tensor field $v$ on $\del M$,
lifting by $\pi$ yields a 
$C^{l,\beta}$ section $\pi^*v$
of $\Sigma^2\overline M|_{\del M}$:
\begin{displaymath}
\pi^*v(X,Y) = v(\pi X, \pi Y).
\end{displaymath}
Define a linear
map $E\colon C^{l,\beta}(\del M;\Sigma^2\del M)
\to C^{l,\beta}_{(0)}(\overline M;\Sigma^2 \overline M)$ by
letting $E(\hat u)=\phi \Pi (\pi^*u)$, where $\phi$ is a 
fixed cutoff function that is equal to $1$ along
$\del M$ and is supported in a small collar neighborhood
of the boundary, and $\Pi$ denotes parallel translation
along normal geodesics with respect to some smooth
background metric.

For any $C^{l,\beta}$ Riemannian metric $\hat g$ 
on $\del M$, we define a metric $T(\hat g)$ on $M$
by
\begin{equation}\label{eq:def-T}
T(\hat g) = h + \rho^{-2}E(\hat g - \hat h).
\end{equation}
Then $T\colon C^{l,\beta}(\del M;\Sigma^2\del M)\to 
C^{l,\beta}(M;\Sigma^2 M)$ defines a smooth (in fact affine) map
of Banach spaces.
It is easy to 
see that if 
$\hat g$ is sufficiently close to $\hat h$ in the
$C^{l,\beta}$ norm,
then $T(\hat g)$ will be an asymptotically
hyperbolic 
metric on $M$ of class $C^{l,\beta}$,
whose conformal infinity
is $[\hat g]$.
Moreover, since $E(0)=0$,
our construction guarantees that
$T(\hat h) = h$.

For higher-order asymptotics, we will
use an extension lemma
due to L. Andersson and P. Chru{\'s}ciel
\cite[Lemma 3.3.1]{AC}.  Translated into our notation, 
it reads as follows.

\begin{lemma}\label{lemma:extension}
Suppose $0\le\alpha<1$ and $k\ge 0$.  Given an integer
$s$ such that $0\le s \le k$ and any $\psi\in C^{k-s,\alpha}(\del M)$,
there exists a function $u\in C^{k,\alpha}_{(s)}(\overline M)$
such that  $(\rho^{-s} u)|_{\del M} = \psi$.
\end{lemma}

The formula \cite[(3.3.6)]{AC} for the derivatives
of $u$ makes it clear that the mapping $\psi\mapsto u$ is
 a bounded linear map from $C^{k-s,\alpha}(\del M)$
to $C^{k,\alpha}_{(s)}(\overline M)$.

Recall the operator $Q$ defined by \eqref{regularized-einstein}.
The next lemma gives a construction of asymptotic
solutions to $Q(g,g_0)=0$  analogous to 
Theorem 2.11 of \cite{Graham-Lee}.

\begin{lemma}\label{lemma:asymp-expansion}
Suppose $0<\beta<1$, 
$2\le l\le n-1$, and $h$ is an asymptotically
hyperbolic metric on $M$ of class $C^{l,\beta}$.
Let $\hat g$ be any metric on $\del M$ of class $C^{l,\beta}$,
and set $g_0 = T(\hat g)$.
There exists an asymptotically hyperbolic 
metric $g$ of class $C^{l,\beta}$
on $M$ such that $\rho^2 g|_{\del M} = \hat g$ and
\begin{equation}\label{eq:asymptotic-Q}
Q(g,g_0) \in C^{l-2,\beta}_{(l-2+\beta)}(\overline
M;\Sigma^2\overline M)\subset
C^{l-2,\beta}_{l+\beta}(M;\Sigma^2 M).
\end{equation}
The mapping $S\colon C^{l,\beta}(\del M;\Sigma^2\del M) \to
\rho^{-2}C^{l,\beta}_{(0)}(\overline M;\Sigma^2\overline M)$
given by $S(\hat g) =  g$ is a smooth map of Banach spaces.
\end{lemma}

\begin{proof}
Define a nonlinear
operator 
\begin{displaymath}
\overline Q\colon C^{l,\beta}_{(0)}(\overline
M;\Sigma^2\overline M) \cross C^{l,\beta}_{(0)}(\overline
M;\Sigma^2\overline M) \to 
C^{l-2,\beta}_{(0)}(\overline
M;\Sigma^2\overline M)
\end{displaymath}
 by
\begin{displaymath}
\overline Q(\overline g,\overline g_0) = \rho^2
Q(\rho^{-2}g,\rho^{-2}g_0).
\end{displaymath}
It is clear that \eqref{eq:asymptotic-Q} is equivalent
to $\overline Q(\rho^2 g, \rho^2 g_0) = O(\rho^{l-2+\beta})$.
It follows from formula (2.19) of 
\cite{Graham-Lee} that 
$\overline Q$ has  the following expression in background coordinates:
\begin{displaymath}
\overline Q(\overline g,\overline g_0) 
= \scr E^0(\overline g,\overline g_0)
+ \rho \scr E^1(\overline g, \overline g_0)
+ \rho^2 \scr E^2(\overline g, \overline g_0),
\end{displaymath}
where $\scr E^j$ is a universal polynomial in 
the components of $\overline g$,
$\overline g_0$, their inverses,
and their coordinate derivatives of order
less than or equal to $j$.
It follows that $\overline Q$ takes its values in
\begin{equation}\label{eq:Q-target}
 C^{l,\beta}_{(0)}(\overline M;\Sigma^2\overline M)
+ \rho C^{l-1,\beta}_{(0)}(\overline M;\Sigma^2\overline M)
+ \rho^2 C^{l-2,\beta}_{(0)}(\overline M;\Sigma^2\overline M),
\end{equation}
and is a smooth map of Banach spaces.

We will recursively construct a sequence of  metrics
$\overline g_0,\dots,\overline g_l\in C^{l,\beta}_{(0)}
(\overline M;\Sigma^2\overline M)$ satisfying
$\overline Q(\overline g_k,\overline g_0)= o(\rho^{k})$.
In fact, we will prove a bit more, namely that
\begin{equation}\label{eq:Qbar-good}
\overline Q(\overline g_k,\overline g_0)\in 
\begin{cases}
\rho C^{l-1,\beta}_{(0)}(\overline M;
\Sigma^2\overline M) +
\rho^2 C^{l-2,\beta}_{(0)}(\overline M;
\Sigma^2\overline M),& k=0,\\
\rho^2 C^{l-2,\beta}_{(k-1)}(\overline M;
\Sigma^2\overline M), & 1\le k\le l-1,\\
\rho^2 C^{l-2,\beta}_{(l-2+\beta)}(\overline M;
\Sigma^2\overline M), & k=l.
\end{cases}
\end{equation}

Begin with $\overline g_0 =\rho^2 T(\hat g)$.
It follows from Corollary
2.6 of \cite{Graham-Lee} that 
\begin{displaymath}
\overline Q(\overline g_0,\overline g_0) = O(\rho).
\end{displaymath}
Since the intersection of \eqref{eq:Q-target}
with $O(\rho)$ is 
\begin{align*}
&C^{l,\beta}_{(1)}(\overline M;
\Sigma^2\overline M) +
\rho C^{l-1,\beta}_{(0)}(\overline M;
\Sigma^2\overline M) +
\rho^2 C^{l-2,\beta}_{(0)}(\overline M;
\Sigma^2\overline M)\\
&\qquad \subset 
\rho C^{l-1,\beta}_{(0)}(\overline M;
\Sigma^2\overline M) +
\rho^2 C^{l-2,\beta}_{(0)}(\overline M;
\Sigma^2\overline M),
\end{align*}
 we have 
\eqref{eq:Qbar-good} in the $k=0$ case.

Assume by induction that for some $k$, $1\le k \le l$,
 we have constructed 
$\overline g_{k-1}\in C^{l,\beta}_{(0)}(\overline M;
\Sigma^2\overline M)$ 
satisfying the analogue of \eqref{eq:Qbar-good}.
Then
letting $\hat v = \rho^{-k}\overline Q(\overline g_{k-1},\overline g_0)|_{\del
M}$, we see from  Lemma 
\ref{lemma:properties-of-vanishing-spaces}\eqref{part:vanishing-mult-neg} that
$\hat v$ is a $C^{l-k,\beta}$ section of $T\overline M|_{\del M}$.

If $\overline r\in C^{l,\beta}_{(k)}(\overline
M;\Sigma^2\overline M)$, the
same argument as in the proof of \cite[Theorem 2.11]{Graham-Lee}
shows that 
\begin{equation}\label{eq:asymptotic-linearization}
\overline Q(\overline g_{k-1}+\overline r,\overline g_0)
= \overline Q(\overline g_{k-1},\overline g_0)
+ \rho^2(\Delta_{\text{L}} +2n) (\rho^{-2}\overline r)
+ o(\rho^{k}).
\end{equation}
By Corollary 
\ref{cor:indicial-op-irregular}, we have
\begin{displaymath}
\rho^{2-k}(\Delta_{\text{L}} +2n) (\rho^{-2}\overline r)|_{\del M}
 = I_{k-2}(\Delta_{\text{L}} +2n)(\rho^{-k} \overline r)|_{\del M}
+ o(1).
\end{displaymath}
Because the indicial radius of $\Delta_{\text{L}}+2n$ is
$R=n/2$, 
$I_{s}(\Delta_{\text{L}}+2n)$ is invertible 
provided $-2<s<n-2$.  Thus as long as $1\le k \le l \le n-1$,
there is a unique $C^{l-k,\beta}$ tensor field 
$\psi$ along $\del M$ such that 
$I_{k-2}(\Delta_{\text{L}} +2n)\psi = -\hat v$.
By Lemma \ref{lemma:extension}, there is a 
tensor field $\overline r\in C^{l,\beta}_{(k)}
(\overline M;\Sigma^2\overline M)$
such that $(\rho^{-k}\overline r)|_{\del M} = \psi$.
Inserting this back into \eqref{eq:asymptotic-linearization},
we conclude that 
$\overline Q(\overline g_{k-1} + \overline r,\overline g_0)
= o(\rho^k)$,
and thus actually satisfies
\eqref{eq:Qbar-good}.
Setting $\overline g_k = \overline g_{k-1} + \overline r$ completes the
inductive step.

After the $k=l$ step, we set $g=g_l$, and 
\eqref{eq:asymptotic-Q} is satisfied.
The smoothness of the operator $S:\hat g\mapsto g$
is immediate from the remark following Lemma
\ref{lemma:extension} above.
\end{proof}

\begin{proof}[Proof of Theorem \ref{thm:einstein}]
The proof follows closely that of Theorem 4.1 in 
\cite{Graham-Lee}.  We define an open subset
$\scr B\subset C^{l,\beta}(\del M;\Sigma^2\del M)
\cross C^{l,\beta}_{l+\beta}(M;\Sigma^2M)$ by
\begin{displaymath}
\scr B = \{(\hat g, r): \text{$\hat g$, $T(\hat g)$, and $S(\hat g) + r$ are
all positive definite}\}.
\end{displaymath}
Define a map $\scr Q\colon \scr B \to 
C^{l,\beta}(\del M;\Sigma^2\del M)
\cross C^{l-2,\beta}_{l+\beta}(M;\Sigma^2M)$
by
\begin{displaymath}
\scr Q(\hat g,r) = (\hat g, Q(S(\hat g) + r, T(\hat g))),
\end{displaymath}
where $T$ is defined in \eqref{eq:def-T}
and $S$ in Lemma \ref{lemma:asymp-expansion}.
It follows just as in 
\cite{Graham-Lee} that $\scr Q$ is a smooth map
of Banach spaces.

Since $T(\hat h)=S(\hat h)=h$ and $h$ is an Einstein metric, 
$\scr Q(\hat h,0) = (\hat h,0)$.
As is shown in \cite{Graham-Lee}, 
the linearization of $\scr Q$ about $(\hat h,0)$ is
the linear map $D\scr Q_{(\hat h,0)}$
from $C^{l,\beta}(\del M;\Sigma^2\del M)\cross
C^{l,\beta}_{l+\beta}(M;\Sigma^2M)$ to $C^{l,\beta}(\del M;\Sigma^2\del M)\cross
C^{l-2,\beta}_{l+\beta}(M;\Sigma^2M)$ given by
\begin{align*}
D\scr Q_{(\hat h,0)} (\hat q,r)
&= (\hat q,
D_1 Q_{(h,h)} (DS_{\hat h}\hat q + r) + D_2Q_{(h,h)}(DT_{\hat h}\hat q))\\
&= (\hat q, (\DeltaL + 2n)r + K\hat q),
\end{align*}
where 
\begin{displaymath}
K\hat q = D_1Q_{(h,h)}(DS_{\hat h}\hat q) 
+ D_2Q_{(h,h)}(DT_{\hat h}\hat q).
\end{displaymath}

Because $\DeltaL+2n$ preserves the 
splitting $\Sigma^2M = \Sigma^2_0M \oplus \R g$,
to show that it satisfies the hypotheses of Theorem
\ref{thm:main-fredholm}, 
we need to show that it has trivial $L^2$
kernel on each
of these bundles separately.  
Since $\mathring{Rm}(g)=
\mathring{Rc}(g)$, the action of $\DeltaL$
on sections of $\R g$ is just given by the scalar 
Laplacian:
\begin{displaymath}
\DeltaL(ug) = (\Del^*\Del u)g.
\end{displaymath}
Since $2n>0$, it follows easily from integration
by parts that $\Del^*\Del u+2n$ acting
on scalar functions has trivial $L^2$ kernel.  
Thus the assumption that 
$\DeltaL+2n$
has trivial $L^2$ kernel on $\Sigma^2_0 M$
is sufficient to allow us to apply the results of
Theorem \ref{thm:main-fredholm}.
In particular, 
$\DeltaL+2n\colon C^{l,\beta}_\delta(M;\Sigma^2M)
\to C^{l-2,\beta}_\delta(M;\Sigma^2M)$ is an isomorphism
for $|\delta-n/2|<n/2$, which is to say for $0<\delta<n$.
This is true for $\delta=l+\beta$, so
the linearization of $\scr Q$ has a bounded
inverse given by
\begin{displaymath}
(D\scr Q_{(\hat h,0)})^{-1} (\hat w,v) = (\hat w, (\DeltaL + 2n)^{-1}(v-K\hat w)).
\end{displaymath}
Therefore, by the Banach space 
inverse function theorem, there is a neighborhood 
of $(\hat h,0)$ on which $\scr Q$ has a smooth inverse.
In particular, for $\hat g$ sufficiently $C^{l,\beta}$
close to 
$\hat h$, there is a solution $r\in C^{l,\beta}_{l+\beta}(M;\Sigma^2M)$ to 
$\scr Q(\hat g,r) = (\hat g,0)$.

Putting $g=S(\hat g)+r$ and $g_0=T(\hat g)$, we have $Q(g,g_0)=0$.
Moreover, if the $C^{l,\beta}$ neighborhood of $\hat h$
is sufficiently small, then $g$ will be uniformly $C^2$ close
to $h$, and therefore $g$ will have strictly negative
Ricci curvature.  By Lemma 2.2 of \cite{Graham-Lee}, this implies that
$g$ is Einstein.
Note that 
$\rho^2 r \in C^{l,\beta}_{l+2+\beta}(M;\Sigma^2 M)$,
which is contained in $C^{l,\beta}_{(0)}(\overline M;
\Sigma^2\overline M)$ by Lemma 
\ref{lemma:spaces-m-and-mbar}.
Since $\rho^2 S(\hat g)\in C^{l,\beta}_{(0)}(\overline M;
\Sigma^2\overline M)$ by construction, it follows that
$g$ is asymptotically hyperbolic of class $C^{l,\beta}$
as claimed.

It remains only to prove that $\DeltaL+2n$ has
trivial $L^2$  kernel on
$\Sigma^2_0M$ under the assumptions stated
in Theorem \ref{thm:einstein}.
First suppose that $h$ has nonpositive sectional
curvature.  A simple algebraic argument (see
\cite[Lemma 12.71]{Besse}) shows that 
if $h$ is an Einstein metric on an $(n+1)$-manifold
with scalar curvature $-n(n+1)$,
its Riemann curvature operator
$\mathring{Rm}$ acting on trace-free symmetric $2$-tensors
satisfies the following 
estimate at each point $p\in M$:
\begin{equation}\label{eq:Koiso-est}
\< \mathring{Rm}(u),u\>_h \le \big( n + (n-1)K_{\max}(p)\big)|u|_h^2,
\end{equation}
where 
$K_{\max}(p)$ is the maximum of the sectional curvatures
of $h$ at $p$.
(In \cite{Besse}, this is
attributed to a hard-to-find
1979 paper of T. Fujitani; however,
the argument was already given in 1978 by N. Koiso 
\cite[Prop.~ 3.4]{Koiso-nondeform}.)
Since the Einstein assumption implies that
$\mathring{Rc}(u) = - n u$,
we have
\begin{displaymath}
\DeltaL + 2n = \Del^*\Del - 2 \mathring{Rm} = DD^* + D^*D +n 
-\mathring{Rm}
\end{displaymath}
(see
\eqref{weitzenbock}).
Therefore, if $h$ has sectional curvatures everywhere
bounded above by $-\kappa\le 0$, 
for any $u\in C^{\infty}_c(M;\Sigma^2_0M)$ we have
\begin{align*}
(u,(\DeltaL+2n)u)
&= \|D^*u\|^2 + \|Du\|^2 + n \|u\|^2 - (u,\mathring{Rm}(u))\\
&\ge ( u, (n-\mathring{Rm})u)\\
&\ge (n-1)\kappa\|u\|^2.
\end{align*}
The same is true for $u\in H^{2,2}(M;\Sigma^2_0M)$
because $C^\infty_c(M;\Sigma^2_0M)$ is dense in that
space.  If $\kappa>0$, it follows immediately 
that $\DeltaL+2n$ has trivial $L^2$
kernel.  On the other hand, if $\kappa=0$, 
the sequence of inequalities above implies that
if $u\in H^{2,2}(M;\Sigma^2_0M)$ is a solution 
to $(\DeltaL+2n)u=0$,  
the nonnegative quantity $\< u, (n-\mathring{Rm})u\>_h$
must vanish identically on $M$.  Since the sectional
curvatures of $h$ approach $-1$ at infinity,
there is some compact set $K\subset M$ such that
$K_{\max}(p)\le -1/2$ for $p\in M\setminus K$,
and then \eqref{eq:Koiso-est}
implies that $u\equiv 0$ on $M\setminus K$. 
Since $\DeltaL$ is a Laplace operator,
it satisfies the weak unique continuation property
\cite{Aronszajn,NUW}, and therefore $u$ is identically zero.

Finally, suppose that the conformal infinity $[\hat h]$
has nonnegative Yamabe invariant.  Then the result of 
\cite{Lee} shows that the 
Laplacian satisfies the following $L^2$ estimate for 
smooth, compactly supported scalar
functions $u$:
\begin{displaymath}
(u,\Del^*\Del u) \ge \frac{n^2}{4}\|u\|^2.
\end{displaymath}
(The proof in \cite{Lee} required $h$ to have a $C^{3,\alpha}$
conformal compactification.  However, using 
Lemma 3.3.1 of \cite{AC}, it is easy to reduce that to
$C^{2,\alpha}$.
See also \cite{Wang} for a different proof.)
By Lemma \ref{lemma:crude-tensor-est}, therefore, the same is true
when $u$ is a smooth, compactly supported tensor field,
and by continuity for all $u\in H^{2,2}(M;\Sigma^2_0M)$.

Suppose $h$ has sectional curvatures bounded
above by 
$(n^2-8n)/(8n-8)$.
Then \eqref{eq:Koiso-est} gives
\begin{displaymath}
2\<u,\mathring{Rm}(u) \>_h 
\le
 2\left( n + \frac{(n-1)(n^2-8n)}{8n-8}\right)|u|_h^2 
= 
\frac{n^2}{4}|u|_h^2.
\end{displaymath}
If $u\in L^2(M;\Sigma^2_0M)$ is a solution to $(\DeltaL+2n)u=0$,
therefore,
\begin{align*}
0 
&=(u,(\DeltaL+2n)u)\\
&=(u,(\Del^*\Del - 2 \mathring{Rm})u)\\
&= \|\Del u\|^2 - \frac{n^2}{4}\|u\|^2
+ \left(u,({n^2}/{4}- 2 \mathring{Rm})u\right)\\
&\ge \left(u,({n^2}/{4}- 2 \mathring{Rm})u\right)\\
&\ge 0.
\end{align*}
It follows as before that the nonnegative function
$\<u,(n^2/4 -2\mathring{Rm})u\>_h$ must be identically
zero, and since the operator $(n^2/4 -2\mathring{Rm})$
is positive definite outside a compact set,
$u$ must be identically zero by analytic continuation.
\end{proof}

\backmatter

\bibliographystyle{amsplain}
\bibliography{fred}

\end{document}